\documentclass[11pt, oneside]{article}   	% use "amsart" instead of "article" for AMSLaTeX format
\usepackage{geometry}                		% See geometry.pdf to learn the layout options. There are lots.
\usepackage{graphicx}				% Use pdf, png, jpg, or eps§ with pdflatex; use eps in DVI mode

\graphicspath{{graphs_1/}}
\usepackage{subfigure}
\usepackage{dsfont}					% TeX will automatically convert eps --> pdf in pdflatex		
\usepackage{amssymb}
\usepackage{amsmath}
\usepackage{amsthm}
\usepackage{amsbsy,mathtools, bm}
\usepackage{xcolor}
\usepackage[normalem]{ulem}
\usepackage{hyperref}

\usepackage{geometry}                		% See geometry.pdf to learn the layout options. There are lots.
\usepackage{graphicx}				% Use pdf, png, jpg, or eps§ with pdflatex; use eps in DVI mode

\graphicspath{{graphs_1/}}
\usepackage{subfigure}
\usepackage{dsfont}					% TeX will automatically convert eps --> pdf in pdflatex		
\usepackage{amssymb}
\usepackage{amsmath}
\usepackage{amsthm}
\usepackage{amsbsy,mathtools, bm}
\usepackage{xcolor}
\usepackage[normalem]{ulem}
\usepackage{hyperref}

 \usepackage[disable,textsize=tiny]{todonotes}

\usepackage[normalem]{ulem}
\usepackage{enumitem}

\newtheorem{theorem}{Theorem}[section]

\newtheorem{lemma}[theorem]{Lemma}
\newtheorem{proposition}{Proposition}[section]
\newtheorem{assumption}{Assumption}[section]
\newtheorem{example}{Example}[section]
\newtheorem{definition}{Definition}[section]
\newtheorem{remark}{Remark}[section]

\newcommand{\dd}{\mathrm{d}}
\newcommand{\Span}{\mathrm{span}}

\newcommand{\W}{\mathds{W}}
\newcommand{\PSD}{\mathbb{S}_+}
\newcommand{\PD}{\mathbb{S}_{++}}
\newcommand{\diag}{\mathrm{diag}}
\newcommand{\coef}{\mathrm{coef}}

\newcommand{\eps}{\varepsilon}
\newcommand{\R}{\mathbb{R}}
\newcommand{\Obj}{\mathrm{Obj}}
\newcommand{\review}[1]{{\color{black}{#1}}}

\usepackage[refpage]{nomencl}
\makeatletter

\def\@@@nomenclature[#1]#2#3#4{%
 \def\@tempa{#2}\def\@tempb{#3}\def\@tempc{#4}%
 \protected@write\@nomenclaturefile{}%
  {\string\nomenclatureentry{#1\nom@verb\@tempa @[{\nom@verb\@tempa}]%
      \begingroup\nom@verb\@tempb\protect\nomeqref{\theequation}%
        |nomlabelref}{\@tempc}}%
 \endgroup
 \@esphack}
\makeatother

\makenomenclature 
 \usepackage[disable,textsize=tiny]{todonotes}

\usepackage[normalem]{ulem}
\usepackage{enumitem}
\usepackage{dsfont}

\title{Distributionally Robust Gaussian Process Regression and Bayesian Inverse Problems}

\author{Xuhui Zhang$^1$, Jose Blanchet$^1$, Youssef Marzouk$^2$, Viet Anh Nguyen$^3$, Sven Wang$^4$}
\definecolor{ssw}{rgb}{0.1,0.45,0.1}
\newcommand{\ssw}[2][]{\todo[color=ssw!40,#1]{#2}}

\definecolor{jb}{rgb}{0,0,1}
\definecolor{ymm}{rgb}{1,0.53,0.0}
\newcommand{\ymmtd}[2][]{\todo[color=ymm!40,#1]{(ymm) #2}}

\date{%
    $^1$Stanford University\\%
    $^2$Massachusetts Institute of Technology\\%
    $^3$Chinese University of Hong Kong\\
    $^4$Humboldt University Berlin\\[2ex]%
}						% Activate to display a given date or no date

\begin{document}
\maketitle
%\section{}
%\subsection{}

\begin{abstract}
We study a distributionally robust optimization formulation (i.e., a min-max game) for two representative problems in Bayesian nonparametric estimation: Gaussian process regression and, more generally, linear inverse problems. Our formulation seeks the best mean-squared error predictor in an infinite-dimensional space against an adversary who chooses the worst-case model in a Wasserstein ball around a nominal infinite-dimensional Bayesian model. The transport cost is chosen to control features such as the degree of roughness of the sample paths that the adversary is allowed to inject. We show that strong duality holds in the sense that max-min equals min-max, and that there exists a unique Nash equilibrium that can be computed by a sequence of finite-dimensional approximations. Crucially, the worst-case distribution is itself Gaussian. We explore the properties of the Nash equilibrium and the effects of hyperparameters through numerical experiments, demonstrating the versatility of our modeling framework.
\end{abstract}

%%%%%%%%%%%%%%%%%%%%%%%%%%%%%%%%%%%%%%%%%%%%%%
%% Please use \tableofcontents for articles %%
%% with 50 pages and more                   %%
%%%%%%%%%%%%%%%%%%%%%%%%%%%%%%%%%%%%%%%%%%%%%%
%\tableofcontents

%%%%%%%%%%%%%%%%%%%%%%%%%%%%%%%%%%%%%%%%%%%%%%
%%%% Main text entry area:
\section{Introduction}

Bayesian nonparametric estimation~\cite{ghosh2003bayesian}
is used ubiquitously in science, engineering, and other statistical application areas for both `direct' nonparametric regression and solving inverse problems. The computation of the posterior (or conditional) mean estimators, which are the most commonly used Bayesian point estimators, involves the solution of an infinite-dimensional optimization problem whose specification requires knowledge of the distributions at hand. In general, this problem has no closed-form solution. If the observations and the parameter of interest are jointly Gaussian, however, then the problem immediately becomes tractable~\cite{ref:gine2015mathematical,ref:stuart2010inverse}. The conditional expectation (the best mean-square estimator of the parameter) is an affine function of the observations. Conditional covariances also can be evaluated in closed forms, enabling quantification of prediction uncertainty.

Gaussianity is a strong assumption that is easily violated in reality. More generally, the joint probabilistic model for the parameter of interest and the observations is often---perhaps inevitably---\emph{misspecified}. This misspecification can take myriad forms, including incorrect prior assumptions on the smoothness or dependence structure of the unknown parameter and incorrect assumptions on the nature of the data-generating process---which in turn involves assumptions on both the observational noise and, in the inverse problem setting, on the ``forward'' operator relating the parameters to the observations. In these situations, it is desirable to ensure some form of \emph{robustness}, e.g., to construct a nonparametric estimator that hedges against the impact of model misspecification on mean-square error. One also would like to represent possible modeling errors nonparametrically and in a way that fully reflects the infinite-dimensional nature of the regression and inverse problem settings.

This paper addresses model misspecification in nonparametric settings by adapting and extending ideas from distributionally robust optimization (DRO)~\cite{ref:delage2010distributionally, ref:rahimian2019distributionally}. \review{DRO formulations create a min-max game between a decision-maker and a fictitious adversary: the decision-maker finds a decision or an estimator that minimizes the worst-case expected loss over all possible distributions the adversary can choose. The set containing the adversary's feasible choice is the \textit{ambiguity set}. In Wasserstein DRO, this ambiguity set is prescribed as a Wasserstein neighborhood around a nominal distribution~\cite{ref:kuhn2019wasserstein, ref:blanchet2021statistical, ref:blanchet2024distributionally}. By choosing appropriately the nominal distribution and the Wasserstein distance, the Wasserstein DRO solutions can exhibit many desirable properties, e.g. in the setting of (group) LASSO \cite{ref:blanchet2016robust}, norm regularization~\cite{ref:shafieezadeh2019regularization}, shrinkage~\cite{ref:nguyen2018distributionally}, ridge regression~\cite{ref:li2022tikhonov}, etc.

The first line of work on Wasserstein DRO set the nominal distribution as the empirical distribution of the available data~\cite{ref:esfahani2018data,ref:gao2017wasserstein}. }\review{Subsequently, finite-dimensional parametric information under the form of Gaussianity is injected into the nominal distribution for the estimation task~\cite{ref:nguyen2018distributionally, ref:shafieezadeh2018wasserstein}. To capture more structural information about the problem, independence has also added to the ambiguity set to model the independence between the finite-dimensional signal and the finite-dimensional error in the minimum mean square error estimation problem~\cite{ref:nguyen2019bridging}.} \review{In Wasserstein DRO when the nominal distribution is a generic stochastic process, the papers~\cite{ref:gao2016distributionally,ref:blanchet2019quantifying} derived some general duality results quantifying the worst-case expected loss over the Wasserstein ambiguity set (i.e., the inner maximum of a DRO formulation), while the full min-max game was not studied in this setting.}

This paper significantly extends the existing literature by introducing a DRO formulation for \textit{in}finite-dimensional Gaussian models. In our formulation, the decision-maker seeks an estimator of the unknown parameter that minimizes a certain Bayes risk, and the adversary chooses a probabilistic model that departs from the baseline/nominal Gaussian model assumed by the decision maker. Specifically, we will allow the adversary to select a model in a nonparametric way within a $\delta$-Wasserstein ball around the nominal model, where $\delta$ is the ambiguity size. The particular Wasserstein geometry that we impose allows the adversary to select models that depart significantly not only \textit{from the Gaussian assumption} but also \textit{from the smoothness properties} dictated by the decision-makers' choice of prior. Consequently, our min-max formulation allows us to explore and assess the impact of model misspecification effectively.

We briefly contrast the nominal and robust \review{inverse} problems to make these ideas concrete; a precise and technical presentation is deferred to Section~\ref{sec:notations}. Let $b^0$ represent the parameter of interest, modeled as a real-valued random process with continuous sample paths on some compact domain $\mathcal{D} \subseteq \mathbb{R}^d$. Suppose we have a finite number of real-valued observations $(Y_1, \ldots, Y_m)$, specified as 
\[Y_i = T(b^0)(x_i) + \epsilon_i^0,\] where $T$ is a bounded linear operator, $(x_i)_{i=1}^m \in \mathcal{D}$ are collocation or design points, and $(\epsilon_i^0)_{i=1}^m$ are real-valued random variables with independent entries representing observational noise. In the nominal case, $b^0$ and $\epsilon_i^0$ are endowed with a Gaussian prior measure $P_0$ on $L^2(\mathcal D)\times \R^m$, and under $P_0$, $\epsilon_i^0$ are independent of $b^0$. We seek an estimator $\phi: \mathbb{R}^m \to L^2(\mathcal{D})$ that minimizes the Bayes risk
\begin{equation}
\min_{\phi}~\mathbb{E}_{P_0}\left[ \left \Vert b^0-\phi(Y_1,\ldots,Y_m) \right \Vert ^2_{L^2(\mathcal{D})}\right].
\label{eq:nominalintro}
\end{equation} 
It is well known that this estimator is given by the conditional expectation $\phi(Y_1,\ldots,Y_m)  = \mathbb{E} [ b^0 \, \vert \, Y_1, \ldots, Y_m ]$ of the posterior distribution. We build the Wasserstein distributionally robust counterpart of this problem by introducing an ambiguity set of probability measures, $\{ P  :  \W(P,P_0) \leq \delta \}$, where $\W$ is a Wasserstein distance on the space of (Borel) probability measures on $C(\mathcal{D}) \times \mathbb{R}^m$.  We then seek an estimator $\phi$ that minimizes the worst-case Bayes risk over this ambiguity set:
\begin{equation}
\inf \limits_{\phi} \sup \limits_{P \, : \, \W(P,P_0)\leq\delta} \mathbb{E}_P\left[ \left \Vert b - \phi(Y_1,\ldots,Y_m) \right \Vert^2_{L^2(\mathcal{D})}\right ],
\label{eq:defminmaxintro}
\end{equation}
where now $Y_i = T(b)(x_i) + \epsilon_i$, and  $b$ and $\epsilon_1, \ldots, \epsilon_m$ are jointly distributed according to $P$. The robust formulation~\eqref{eq:defminmaxintro} thus adds a fictitious adversary to the nominal formulation~\eqref{eq:nominalintro}, and the adversary is the inner player solving the supremum subproblem.

% \viet{Here we start to mention roughness, but it has not been well defined. Is it possible to add one sentence to define roughness of a sample path of $b$? ``Given two sample path of $b$, denoted momentarily as $\hat b_1$ and $\hat b_2$, we say that $\hat b_1$ is rougher than $\hat b_2$ if [xxxx]}{\color{red}xz:added a sentence below, we can talk about the ``spaces'' that rough sample path live in} 

By choosing among distributions in the ambiguity set, the adversary may inject additional roughness to sample paths of $b$, beyond what is specified by the nominal distribution $P_0$. The ability to add roughness is directly tied to our specification of $\W$: we consider perturbations to $b$ that are elements of an RKHS $\mathcal{H}_w$ whose norm is parameterized by a sequence of weights $(w_n)_{n \geq 1}$. Using these weights, we obtain a flexible class of ground cost defining our Wasserstein distance which can adjust the penalty for transportation along different ``modes'' of the sequence decomposition of $b$ in the Karhunen-Loeve expansion. \review{In particular, the choice of weights $(w_n)_{n \geq 1}$ induces a hierarchy of the Hilbert space norms so that rougher perturbations correspond to spaces whose weights $(w_n)_{n \geq 1}$ grow to infinity in a slower rate.} This infinite-dimensional construction is an important novelty from the DRO perspective. \review{To guide intuition, consider the Mat\'ern process (see Example~\ref{ex:matern} in the main body). The smoothness of the Mat\'ern kernel is controlled by a parameter $\alpha$, which determines the regularity of the sample paths. When $\alpha$ is large, the paths are smooth, while a smaller $\alpha$ results in rougher paths. In our framework, optimal transport perturbation of $b$ with a slower-growing weights $(w_n)_{n \geq 1}$ is analogous to reducing the value of $\alpha$, which introduce more roughness into the sample paths of $b$.}
The adversary also can modify the distribution of the additive noise $\epsilon_i$, which can be understood in part as compensating for the misspecification of the nominal forward operator $T$ \cite{kennedy2001bayesian}. Moreover, the adversary can replace $P_0$ (in any of its marginals or jointly) with a non-Gaussian distribution.

Our assumptions on the operator $T$ in the formulation above will encompass many linear inverse problems \cite{ref:stuart2010inverse,ref:Dashti2017inverse}, e.g., learning the initial or boundary conditions of a heat equation, or canonical problems in computerized tomography \cite{ref:natterer2001mathematics}. By setting $T = \text{Id}$, however, we recover the important case of Gaussian process regression, to which our main results immediately apply. For non-identity $T$, our results apply both to recovering $b$ (e.g., solving the inverse problem) and to estimating $u$ (e.g., PDE-constrained regression).

Wasserstein-type distances defined on the space of stochastic processes were recently studied by~\cite{ref:backhoff2020adapted,ref:acciaio2017causal,ref:nadal2019wasserstein}. The focus therein is on processes indexed by a one-dimensional parameter (representing, for instance, time). A typical setting involves price processes in finance~\cite{ref:backhoff2020adapted, ref:acciaio2017causal}, where one needs to define a Wasserstein distance that respects causal structure (i.e., filtrations). In contrast, we consider in this work a Wasserstein distance on multi-dimensional fields.

Note also that in our formulation \eqref{eq:defminmaxintro}, we have an infinite-dimensional action space for the outer player and an infinite-dimensional action space for the inner player. Moreover, the actions of the inner player are themselves probability measures on infinite-dimensional spaces. These features differentiate our analysis from prior work, such as~\cite{ref:blanchet2016robust,ref:esfahani2018data,ref:shafieezadeh2018wasserstein}.
To our knowledge, previous work in the DRO literature assumes either the action set of the decision maker to be finite-dimensional or the probability measures (the action set of the fictitious adversary) to be supported on finite-dimensional spaces. 

Alternative to the Wasserstein distance, one can also naturally formulate problem~\eqref{eq:defminmaxintro} with an information divergence. Distributionally robust conditional mean estimation with a relative entropy ambiguity set centered on a 
multivariate Gaussian prior has been studied in \cite{ref:levy2004robust,ref:levy2013robust,ref:yi2022robust}. As shown in~\cite[Theorem~1]{ref:levy2004robust} and \cite[Theorem~1]{ref:levy2013robust}, however, the resulting robust estimator coincides exactly with its non-robust counterpart, and only the posterior covariance is inflated. The same conclusion holds for ambiguity sets constructed from the $\tau$-divergence family~\cite{ref:zorzi2017robust,ref:zorzi2017ontherobust}. 
In strong contrast to these results, the robust estimator of our formulation~\eqref{eq:defminmaxintro} typically differs from its non-robust counterpart.

While the results just discussed were derived in the finite-dimensional setting, it is reasonable to expect that similar properties would be preserved under reasonable assumptions in the infinite-dimensional setting, which is our concern here. Distributionally robust formulations of non-causal filtering that employ a relative entropy or $\tau$-divergence ambiguity set centered on a stationary Gaussian process prior have been studied in~\cite{ref:levy2013robust,ref:zorzi2017ontherobust}, where it was shown that the nominal non-causal Wiener filter remains optimal. Moreover, an infinite-dimensional DRO formulation based on, e.g., relative entropy or any other criterion that requires the existence of a likelihood ratio will typically restrict the adversary to preserve sample path properties, such as the degree of sample path smoothness under the baseline model (i.e., the prior $P_0$). Choosing instead a Wasserstein ambiguity set, as we will explain, permits adversarial distributions with rougher sample paths than the prior, so that \textit{smoothness misspecification} is naturally addressed by our robust estimation problem.

%% SUMMARY OF CONTRIBUTIONS
We now summarize our main contributions. Under reasonable assumptions to be made precise later:
\begin{itemize}[leftmargin=5mm]
\item  We analyze \review{the Wasserstein distributionally robust inverse} problem~\eqref{eq:defminmaxintro} and show that strong duality holds, in the sense that the \review{infimum} and the \review{supremum} operators can be switched without any loss of optimality.
    
\item We show that there exists an upper bound $\delta_0>0$ such that if $0<\delta<\delta_0$, problem \eqref{eq:defminmaxintro} also admits a unique Nash equilibrium pair $(\phi_\infty^\star, P_\infty^\star)$. Moreover, the worst-case distribution $P_\infty^\star$ involves a modified Gaussian process with potentially rougher paths than the prior. Consequently, the robustified conditional mean remains affine in the observations, and the optimization problem is tractable.
    
\item We numerically approximate problem~\eqref{eq:defminmaxintro} by a sequence of finite-dimensional counterparts and, therefore, obtain a procedure to compute the associated Nash equilibrium. Our numerical algorithm adapts the finite-dimensional Frank-Wolfe algorithm in~\cite{ref:shafieezadeh2018wasserstein}. 
    
\end{itemize}

%% INTERPRETATION OF WORST-CASE (keep)
One way to interpret our results is that the \review{solution of linear inverse problems and Gaussian process regression} can be made robust in a nonparametric sense. The interpretation of the worst-case covariance function (i.e., the covariance of the Gaussian process $P_\infty^\star$) enables the computation of the upper bound on the worst-case mean square error and thus an upper bound on the quality of the robust solution. In our numerical examples, we explore the structural properties of this worst-case covariance function and the associated robust estimator.
In particular, we explore the prior and posterior covariances of $b$ under $P_0$ and $P_\infty^\star$ (i.e., $\mathbb{C}\text{ov}_{P}[b]$ and $\mathbb{C}\text{ov}_{P}[b \, \vert \, y_1, \ldots, y_m]$ for $P=P_0$ and $P=P_\infty^\star$) and find that the worst-case distributions (both prior and posterior) have greater uncertainty in regions of $\mathcal{D}$ where information is limited, which intuitively guarantees greater robustness of the predictions. Moreover, we observe that in cases where there is (i) a smoother nominal prior, (ii) a smaller transport penalty in basis directions that induce roughness, (iii) a larger $\delta$, or (iv) a smaller nominal observational noise, the worst-case distributions induce sharper contrasts between the observed and unobserved locations along both the prior and posterior sample paths. One criticism of our results is that the distributionally robust solution may be pessimistic. While our framework provides the flexibility of tuning hyperparameters such as the weights in defining the Wasserstein distance $\W$, we leave the precise specification of these hyperparameters for future work.
    
%% PAPER PLAN
The remainder of this paper is organized as follows. We introduce our problem setup in Section~\ref{sec:notations}. We present our main theoretical results in Section~\ref{sec:mainresults}, illustrate the applicability of our general framework by highlighting several examples in Section~\ref{sec:examples}, and present simple numerical experiments in Section~\ref{sec:numerics}. All proofs are deferred to Appendix A. \review{For convenience, we provide a glossary of notations in Appendix B (list of important notations, their brief description, and page of reference).}
%xz: note that the A B is hardcoded, be careful

\section{Problem Statement}\label{sec:notations}

Let $\mathcal{D}\subset\mathbb{R}^d$\label{nomencl:D}\nomenclature{$\mathcal{D}$}{compact set}{nomencl:D}, $d\geq1$ be a compact set. We write $C(\mathcal{D})$ to denote the space of real-valued continuous functions on $D$, which is naturally endowed with the sup-norm $\|\cdot\|_{C(\mathcal{D})}$. We denote by $L^2(\mathcal{D})$ the space of real-valued square-integrable functions on $\mathcal{D}$. Since $\mathcal{D}$ is compact, we have $C(\mathcal{D}) \subseteq L^2(\mathcal{D})$. 

% \xuhui{simple general formula involving $u$}\ymmtd{is this blue text an old/defunct comment?}

We introduce a probability measure $P_0$\label{nomencl:P0}\nomenclature{$P_0$}{nominal probability measure}{nomencl:P0} under which the so-called prior input process $b^0$\label{nomencl:b0}\nomenclature{$b^0$}{prior input process}{nomencl:b0} is a $C(\mathcal{D})$-valued centered Gaussian random field.\footnote{For the theory of Banach space-valued Gaussian random variables, we refer to~\cite[Chapter~2]{ref:gine2015mathematical} and~\cite{ref:vaart2008reproducing}.}  Further, under the inclusion $C(\mathcal{D})\hookrightarrow L^2(\mathcal{D})$ where $\hookrightarrow$ denotes the inclusion map, the random field $b^0$ can be viewed as $L^2(\mathcal{D})$-valued. This random field generates a positive definite kernel, namely, $\mathcal{K}(x,x')=\mathbb{E}_{P_0}[b^0(x)b^0(x')]$ and thus an associated reproducing kernel Hilbert space (RKHS) which is obtained as the closure of functions of the form $f(x)=\mathbb{E}_{P_0}[(a_1b^0(x_1)+ \cdots +a_nb^0(x_n))b^0(x)]$. The closure can be taken relative to the norms $\|\cdot\|_{C(\mathcal{D})}$ and $\|\cdot\|_{L^2(\mathcal{D})}$; both limiting procedures coincide~\cite[Lemma~8.1]{ref:vaart2008reproducing}. By the spectral decomposition of the covariance operator of $b^0$ (i.e., $\mathcal{K}(\cdot,\cdot)$, \cite[Example~2.6.15]{ref:gine2015mathematical}), there exists a complete orthonormal system $\{e_n\}_{n=1}^\infty$ of $L^2(\mathcal{D})$\label{nomencl:en}\nomenclature{$\{e_n\}_{n=1}^\infty$}{orthonormal system of $L^2(\mathcal{D})$}{nomencl:en} where $e_n\in C(\mathcal{D})$), an i.i.d.~sequence of standard univariate normal random variables $\{g_n\}_{n=1}^\infty$, and a non-negative sequence of ``eigenvalues'' $\{\kappa_n^2\}_{n=1}^\infty$ satisfying $\sum_{n=1}^\infty \kappa_n^2<\infty$, such that under $P_0$,
\[
b^0 = \sum_{n\geq1}\kappa_ng_ne_n,
\]
\label{nomencl:kappan}\nomenclature{$\kappa_n$}{eigen-values in Karhunen Loeve expansion}{nomencl:kappan}where the convergence of the above infinite sum occurs in $C(\mathcal{D})$, and thus also in $L^2(\mathcal{D})$, almost surely. We impose a full-rank assumption on the prior $b^0$ in the following sense, which is necessary to ensure that the worst-case distribution is unique in the proof of Theorems~\ref{thm:swap} and~\ref{thm:existnash} below.
\begin{assumption}[Full rank]\label{assmp:fullrank}
We assume that $\kappa_n\neq0$ for all $n\geq1$.
\end{assumption}
% \viet{what do you mean by the ``RKHS of $b^0$?}\xuhui{xz:In the sense of~\cite[Chapter~2]{ref:gine2015mathematical} or~\cite{ref:vaart2008reproducing}?}
% Assumption~\ref{assmp:fullrank} implies that the support of $b^0$ is not contained in a proper subspace of $C(\mathcal{D})$. This assumption is necessary to ensure that the worst-case distribution is unique in the proof of Theorems~\ref{thm:swap} and~\ref{thm:existnash} below. 

A general class of Gaussian smoothness priors can be constructed via the Laplace operator on $\mathcal{D}$.  Specifically, we will use the following class of Mat\'ern processes in Sections~\ref{sec:examples} and~\ref{sec:numerics}, which provides natural prior distributions for $\alpha$-regular functions vanishing at the boundary $\partial \mathcal D$\footnote{\review{$\alpha$-regular functions refer to functions that are H\"{o}lder (and Sobolev)-continuous of order $\alpha-d/2-q$ for any $q>0$}.}.
\begin{example}[Mat\'ern prior~{\cite[Equation~(2)]{ref:lindgren2011explicit}}]\label{ex:matern}
Suppose $\mathcal{D}$ has a smooth boundary. The prior with a Mat\'ern covariance function with parameters $\kappa\geq0$ and $\alpha>\frac{d}{2}$\label{nomencl:matalpha}\nomenclature{$\alpha$}{smoothness parameter in Example~\ref{ex:matern}}{nomencl:matalpha} controlling the smoothness can be expressed as
\[
b^0 = \sum_{n\geq1}\left(\kappa^2 + \lambda_n\right)^{-\frac{\alpha}{2}} g_n e_n,
\]
where the eigenvalues $\lambda_n$\label{nomencl:lambdan}\nomenclature{$\lambda_n$}{Dirichlet-Laplacian eigenvalues in Example~\ref{ex:matern}}{nomencl:lambdan} and 
eigenfunctions $e_n$ correspond to the Dirichlet-Laplace operator on $\mathcal{D}$~\cite[Corollary~5.1.5]{ref:taylor1996partialone}. The eigenvalues $\lambda_n$ satisfy Weyl's law $\lambda_n = \Theta(n^{2/d})$~\cite[Corollary~8.3.5]{ref:taylor1996partial} and the eigenfunctions $e_n\in  C^\infty(\bar{\mathcal{D}})$, $n\geq1$, where $\bar{\mathcal{D}}$ denotes the closure of $\mathcal{D}$ and $C^\infty(\bar{\mathcal{D}})$ is the space of infinitely smooth functions on $\bar{\mathcal{D}}$.
\end{example}

We consider perturbations to the nominal prior $b^0$ by borrowing ideas from the field of distributionally robust optimization (DRO). We assume that the perturbations are supported in a space of continuous functions $\mathcal{H}_w$\label{nomencl:hw}\nomenclature{$\mathcal{H}_w$}{RKHS for perturbations}{nomencl:hw} that is also an RKHS; in particular $\mathcal{H}_w$ is also a Polish space, which is important to invoke key duality results to study the maximization in our DRO formulation. The useful feature of RKHS is that point evaluation functionals are well-defined and continuous with respect to the Hilbert space norm~\cite[
Example 2.2]{{ref:vaart2008reproducing}}. 
% [\sven{Here or in the introduction, motivate why this assumption is needed. Is it essentially the Sobolev embedding $H^\alpha\subseteq L^\infty$?}] 
Specifically, we define the space 
\begin{equation}\label{eq:rkhsdefnotation}
\mathcal{H}_w = \left\{f\in L^2(\mathcal{D}):\sum_{n\geq1}w_n \langle f,e_n\rangle^2 <\infty\right\},
\end{equation}
\label{nomencl:wwn}\nomenclature{$w_n$}{Sequence of weights parametrizing the Hilbert space norm}{nomencl:wwn}which is parameterized by a positive sequence $w=(w_n)_{n\geq1}$. Notice that from this point, we abbreviate the inner product $\langle \cdot,\cdot\rangle_{L^2(\mathcal{D})}$ as $\langle\cdot,\cdot\rangle$. Typically, the Hilbert norm on $\mathcal{H}_w$ is stronger than the usual $L^2(\mathcal{D})$ norm. More precisely, we impose the following assumption on $w$ and the basis $\{e_n\}_{n=1}^\infty$:

\begin{assumption}[RKHS conditions]\label{assmp:perturbation}
 Assume that $\lim_{n\to\infty}w_n=\infty$ and $\mathcal{H}_w$ is endowed with the inner product 
\[
\langle f,\tilde f\rangle_{\mathcal{H}_w} =\sum_{n\geq1}w_n\langle f,e_n\rangle\langle \tilde f ,e_n\rangle \qquad \forall f,\tilde f \in \mathcal{H}_w.
\]
Further, the space $\mathcal{H}_w$ equipped with the Hilbert space norm $\|\cdot\|_{\mathcal{H}_w}$ is an RKHS compactly embedded in $C(\mathcal{D})$.
% , identifying $f\in\mathcal{H}_w$ with its continuous version in $L^2(\mathcal{D})$.
\end{assumption}
Under Assumption~\ref{assmp:perturbation},  $\mathcal{H}_w$ is a separable Hilbert space such that point evaluation functionals are well-defined and continuous. Note that the sequence $w$ controls the roughness (or equivalently the smoothness) of functions in $\mathcal{H}_w$. Throughout the rest of the paper, we will fix a given sequence $w$ satisfying Assumption~\ref{assmp:perturbation}. The norm $\|\cdot\|_{\mathcal{H}_w}$ will define our adversarial perturbations. The next example, a continuation of Example~\ref{ex:matern}, provides intuition about the interpretation of $w$ in terms of roughness.

\begin{example}[RKHS space]\label{ex:perturbation}
Let the eigenvalues $\lambda_n$ and 
eigenfunctions $e_n$ correspond to the Dirichlet-Laplacian operator on $\mathcal{D}$. Consider $w_n = \Theta(\lambda_n^\beta)$.\label{nomencl:matbeta}\nomenclature{$\beta$}{roughness parameter in Example~\ref{ex:perturbation}}{nomencl:matbeta}\footnote{For positive sequences $\{a_n\},\{b_n\}$, the notation $a_n=\Theta(b_n)$ means that $1/c_0 \leq a_n/b_n \leq c_0$, for some $c_0\in(0,\infty)$.} Then $\beta$ controls the roughness of functions in $\mathcal{H}_w$. Note that the ``spectrally defined'' spaces $\mathcal H_w$ are subspaces of the classical Sobolev spaces on $\mathcal D$.
% [\sw{Comment on the `spectrally defined' spaces $\mathcal H_w$ being a subspace the classical Sobolev spaces on $\mathcal D$ such that we may actually apply the Sobolev embedding.}]
For any $\beta>\frac{d}{2}$, by the Sobolev embedding theorem~\cite[Proposition 4.1.3]{ref:taylor1996partialone}, we can identify $f\in\mathcal{H}_w$ with its continuous version. Under this identification $\mathcal{H}_w$ is a RKHS in $C(\mathcal{D})$.
\end{example}

Now, suppose that we have a linear forward operator $T$\label{nomencl:T}\nomenclature{$T$}{linear forward operator}{nomencl:T} that maps sample paths of the prior input process $b^0$ to sample paths of another process $u$\label{nomencl:u}\nomenclature{$u$}{observation process $T(b)$}{nomencl:u}, which also takes values in $C(\mathcal{D})$. We make the following assumptions on this operator.

\begin{assumption}[Operator]\label{assmp:operator} We assume the following:
\begin{enumerate}[label=(\roman*)]
    \item The forward map $T:C(\mathcal{D})\to C(\mathcal{D})$ is linear and bounded (with operator norm $C_T>0$).
%    \item There exists a positive sequence $\tilde w= \{\tilde w_n\}_{n=1}^\infty$, $\lim_{n\to\infty}\tilde w_n\to\infty$, $\tilde w_n = o(w_n)$, and the space $\mathcal{H}_{\tilde w}$ defined in the same manner as~\eqref{eq:rkhsdefnotation} with its Hilbert space norm $\|\cdot\|_{\mathcal{H}_{\tilde w}}$ constitutes a RKHS in $C(\mathcal{D})$.

    \item There exists a positive sequence $\tilde w= \{\tilde w_n\}_{n=1}^\infty$ and a corresponding space $\mathcal{H}_{\tilde w}$ as in \eqref{eq:rkhsdefnotation}, with Hilbert space norm $\|\cdot\|_{\mathcal{H}_{\tilde w}}$, which constitutes a RKHS continuously embedded in \ssw{Continuous embedding (or variation of this phrase)} $C(\mathcal{D})$ such that for some positive constant $C_{\tilde w}$,
 \[
 \|T(f)\|_{\mathcal{H}_{\tilde w}}\leq C_{\tilde w}\|f\|_{\mathcal{H}_{\tilde w}}\qquad\forall f\in\mathcal{H}_{\tilde w}.
 \]
 In other words, $T$ is bounded when restricted to $\mathcal{H}_{\tilde w}$.
\end{enumerate}
\end{assumption}

The forward operator $T$ defines our data-generating process. In particular, let $x_i\in\mathcal{D}$, $i=1,\ldots,m$\label{nomencl:xi}\nomenclature{$x_i$}{design points}{nomencl:xi} be the design points.\label{nomencl:numm}\nomenclature{$m$}{number of design points}{nomencl:numm} Also under $P_0$, let $\epsilon^0= (\epsilon_1^0,\ldots,\epsilon_m^0)$ be a vector of errors\label{nomencl:eps0}\nomenclature{$\epsilon^0$}{nominal vector of errors}{nomencl:eps0}. Then we observe a single, noisy path of $u$ at the design points, i.e.,
\[
Y_i = u^0(x_i) + \epsilon^0_i,\quad  i =1,\ldots,m,
\]
\label{nomencl:yi}\nomenclature{$Y_i$}{observations at design point $x_i$}{nomencl:yi}where $u^0(x_i) = T(b^0)(x_i),\ i=1,\ldots,m$ are point evaluations of a single sample path. 
\review{We are now ready to state the core assumption about the nominal distribution $P_0$. 
\begin{assumption}[Independence]\label{assmp:independence_revision} We assume that under the nominal distribution $P_0$, $\epsilon^0$ is a vector of independent $\mathcal{N}(0,\sigma^2)$ errors,
and that $b^0$ and $\epsilon^0$ are independent. Consequently, under $P_0$, the pair $(b^0,\epsilon^0)$ constitutes a Gaussian random variable on the product space $C(\mathcal{D})\times\mathbb{R}^m$, which is Banach with the product norm.
\end{assumption}
}

% so that any bounded linear functional (chosen from the topological dual space of $C(\mathcal{D})\times\mathbb{R}^m$) of $(b^0,\epsilon^0)$ is univariate Gaussian~\cite[Example 2.1.16]{ref:gine2015mathematical}.

% \begin{assumption}\label{assmp:operator}
%  The forward map $T:C(\mathcal{D})\to C(\mathcal{D})$ is linear and bounded (with operator norm $C_T>0$, and therefore not trivial). Moreover, define 
%  \[
%  \mathcal{H}_s= \left\{f\in L^2(\mathcal{D}):\sum_{n\geq1}\lambda_n^s \langle f,e_n\rangle^2 <\infty\right\},
%  \]
%  with its Hilbert space norm $\|\cdot\|_{\mathcal{H}_s}$ in the same manner as~\eqref{eq:mainhilbnorm}, we assume there exists $\frac{d}{2}<s<\beta$ and a positive constant $C_s$, such that
%  \[
%  \|T(f)\|_{\mathcal{H}_s}\leq C_s\|f\|_{\mathcal{H}_s},\qquad\forall f\in\mathcal{H}_s.
%  \]
% \end{assumption}

We will see examples in Section~\ref{sec:examples} where Assumption~\ref{assmp:operator} on the forward operator is satisfied; these include Gaussian process regression (where $T =\text{Id}$) and several canonical linear inverse problems. Assumption~\ref{assmp:operator} entails certain compatibility of the forward operator with the prior basis $\{e_n\}_{n=1}^\infty$, which is similar to the assumption of ``norm equivalence on regularity scales'' in the literature on \textit{linear} Bayesian inverse problems~\cite{ref:gugushvili2020baeysian,ref:agapious2013posterior}. We note that other assumptions in the literature exist, e.g., the ``band-limited'' assumption in~\cite{ref:ray2013bayesian}. In our framework, we will consider Wasserstein perturbations to the prior $P_0$, which encompass a much richer collection of prior families than typically considered in the literature. 

\subsection{The Nominal \review{Inverse} Problem} 

% \viet{is it possible to rename the section title to ``the nominal inverse problem"?}
% We assume the ability to sample noisy observations of $u$ at $m$ design points in $\mathcal{D}$. Let  $\epsilon^0= (\epsilon_1^0,\ldots,\epsilon_m^0)$ be independent $\mathcal{N}(0,\sigma^2)$ errors, and denote the joint measure on $b^0$ and $\epsilon^0$ as $P_0$. Note that $b^0$ and $\epsilon^0$ are independent, then $(b^0,\epsilon^0)$ is a Gaussian random variable on the product space $C(\mathcal{D})\times\mathbb{R}^m$, so that any bounded linear functional (chosen from the topological dual space of $C(\mathcal{D})\times\mathbb{R}^m$) of $(b^0,\epsilon^0)$ is univariate Gaussian.
% \viet{stupid question, but is it possible to write down specifically what it means by ``a Gaussian on $C(\mathcal{D})\times\mathbb{R}^m$"? Does that require specifying a mean function and a covariance function?}. 
Let $\mathcal M$\label{nomencl:M}\nomenclature{$\mathcal M$}{space of measurable maps}{nomencl:M} denote the space of measurable maps from the data space $\mathbb{R}^m$ to the signal space $C(\mathcal{D})$ and let $\mathcal{P}$\label{nomencl:mcP}\nomenclature{$\mathcal P$}{space of (Borel) probability measures}{nomencl:mcP} denote the space of (Borel) probability measures on $C(\mathcal{D})\times\mathbb{R}^m$.  The classical Bayes risk minimization for $L^2(\mathcal{D})$ loss is
\begin{equation}
\min_{\phi\in\mathcal{M}}~\mathbb{E}_{P_0}\left[\|b^0-\phi(Y_1,\ldots,Y_m)\|^2_{L^2(\mathcal{D})}\right],
\label{eq:nominalregression}
\end{equation}
where $\phi(Y_1,\ldots,Y_m)$\label{nomencl:phi}\nomenclature{$\phi$}{non-parametric predictor for $b$}{nomencl:phi} is the estimator for $b^0$ given observations $(Y_1,\ldots,Y_m)$. We regard the problem as the nominal estimation problem since the Bayes risk is evaluated under the nominal measure $P_0$. The solution of the nominal problem is the posterior mean (or conditional expectation), i.e., 
\[
\phi_0(Y_1,\ldots,Y_m)(x) =\mathbb{E}_{P_0}\left[b^0(x)|Y_1,\ldots,Y_m\right].
\]
This follows by noting that the $L^2(\mathcal{D})$-norm is a Lebesgue integral and by using Fubini's theorem. Since we assume that the nominal distribution $P_0$ is Gaussian, this estimator corresponds to the linear prediction rule,
\[
\mathbb{E}_{P_0}\left[b^0(x)|Y_1,\ldots,Y_m\right]= \left(k^0(x,x_1),\ldots,k^0(x,x_m)\right)\cdot (K^0)^{-1} \cdot \left(Y_1,\ldots,Y_m\right)^\top,
\]
where $k^0(x,x_i) = \mathbb{E}_{P_0}[b^0(x)Y_i]$, $K^0_{ij} = \left(\mathbb{E}_{P_0}[Y_i \, Y_j]\right)$, and $K^0 \in \PD^{m}$, where $\PD^{m}$ denotes the set of strictly positive definite matrices. The invertibility of $K^0$ can be directly deduced from Assumption~\ref{assmp:independence_revision}.

% {\color{red}\sout{In the introduction, we wrote in \eqref{eq:nominalintro} an analogous estimation problem for $b^0$, which we shall revisit below.
% % in Section~\ref{sec:inverserobust}. 
% For non-identity $T$, \eqref{eq:nominalintro} corresponds to solving a linear inverse problem while \eqref{eq:nominalregression} is regression (e.g., imputing the rest of $u^0$ given noisy pointwise observations). \ymmtd{I added this so as not to forget about the inverse problems, particularly since we wrote it first in the intro. Please adjust as needed.{\color{red}xz:this makes sense to me!}} }  xz:no longer needed, moved to section 3.3}

% [\sven{Undefined notation.}]

% \viet{not sure if we need to mention this here-- it looks quite distracting}Alternatively, we also consider the inverse problem
% \[
% \min_{\phi_b\in\mathcal{M}}\mathbb{E}_{P_0}\left[\|b(\cdot)-\phi_b(\cdot;Y(x_1),\ldots,Y(x_m))\|^2_{L^2(\mathcal{D})}\right],
% \]
% whose linear prediction rule is given by
% \[
% \phi_b(x;Y(x_1),\ldots,Y(x_m)) = \left(k_b(x,x_1),\ldots,k_b(x,x_m)\right)\cdot (K)^{-1} \cdot \left(Y(x_1),\ldots,Y(x_m)\right)^\top,
% \]
% where $k_b(x,x_j) = \mathbb{E}_{P_0}[b(x)Y(x_j)]$ due to Gaussianity.

\subsection{The Distributionally Robust Optimization Inverse Problem}
\label{sec:droformulation}
% Let  $\epsilon^0= (\epsilon_1^0,\ldots,\epsilon_m^0)$ be independent $\mathcal{N}(0,\sigma^2)$ errors. Denote the joint measure on $b^0$ and $\epsilon^0$ as $P_0$. Note that $b^0$ and $\epsilon^0$ are independent, then $(b^0,\epsilon^0)$ is a Gaussian random variable on $C(\mathcal{D})\times\mathbb{R}^m$. Let $\mathcal M$ denote the space of measurable functions on $\mathcal{D}\times \mathbb{R}^m$ and $\mathcal{P}$ denote the space of probability measures on $C(\mathcal{D})\times\mathbb{R}^m$. For design points $x_i\in\mathcal{D},i=1,\ldots,m$, 

Instead of considering the nominal measure $P_0$ on the prior and noise distributions, we postulate a min-max game where an adversary chooses a measure $P$ in opposition to the decision-maker's choice of estimator.  In particular, for some forward map $T$ describing the relationship between the parameter $b$ and the function $u = T(b)$, we assume that the observations $Y_1,\dots, Y_m$ arise as
\begin{equation}\label{eq:system}
\begin{cases} 
(b,\epsilon_1,\ldots,\epsilon_m) \sim P,\\
Y_1=u(x_1)+\epsilon_1, \ldots ,Y_m=u(x_m)+\epsilon_m,
\end{cases}
\end{equation}
where $u(x_1),\ldots,u(x_m)$ are point evaluations of a single sample path. Note that we consider a fixed-design setting where $x_1,\ldots,x_m$ is fixed throughout. Instead of the nominal Bayes risk, we consider a `worst-case Bayes risk' with respect to all possible mis-specification on both $b$ and $\epsilon$, i.e., the whole data-generating process.
Our goal is to minimize the worst-case Bayes risk by solving
\begin{equation} \label{eq:dro}
\inf \limits_{\phi\in\mathcal M}\sup\limits_{P\in\mathcal{P},\W(P,P_0)\leq\delta} \mathbb{E}_P\left[\| b-\phi(Y_1,\ldots,Y_m)\|^2_{L^2(\mathcal{D})}\right],
\end{equation}
\label{nomencl:delta}\nomenclature{$\delta$}{radius of the Wasserstein Ball}{nomencl:delta}\label{nomencl:genb}\nomenclature{$b$}{denotes a generic choice of input process by nature}{nomencl:genb}where the adversary's admissible choice of $P$\label{nomencl:genP}\nomenclature{$P$}{denotes a generic choice of probability measure by nature}{nomencl:genP} is constrained by the Wasserstein distance $\mathds{W}(P,P_0)$ relative to the nominal measure $P_0$. 
%
%{\color{red}\sout{We will also consider the analogous robust formulation for estimating $b$ itself, as written in \eqref{eq:defminmaxintro}; we revisit this formulation specifically in Section~\ref{sec:inverserobust}. \ymmtd{Added signposting for the inverse problem, to keep this thread alive. Do you think there is any better way to do this?{\color{red}xz: this way is good to me!}} } xz: no longer needed}

The Wasserstein distance $\mathds{W}$\label{nomencl:wass}\nomenclature{$\mathds{W}$}{Wasserstein distance}{nomencl:wass} is constructed based on the optimal transport cost between $P$ and $P_0$, which is defined as follows.
\begin{definition}[Optimal transport cost]
The optimal transport cost between the probability measures $P$ and $P_0$ on $C(\mathcal{D}) \times \R^m$ is defined as
\begin{equation}\label{eq:wassdistdef}
    D_c(P,P_0)= \inf_{\pi}~\big\{\mathbb{E}_{\pi}[c((b,\epsilon),(b^0,\epsilon^0))]~:~\pi_{(b,\epsilon)} = P,~\pi_{( b^0,\epsilon^0)} = P_0 \big\},
\end{equation}
\label{nomencl:geneps}\nomenclature{$\epsilon$}{denotes a generic choice of observation errors by nature}{nomencl:geneps}where the infimum is taken over all couplings $\pi$ between $(b,\epsilon)$ and $(b^0,\epsilon^0)$ with marginals $P$ and $P_0$, and $c$ is some ground cost on $C(\mathcal{D}) \times \R^m$. 
\end{definition}
\label{nomencl:pi}\nomenclature{$\pi$}{coupling of probability measures}{nomencl:pi}\review{Here, our definition of optimal transport cost follows the general literature, see~\cite{ref:villani2008optimal}, with specification that the probability measure consists of a process component (namely $b$) and a finite-dimensional vector component (namely $\epsilon$).} The existence of an optimal coupling (i.e., an optimal solution to problem~\eqref{eq:wassdistdef}) is guaranteed whenever $c$ is non-negative and lower semi-continuous with respect to the product norm on the Polish space $(C(\mathcal{D})\times \mathbb{R}^m)^2$; see e.g.,~\cite[Theorem~4.1]{ref:villani2008optimal}. \review{The particular choice of $c$ is typically motivated by the geometric properties of the problem at hand. For instance, the geometry of finite-dimensional spaces differs from that of infinite-dimensional spaces.} In this paper, we use the ground cost function $c$ defined as\label{nomencl:c}\nomenclature{$c$}{ground cost function}{nomencl:c}
\begin{align*}
c((b,\epsilon),(b^0,\epsilon^0)) & =\|\epsilon-\epsilon^0\|_2^2
+ \|b-b^0 \|^2_{\mathcal{H}_w}  \\
& =  \sum_{i=1}^m (\epsilon_i-\epsilon^0_i)^2 + \sum_{n\geq1} \left(\langle b,e_n\rangle - \langle  b^0,e_n\rangle\right)^2w_n.
\end{align*}
To see that $c$ is lower semi-continuous, note that if $\|b_k-b_{\infty}\|_{C(\mathcal{D})}\to0$ and $\|\tilde b_k-\tilde b_{\infty}\|_{C(\mathcal{D})}\to0$ as $k\to\infty$, then $\liminf_{k\to\infty}\left(\langle b_k,e_n\rangle - \langle  \tilde b_k,e_n\rangle\right)^2w_n\geq\left(\langle b_\infty,e_n\rangle - \langle  \tilde b_\infty,e_n\rangle\right)^2w_n$ for each $n$ since $e_n\in C(\mathcal{D})$. Thus, we have
\[
\liminf_{k\to\infty}\sum_{n\geq1} \left(\langle b_k,e_n\rangle - \langle  \tilde b_k,e_n\rangle\right)^2w_n\geq\sum_{n\geq1} \left(\langle b_\infty,e_n\rangle - \langle  \tilde b_\infty,e_n\rangle\right)^2w_n.
\]
With this choice of the ground cost $c$, the optimal coupling between $P$ and $P_0$ exists for any $P$ and $P_0$. Note that $D_c$ is not a distance because it does not satisfy the triangle inequality, but its square root $\sqrt{D_c}$ is a Wasserstein-type distance on its domain of finiteness $\{P\in\mathcal{P}:D_c(P,P_0)<\infty\}$~\cite[Definition~6.1]{ref:villani2008optimal}. Throughout this paper, we will use \[\W \coloneqq \sqrt{D_c}.\]

\review{Our choice of the ground cost function $c$ penalizes discrepancies in both the observation error component (with a usual  $l_2$ norm in finite-dimensional Euclidean space) and the stochastic process component (with a Hilbert space norm), which to our knowledge, is novel in the literature on DRO. Allowing the adversary to modify the distribution of the additive noise $\epsilon_i$ can also be understood in part as compensating for the misspecification of the nominal forward operator $T$ \cite{kennedy2001bayesian}. One technical reason that we chose optimal transport cost instead of e.g. relative entropy or any other criterion that
requires the existence of a likelihood ratio is that $P$ and $P_0$ can have different support, which is needed when dealing with non-parametric models and permits adversarial distributions with different smoothness properties of the sample paths.}

It is important to stress that our definition of the Wasserstein distance $\W$ depends on the specification of the function class $\mathcal{H}_w$, or equivalently on the Hilbert space norm $\|\cdot\|_{\mathcal{H}_w}$. This dependence provides a convenient tool to control features such as the amount of roughness or smoothness that the adversary is allowed to inject in the sample path of the process $b$ (and hence of $u$).  Intuitively speaking, the sequence $w$ puts different penalties on the mass transportation of different ``modes'' of the spectral decomposition of the sample path of $b$. For example, if the sequence $w$ increases to infinity slowly, the adversary is under-penalized for moving mass corresponding to the higher ``modes,'' resulting in rougher sample paths of $b$. The modeling of the adversary's behavior thus conveniently reduces to the specification of the Hilbert norm $\|\cdot\|_{\mathcal{H}_w}$. 

We impose a final compatibility assumption between the operator $T$ and the adversarial cost introduced. 

\begin{assumption}[Operator and adversarial cost]\label{assmp:operator2} We select $\tilde w$ in Assumption~\ref{assmp:operator} such that $\tilde w_n = o(w_n)$ as $n \to \infty$.
\end{assumption}

Intuitively, Assumption~\ref{assmp:operator2} says that the operator is bounded even if the adversarial perturbations are slightly rougher than the adversarial choice. This assumption, we believe, is purely technical. The natural condition to impose is that the operator is bounded only on the chosen adversarial space. Our results hold under this more natural (and weaker) assumption in the case of standard Gaussian process regression, namely when $T$ equals the identity map. 

\begin{remark}
Several comments are in order regarding our framework:
\begin{enumerate}[label=(\roman*)]
    \item It is natural to consider random elements on a general Banach space $\mathbb{B}$ other than $C(\mathcal{D})$, e.g., by embedding the Banach space $\mathbb{B}$ in its second dual $\mathbb{B}^{\star\star}$ and identifying a Borel measurable random element $b$ in $\mathbb{B}$ with the stochastic process $\left(b^\star(b):b^\star\in\mathbb{B}^\star\right)$, but at the expense of additional technicality, see a discussion in~\cite[Section~2.3]{ref:vaart2008reproducing}. We choose $C(\mathcal{D})$ to illustrate our main conceptual contribution of a distributionally robust formulation of nonparametric regression and inverse problems. 
    \item The $L^2(\mathcal{D})$ norm in the objective function of the formulation~\eqref{eq:dro} can be replaced by another member in the hierarchy of the Hilbert space norms. However, the latter norm lacks the Lebesgue integral representation, especially for those with fractional power~\cite[Section~4.1]{ref:taylor1996partialone}, and we leave the extension to future work. 
    \item It is tempting to replace the Hilbert space norm with the $L^2(\mathcal{D})$ norm in the definition of the ground cost function $c$. However, point evaluations are not continuous under the $L^2(\mathcal{D})$ norm. Our proofs for the main results rely crucially on the continuity of point evaluations, and thus, we resort to the RKHS in Assumption~\ref{assmp:operator}.
\end{enumerate}
\end{remark}

\section{Main Results}\label{sec:mainresults}
\subsection{Strong Duality}
Our first main result is a minimax theorem, which states that one may interchange the infimum and supremum operators in the \review{inverse} problem~\eqref{eq:dro}.
\begin{theorem}[Strong duality for the \review{inverse} problem]\label{thm:swap}
Suppose that Assumptions~\ref{assmp:fullrank}--\ref{assmp:operator2} hold. For any $\delta>0$,
 strong duality holds, that is,
 \begin{align}
&\inf_{\phi\in\mathcal M}\sup_{P\in\mathcal{P},\W(P,P_0)\leq\delta} \mathbb{E}_{P}\left[\|b-\phi(Y_1,\ldots,Y_m)\|^2_{L^2(\mathcal{D})}\right]\notag\\
& \qquad = \sup_{P\in\mathcal{P},\W(P,P_0)\leq\delta}\inf_{\phi\in\mathcal M}\mathbb{E}_{P}\left[\|b-\phi(Y_1,\ldots,Y_m)\|^2_{L^2(\mathcal{D})}\right].\label{eq:thmstrongdualityregress}
\end{align}
\end{theorem}
\review{The full proof of Theorem~\ref{thm:swap} is deferred to Appendix A. To prove Theorem~\ref{thm:swap}, our technical contribution consists of a novel asymptotic approximation technique. First, the left hand side of~\eqref{eq:thmstrongdualityregress} can be trivially upper bounded by the min-max problem constraining $\phi$ as \textit{affine estimators} (see Definition~\ref{def:affinepredictor_revision}), namely
\begin{align}
&\inf_{\phi\in\mathcal M \textrm{ and affine}}\sup_{P\in\mathcal{P},\W(P,P_0)\leq\delta} \mathbb{E}_{P}\left[\|b-\phi(Y_1,\ldots,Y_m)\|^2_{L^2(\mathcal{D})}\right]\notag\\
& \qquad \geq  \inf_{\phi\in\mathcal M}\sup_{P\in\mathcal{P},\W(P,P_0)\leq\delta} \mathbb{E}_{P}\left[\|b-\phi(Y_1,\ldots,Y_m)\|^2_{L^2(\mathcal{D})}\right].\label{eq:thmstrongdualityregress_revision_left}
\end{align}
Second, define $\Span\{e_n\}_{n=1}^N$\label{nomencl:spanN}\nomenclature{$\Span\{e_n\}_{n=1}^N$}{linear subspace of spanned by the first $N$ basis vectors}{nomencl:spanN} as the (closed) linear subspace of $C(\mathcal{D})$ spanned by the basis vectors $(e_n:1\leq n\leq N)$\label{nomencl:N}\nomenclature{$N$}{number of basis vectors in finite-dimensional approximation}{nomencl:N}, the right hand side of~\eqref{eq:thmstrongdualityregress} can be lower bounded by constraining the ambiguity set, namely
\begin{align}
&\sup_{P\in\mathcal{P},\W(P,P_0)\leq\delta}\inf_{\phi\in\mathcal M}\mathbb{E}_{P}\left[\|b-\phi(Y_1,\ldots,Y_m)\|^2_{L^2(\mathcal{D})}\right]\notag\\
& \qquad \geq  \sup_{P\in\mathcal W_N(\delta)}\inf_{\phi\in\mathcal M}\mathbb{E}_{P}\left[\|b-\phi(Y_1,\ldots,Y_m)\|^2_{L^2(\mathcal{D})}\right],\label{eq:thmstrongdualityregress_revision_right}
\end{align}
where $\mathcal W_N(\delta)$ denotes the Wasserstein balls around $P_0$ arising from perturbations constrained in the first $N$ coordinates  (see equation~\eqref{eq:revision_Wassdelta} in the appendix). In the end, we argue that the left hand side of the relaxation~\eqref{eq:thmstrongdualityregress_revision_left} and the asympototic limit of the right hand side of the relaxation~\eqref{eq:thmstrongdualityregress_revision_right} coincides, and since they bound problem~\eqref{eq:thmstrongdualityregress} the strong duality follows.
} 

To proceed, we first show a strong duality result for a sequence of finite-dimensional approximations.  In particular, we consider truncating $P$ and $P_0$ into the space $\Span\{e_n\}_{n=1}^N\times\mathbb{R}^m$, and denote the induced measure as $Q^{(N)}$\label{nomencl:QuN}\nomenclature{$Q^{(N)}$}{induced measure of $P$ in $\mathbb{R}^{N+m}$}{nomencl:QuN} and $Q_0^{(N)}$\label{nomencl:Q0uN}\nomenclature{$Q_0^{(N)}$}{induced measure of $P_0$ in $\mathbb{R}^{N+m}$}{nomencl:Q0uN}, respectively. The truncation is through the coordinate projections after expanding functions in the basis $\{e_n\}_{n=1}^\infty$. Since the coordinate projections are bounded linear mappings, $Q_0^{(N)}$ is centered Gaussian.
Notice that the space $\Span\{e_n\}_{n=1}^N\times\mathbb{R}^m$ is isomorphic to $\mathbb{R}^{N+m}$, thus we view the truncated measures $Q^{(N)}$ and ~$Q_0^{(N)}$ as finite-dimensional measures on $\mathbb{R}^{N+m}$. \review{This reduction to finite-dimensional setting enables the use of existing results (e.g., Sion's minimax theorem~\cite{ref:sion1958minimax}).}

For convenience, denote 
\[
\Obj(\phi,P) = \mathbb{E}_{P}\left[\|b-\phi(Y_1,\ldots,Y_m)\|^2_{L^2(\mathcal{D})}\right],
\]
and by a slight abuse of notation we denote $\Obj(\phi,Q^{(N)})$ for the truncated measure. The related strong duality in finite-dimensions reads
\begin{equation}\label{eq:finitedimdualitymain}
\min_{\phi\in\mathcal{M}} \max_{Q^{(N)}:\W_N(Q^{(N)},Q_0^{(N)})\le \delta}~\Obj(\phi,Q^{(N)}) =\max_{Q^{(N)}:\W_N(Q^{(N)},Q_0^{(N)})\le \delta} \min_{\phi\in\mathcal{M}}~\Obj(\phi,Q^{(N)}),
\end{equation}
where $\W_N^2$\label{nomencl:WN2}\nomenclature{$\W_N$}{induced optimal transport cost on $\mathbb{R}^{N+m}$}{nomencl:WN2} is the induced optimal transport cost on $\mathbb{R}^{N+m}$
\[
\W_N^2(Q^{(N)},Q_0^{(N)})= \min_{\pi}\left\{\mathbb{E}_{\pi}[c_N(r,s)]: \pi_{r} = Q^{(N)},\pi_{s} = Q_0^{(N)}\right\}.
\]
Here $\pi_r$ and $\pi_s$ are projections onto the first and second components of the coupling $\pi$. In the definition of $\W_N^2$, $c_N$ is the induced cost function on $\mathbb{R}^{N+m}$ with\label{nomencl:cN}\nomenclature{$c_N$}{induced cost function on $\mathbb{R}^{N+m}$}{nomencl:cN}
\[
c_N(r,s) = \sum_{n=1}^N (r_n-s_n)^2 w_n + \sum_{j=1}^m (r_{N+j} - s_{N+j})^2 \quad \text{for any~} r,~s\in\mathbb{R}^{N+m},
\] 
\review{where recall that $w_n$ are the parameters of the Hilbert space norm $\|\cdot\|_{\mathcal{H}_w}$.} We note that in~\eqref{eq:finitedimdualitymain}, the minimizer in the $\phi$ variable and the maximizer in the $Q^{(N)}$ variable exist, which justifies the min and the max operators. Denote by $\phi^\star_N$ and $Q^\star_N$ the (unique) solution to the outer optimization problem in the left-hand and the right-hand side of~\eqref{eq:finitedimdualitymain}, respectively. Then $(\phi^\star_N,Q^\star_N)$\label{nomencl:nashfiniteopt}\nomenclature{$(\phi^\star_N,Q^\star_N)$}{pair of Nash equilibrium in finite-dimensional approximation}{nomencl:nashfiniteopt} is the (unique) pair of Nash equilibrium for problem~\eqref{eq:finitedimdualitymain} in the sense that
\[
\Obj(\phi^\star_N,Q^\star_N) = \min_{\phi\in\mathcal{M}}~\Obj(\phi,Q^\star_N) = \max_{Q^{(N)}:\W_N(Q^{(N)},Q_0^{(N)})\le \delta} ~\Obj(\phi^\star_N,Q^{(N)}).
\]
The structural properties of the optimal solutions reveal that $Q^\star_N$ is a centered Gaussian distribution, and $\phi^\star_N$ is a linear prediction rule (see details in the full proof). Namely, for any $x \in \mathcal{D}$,
\begin{align*}
\phi^\star_N(Y_1,\ldots,Y_m)(x) & =\mathbb{E}_{Q^\star_N}\left[b(x)\big|Y_1,\ldots,Y_m\right]\\
& =  \left(k_\epsilon^{(N)}(x,x_1),\ldots,k_\epsilon^{(N)}(x,x_m)\right)\cdot (K_\epsilon^{(N)})^{-1} \cdot \left(Y_1,\ldots,Y_m\right)^\top,
\end{align*}
where $k_\epsilon^{(N)}(x,x_i) = \mathbb{E}_{Q^\star_N}[b(x)Y_i]$\label{nomencl:covobsPbyfinite}\nomenclature{$k_\epsilon^{(N)}(x,x_i)$}{covariance function of $b(x)$ and $Y_i$ under $Q^\star_N$}{nomencl:covobsPbyfinite}, and $K_\epsilon^{(N)} = \left(\mathbb{E}_{Q^\star_N}[Y_iY_j]\right)_{ij} \in \PD^{m}$\label{nomencl:covobsPfinite}\nomenclature{$K_\epsilon^{(N)}$}{covariance matrix of observations $Y$ under $Q^\star_N$}{nomencl:covobsPfinite} is invertible.\label{nomencl:pdmat}\nomenclature{$\PD^{m}$}{positive-definite matrix of dimension $m$-by-$m$}{nomencl:pdmat}
Similar finite-dimensional duality results have been established in~\cite{ref:nguyen2019bridging,ref:shafieezadeh2018wasserstein}, but with a different definition of the Wasserstein distance and the ambiguity set. 

\review{Having established the strong duality in finite-dimensions, the rest of the proof then argues that the left hand side of the relaxation~\eqref{eq:thmstrongdualityregress_revision_left} can be upper bounded by the value of the finite-dimensional min-max problem with an error of $o(1)$, and the right hand side of the relaxation~\eqref{eq:thmstrongdualityregress_revision_right} can be lower bounded by the value of the finite-dimensional max-min problem.  Therefore, as the dimension $N$ grows to infinity, the approximation error becomes negligible. This limiting argument is a key technical novelty, which is achieved by a careful analysis of the tail modes of the spectral decomposition.
}
We state one intermediate result that guarantees the approximation of objective values.
\begin{proposition}[Approximation of objective values]\label{prop:fdvalueconverges}
Let $\phi^\star_N$ be the (unique) solution to the min-max problem in~\eqref{eq:finitedimdualitymain}, and $Q^\star_N$ be the (unique) solution to the max-min problem in~\eqref{eq:finitedimdualitymain}. We have
\[
\Obj(\phi^\star_N,Q^\star_N) = \sup_{P\in\mathcal{P},\W(P,P_0)\leq\delta} \Obj(\phi^\star_N,P) + o(1)= \inf_{\phi\in\mathcal M}\sup_{P\in\mathcal{P},\W(P,P_0)\leq\delta}\Obj(\phi,P)+o(1),
\]
asymptotically as the number of basis vectors in the approximation tends to infinity, i.e., as $N\to\infty$.
\end{proposition}

%%%%%%%%%%%%%%%%%%%%%%%%%%%%%%%%%%%%%%%%

\subsection{Existence, Uniqueness, and Construction of the Nash Equilibrium}\label{sec:nashequi}
In this part, we show that~\eqref{eq:thmstrongdualityregress} admits a unique pair of Nash equilibrium under certain conditions. Recall the Nash equilibrium $(\phi_N^\star, Q^\star_N)$ corresponding to the finite-dimensional approximation in the last section. We now add the tail of $P_0$ to $Q^\star_N$, namely, we denote $P_N^\star\in\mathcal{P}$\label{nomencl:PNstar}\nomenclature{$P_N^\star$}{the probability measure composing of $Q^\star_N$ and the tail of $P_0$}{nomencl:PNstar} as the measure under which the two random elements 
\begin{equation}\label{eq:independentassertionclaim2}
(\{\langle b,e_n\rangle\}_{n=1}^N,\epsilon)\quad\textrm{ and }\quad\{\langle b,e_{n}\rangle\}_{n=N+1}^\infty
\end{equation}
are independent. Moreover, we have that
\[
 \mathcal L^{P^\star_N} \big(\{\langle b,e_{n}\rangle\}_{n=N+1}^\infty\big) = \mathcal L^{P_0} \big(\{\langle b^0,e_{n}\rangle\}_{n=N+1}^\infty\big),
\]
where $\mathcal{L}^P$ denotes the law under $P$\label{nomencl:LP}\nomenclature{$\mathcal{L}^P$}{the law under $P$}{nomencl:LP}. 
We can extract a subsequence from $P^\star_N$, $N\geq1$ by a compactness argument. 
\begin{proposition}[Compactness of the ambiguity set]\label{prop:weakcompact}
% The set $\{P:\W(P,P_0)\leq\delta\}$ is sequentially weakly compact \ssw{With respect to which topology? I.e., where does the weak topology at hand come from?}. Namely,
Under the conditions of Theorem~\ref{thm:swap}, for every sequence $P_N\in\mathcal{P},N\geq1$ that satisfies $\W(P_N,P_0)\leq\delta$, there exists a weakly convergent subsequence $P_{N_l}$, $l\geq1$ with $P_{N_l}\Rightarrow P_\infty$, such that the limit $P_\infty\in\mathcal{P}$ also satisfies $\W(P_\infty,P_0)\leq\delta$.
\end{proposition}
By Proposition~\ref{prop:weakcompact}, we find a weakly convergent subsequence of $P^\star_N$, denoted as $P^\star_{N_l}$, $l\geq1$, with a limit $P^\star_\infty$ that belongs to the feasible set. The subsequence $P_{N_l}^\star$ is centered Gaussian, thus the limit $P^\star_\infty$ is also centered Gaussian. To see this, consider any bounded linear functional $F: C(\mathcal{D})\times \mathbb{R}^m\to\mathbb{R}$, and construct the Skorohod representations $Z_{N_l}\sim P_{N_l}^\star$ and $Z_\infty\sim P^\star_\infty$, where $Z_{N_l}$ converges to $Z_\infty$ almost surely~\cite{ref:skorokhod1956limit} on the common probability space. It follows that $F(Z_{N_l})$ converges to $F(Z_\infty)$ almost surely, and we note that the limit of centered univariate Gaussians must also be centered univariate Gaussian.

Define the matrix \review{$K_{\epsilon,P^\star_\infty} = \left(\mathbb{E}_{P^\star_\infty}[Y_iY_j]\right)_{ij} \in \PSD^m$}.\label{nomencl:covobsP}\nomenclature{$K_{\epsilon,P}$}{covariance matrix of observations $Y$ under $P$}{nomencl:covobsP}
Under the condition that $K_{\epsilon,P^\star_\infty}$ is invertible, the solution $\phi^\star_\infty$ to $\min_\phi \Obj(\phi,P^\star_\infty)$ is well-defined as
\begin{align}\label{phi-star-inf}
\phi^\star_{\infty}(Y_1,\ldots,Y_m)(x) &= \mathbb{E}_{P^\star_\infty}\left[b(x)|Y_1,\ldots,Y_m\right] \\
&= \left(k_{\epsilon,P^\star_\infty}(x,x_1),\ldots,k_{\epsilon,P^\star_\infty}(x,x_m)\right)\cdot (K_{\epsilon,P^\star_\infty})^{-1} \cdot \left(Y_1,\ldots,Y_m\right)^\top,\notag
\end{align}
where $k_{\epsilon,P^\star_\infty}(x,x_i) = \mathbb{E}_{P^\star_\infty}[b(x)Y_i]$.\label{nomencl:covobsPby}\nomenclature{$k_{\epsilon,P}(x,x_i)$}{covariance function of $b(x)$ and $Y_i$ under $P$}{nomencl:covobsPby} The main result of this section is the following theorem.
\begin{theorem}[Nash equilibrium] \label{thm:existnash} Suppose that the conditions of Theorem~\ref{thm:swap} hold.
Let $P^\star_\infty$ denote the limit of a weakly convergent subsequence of $P^\star_N,N\geq1$, and assume that $K_{\epsilon,P^\star_\infty}$ is invertible under $P^\star_\infty$.  Then, for $\phi^\star_\infty$ given by~\eqref{phi-star-inf}, the pair $(\phi^\star_\infty,P^\star_\infty)$ forms a Nash equilibrium
to problem~\eqref{eq:thmstrongdualityregress}, i.e.,
\[
\Obj(\phi^\star_\infty,P^\star_\infty) = \min_{\phi\in\mathcal{M}}\Obj(\phi,P^\star_\infty) = \max_{P\in\mathcal{P}:\W(P,P_0)\leq\delta}\Obj(\phi^\star_\infty,P).
\]
Thus, $(\phi^\star_\infty,P^\star_\infty)$ represents a pair of equilibrium strategies where neither the decision-maker nor the adversary, the two players of the zero-sum game, benefits from changing their strategy.
Moreover, the pair $(\phi^\star_\infty,P^\star_\infty)$ is unique, with components respectively given by the pointwise limit \[\lim_{N\to\infty}\phi^\star_N(Y_1,\ldots,Y_m)(x) = \phi^\star_\infty(Y_1,\ldots,Y_m)(x)~~~ \forall~Y_1,\ldots,Y_m\in \R,~~~ x\in\mathcal{D},\]
and the weak limit $P^\star_N \Rightarrow P^\star_\infty$ as $N\to\infty$.
\end{theorem}
\label{nomencl:nashinfopt}\nomenclature{$(\phi^\star_\infty,P^\star_\infty)$}{pair of Nash equilibrium in infinite-dimensions}{nomencl:nashinfopt}One important consequence of Theorem~\ref{thm:existnash} is that the worst-case distribution involves a modified Gaussian process with potentially rougher paths than the prior. Hence, the robustified decision, as the conditional mean of the worst-case distribution, remains affine in the observations and, therefore, is tractable. \review{We formally state the following theorem.
\begin{theorem}[Nash equilibrium property]\label{thm:nashproperty_revision}
Under the notation of Theorem~\ref{thm:existnash}, the robustified decision $\phi^\star_\infty$ is an affine function in the observations, and the worst-case distribution $P^\star_\infty$ is a Gaussian process.
\end{theorem}
}

The  convergence statement in Theorem~\ref{thm:existnash} readily gives rise to an algorithm for computing the associated Nash equilibrium. 

The result of Theorem~\ref{thm:existnash} requires that $K_{\epsilon,P^\star_\infty}$ is invertible for the limit of some weakly convergent subsequence of $P^\star_N$. Fortunately, due to the following proposition, it is not difficult to check whether this condition holds in practice.
\begin{proposition}\label{prop:nasheqparity}
Either one (and only one) of the following cases occurs.
\begin{enumerate}
    \item there exists a weakly convergent subsequence of $P^\star_N$, with the limit denoted by $P^\star_\infty$, such that the matrix $K_{\epsilon,P^\star_\infty}$ is invertible under $P^\star_\infty$, or
    \item the sequence of determinants $\textrm{det}(K^{(N)}_\epsilon)\to0$ as $N\to\infty$.
\end{enumerate}
\end{proposition}
% Therefore, the results of Theorem~\ref{thm:existnash} hold once we rule out that $\textrm{det}(K^{(N)}_\epsilon)\to0$ as $N\to\infty$, which can be verified numerically.
 To conclude this section, we provide a sufficient condition to ensure that the first case of Proposition~\ref{prop:nasheqparity} occurs.
% Note that $
%  \left(\mathbb{E}_{P_0}[Y_iY_j]\right)_{ij}\succeq\sigma^2 I_{m\times m}$, since $\epsilon_i^0$ are independent $\mathcal{N}(0,\sigma^2)$ noise under $P_0$. We have
\begin{lemma}[Invertibility]\label{lem:smalldelta}
There exists a strictly positive constant $\delta_0$ that depends on $(T,m,(x_i)_i,\mathcal{H}_w,\mathcal{H}_{\tilde w},\sigma^2)$ such that for any $\delta<\delta_0$ and $P$ satisfying $\W(P,P_0)\leq\delta$, the matrix $\left(\mathbb{E}_{P}[Y_iY_j]\right)_{ij}$ is invertible. %In particular, $K_\epsilon$ is invertible.
% \begin{equation}\label{eq:matrixineqution}
%  \left(\mathbb{E}_{P}[(u(x_i)+\epsilon_i)(u(x_j)+\epsilon_j)]\right)_{ij}
% \end{equation}
% is invertible.
\end{lemma}

% \begin{lemma}\label{lem:subrkhsofprior}
% If $w_n=\Omega(\kappa_n^{-2})$ as $n\to\infty$, then $K_\epsilon$ is invertible.
% \end{lemma}
%%%%%%%%%%%%%%%%%%%%%%%%%%%%%%%%%%%%%%%
\subsection{Strong Duality for the \review{Regression} Problem}
\label{sec:inverserobust}

\review{Problem~\eqref{eq:nominalregression} is a linear inverse problem where our primary interest lies in recovering the unknown input $b$. Alternatively, we mention in passing a distributionally robust formulation for regressing $u=T(b)$ given noisy pointwise observations (i.e., a regression problem)}. Under the observation system~\eqref{eq:system}, the goal of the decision-maker is to seek a nonparametric predictor $\phi_u\in\mathcal{M}$\label{nomencl:phiu}\nomenclature{$\phi_u$}{non-parametric predicator for $u$}{nomencl:phiu} that minimizes the worst-case objective
% Alternatively, we also consider the inverse problem
% \[
% \min_{\phi_b\in\mathcal{M}}\mathbb{E}_{P_0}\left[\|b(\cdot)-\phi_b(\cdot;Y(x_1),\ldots,Y(x_m))\|^2_{L^2(\mathcal{D})}\right],
% \]
% whose linear prediction rule is given by
% \[
% \phi_b(x;Y(x_1),\ldots,Y(x_m)) = \left(k_b(x,x_1),\ldots,k_b(x,x_m)\right)\cdot (K)^{-1} \cdot \left(Y(x_1),\ldots,Y(x_m)\right)^\top,
% \]
% where $k_b(x,x_j) = \mathbb{E}_{P_0}[b(x)Y(x_j)]$ due to Gaussianity.
\begin{equation} \label{eq:droinverse}
\inf \limits_{\phi_u\in\mathcal M}\sup\limits_{P\in\mathcal{P},\W(P,P_0)\leq\delta} \mathbb{E}_P\left[\| u-\phi_u(Y_1,\ldots,Y_m)\|^2_{L^2(\mathcal{D})}\right],
\end{equation}
where nature's admissible choice of $P$ is constrained by the Wasserstein distance $\mathds{W}(P,P_0)$ constructed from~\eqref{eq:wassdistdef}. We state the strong duality of~\eqref{eq:droinverse}.

\begin{theorem}[Strong duality for the regression problem]\label{thm:swapinverse} 
Suppose that Assumptions~\ref{assmp:fullrank}--\ref{assmp:operator2} hold. For any $\delta>0$,
\begin{align}
&\inf_{\phi_u\in\mathcal M}\sup_{P\in\mathcal{P},\W(P,P_0)\leq\delta} \mathbb{E}_{P}\left[\|u-\phi_u(Y_1,\ldots,Y_m)\|^2_{L^2(\mathcal{D})}\right]\notag\\
&\qquad=\sup_{P\in\mathcal{P},\W(P,P_0)\leq\delta}\inf_{\phi_u\in\mathcal M}\mathbb{E}_{P}\left[\|u-\phi_u(Y_1,\ldots,Y_m)\|^2_{L^2(\mathcal{D})}\right].\label{eq:thmstrongdualityinverse}
\end{align}
\end{theorem}
The proof of Theorem~\ref{thm:swapinverse} works verbatim as that of Theorem~\ref{thm:swap}. Though strong duality holds for both the regression and the inverse problems under our formulations, we note that for ill-posed inverse problems in the Bayesian nonparametrics literature, the minimax rate for estimating $b$ is slower than the minimax rate for estimating $u$ \cite{D95,C11,CT02,KVV11}. 
% {\color{red}The condition in Theorem~\ref{thm:swapinverse} is more restrictive than Theorem~\ref{thm:swap} since $T$ induces a certain degree of smoothing to the sample path on $u$, and thus the problem of recovering $b$ is necessarily harder than the problem of regressing on $u$.} 
%\ssw{yes, we should discuss here that the minimax rates are slower for estimating $b$ from data. So basically, here we need that $b\sim P$ itself is already supported on $C(\mathcal D)$, for any perturbed prior $\mathbb W(P,P_0)\le \delta$.}

Regarding the Nash equilibrium associated with~\eqref{eq:thmstrongdualityinverse}, it is not hard to see, after examining the proof of Theorem~\ref{thm:existnash}, that we can develop the same theory verbatim to that of Section~\ref{sec:nashequi}. For ease of exposition,  we suppress the details here.

%%%%%%%%%%%%%%%%%%%%%%%%%%%%%%%%%%%%%

\section{Some Examples}\label{sec:examples}
In this section, we give several examples that illustrate the applicability of our general framework. We restrict our attention to Mat\'ern process priors and the Sobolev-type space of perturbations given by Examples~\ref{ex:matern} and~\ref{ex:perturbation}, respectively. We assume the relation $\kappa_n =\lambda_n^{-\alpha/2}$ (so that $\kappa=0$ in Example~\ref{ex:matern}) and $w_n = \lambda_n^{\beta}$, where $\alpha>\frac{d}{2}$ and $\beta>\frac{d}{2}$. 

\begin{example}
[Gaussian process regression]\label{ex:gpregress}
By choosing $T$ as the identity operator, we recover Gaussian process regression. Assumptions ~\ref{assmp:operator}--\ref{assmp:operator2} are satisfied for $\tilde w_n=\lambda_n^{\tilde\beta}$ and any $\frac{d}{2}<\tilde\beta<\beta$. If $\mathcal{D}$ is the one-dimensional interval $[0,1]$, the eigenvalues are $\lambda_n = n^2\pi^2$, and the eigenfunctions are
\[
e_n(x) = \sqrt{2}\sin(n\pi x) \quad \forall x\in[0,1].
\]
\end{example}

\begin{example}[Laplace equation]
The Laplace equation with a homogeneous Dirichlet boundary condition is
\[
\begin{cases}
\bigtriangleup u(x) = b(x) & \forall x\in\mathcal{D}^o,\\
u(x) = 0 & \forall x\in\partial\mathcal{D},
\end{cases}
\]
where $\mathcal{D}^o$ is the interior of $\mathcal{D}$, while $\partial\mathcal{D}$ denotes its boundary.
We have that the forward map $T$ is the inverse-Laplacian operator, thus,
\[
T(f) = \sum_{n\geq1}-\lambda_n^{-1}\langle f,e_n\rangle e_n.
\]
It is straightforward to see that Assumptions~\ref{assmp:operator}--\ref{assmp:operator2} are satisfied for $\tilde w_n=\lambda_n^{\tilde\beta}$ and any $\frac{d}{2}<\tilde\beta<\beta$. Similarly, one could treat other elliptic inverse problems.
\end{example}

\begin{example}[Heat equation] The one-dimensional homogeneous heat equation without source is
\[
\begin{cases}
u_t =  u_{xx} &  0<x<1,\\
u(x,0) = b(x) & 0<x<1,\\
 u(0,t)=u(1,t) = 0 & t\geq0,
\end{cases}
\]
where $u(\cdot,t)$ is the temperature profile at time $t$, and $b$ is the initial condition.
By separation of variables, the solution to the heat equation is
\[
u(x,t) = \sum_{n\geq1}e^{-n^2\pi^2t}\langle b,e_n\rangle e_n.
\]
Therefore, the (time-dependent) forward map $T$ satisfies, for any $t\geq0$, 
\[
T(f) = \sum_{n\geq1}e^{-n^2\pi^2t}\langle f,e_n\rangle e_n.
\]
Assumptions~\ref{assmp:operator}--\ref{assmp:operator2} are satisfied for $\tilde w_n=\lambda_n^{\tilde\beta}$ and any $\frac{d}{2}<\tilde\beta<\beta$.
\end{example}

\begin{example}[Radon transform in the plane] 
The Radon transform of a function $f$ is the function
\[
T(f)(s,\omega) = \int_{-\infty}^{\infty} f(s\omega + t\omega^\perp)\mathrm{d}t,\quad s\in\mathbb{R},~ \omega\in\mathcal{S}^1,
\]
where $\mathcal{S}^1$ is the unit circle and $\omega^\perp$ is the vector in $\mathcal{S}^1$ obtained by rotating $\omega$ counterclockwise by $90^\circ$. Recall we consider $f$ to be supported in a compact domain $\mathcal{D}\subset\mathbb{R}^2$, and thus $T(f)$ vanishes outside a compact subset of $\mathbb{R}\times \mathcal{S}^1$. It is straightforward to see that $T:C(\mathcal{D})\to C(\mathbb{R}\times \mathcal{S}^1)$ is linear and bounded, where we allow a slight modification of our framework since $T$ takes value in a different space from its domain of definition. By~\cite[Theorems II.5.1 and II.5.2]{ref:natterer2001mathematics}, the Radon transform satisfies the following Sobolev estimate: 
\[
\|T(f)\|_{H^{\tilde\beta}(\mathbb{R}\times\mathcal{S}^1)}\leq C_{\tilde w}\|f\|_{\mathcal{H}_{\tilde w}}\quad\forall f\in \mathcal{H}_{\tilde w},
\]
for $\tilde w_n=\lambda_n^{\tilde\beta}$ and any $1=\frac{d}{2}<\tilde\beta<\beta$, where $H^{\tilde\beta}(\mathbb{R}\times\mathcal{S}^1)$ is the usual order-$\tilde\beta$ Sobolev space. Identifying $\mathcal{S}^1$ with $[0,2\pi)$, and by the Sobolev embedding theorem~\cite[Proposition 4.1.3]{ref:taylor1996partialone}, we see that point evaluations in $H^{\tilde\beta}(\mathbb{R}\times\mathcal{S}^1)$ are continuous. Our theory in Section~\ref{sec:mainresults} applies with this slight modification of Assumptions \ref{assmp:operator}--\ref{assmp:operator2}, after inspecting the proofs.
\end{example}

%%%%
\section{Numerical Experiments}\label{sec:numerics}
Our focus in this paper is formulating a min-max framework for regression and linear inverse problems in an infinite-dimensional setting and elucidating key theoretical properties of this formulation. We now present some numerical experiments that offer further insights into the properties of the Nash equilibrium. 
In particular, we compute the Nash equilibrium $(\phi^\star_N,Q^\star_N)$ of a finite-dimensional approximation of the robust estimation problem (see definitions in Section~\ref{sec:mainresults}) with $N=200$, by an adaptation of the Frank-Wolfe algorithm in~\cite{ref:shafieezadeh2018wasserstein}. \review{The concrete algorithm as developed in~\cite{ref:shafieezadeh2018wasserstein} and the details of our adaptation is discussed in Appendix~\ref{sec:algodetail_revision}.} Under the conditions of Theorem~\ref{thm:existnash}, the use of the finite-dimensional Nash equilibrium is justified. When the condition fails, Proposition~\ref{prop:fdvalueconverges} still guarantees that the objective values of the finite-dimensional games converge to those of the \textit{in}finite-dimensional games. We illustrate our results for Gaussian process regression on the unit interval $\mathcal{D}=[0,1]$, as described in Example~\ref{ex:gpregress}. Throughout this section we set $\kappa_n =\lambda_n^{-\alpha/2}$ and $w_n = \lambda_n^{\beta}$. Recall that $(\kappa_n)_{n \geq 1}$, and hence $\alpha$, control the smoothness of the nominal prior, while $(w_n)_{n \geq 1}$, and hence $\beta$, control the roughness of the adversarial perturbations to the prior via the $\Vert \cdot \Vert_{\mathcal{H}_w}$ component of the transport cost.

First we fix a set of baseline parameters: $\alpha=2$ and $\beta = 0.51$; $\delta^2 = 0.1$ (size of the Wasserstein ball); and $\sigma = 0.1$ (observational noise magnitude). We choose $m=10$ design points $(x_i)_{i=1}^{10}$ equispaced on either the $(0,1)$ or $(0,0.5)$ intervals, excluding the endpoints. More specifically, we choose $x_i = i / 11$ and $x_i = i/22$ respectively. 
Below we visualize the prior and posterior covariances of $b$ under the nominal and worst-case measures $Q_0^{(N)}$ and $Q_N^\star$; i.e., $\mathbb{C}\text{ov}_{P}[b]$ and $\mathbb{C}\text{ov}_{P}[b \, \vert \, y_1, \ldots, y_m]$ for $P \in \{Q_0^{(N)}, Q_N^\star \}$. Since both the nominal and the worst-case measures are Gaussian, the posterior covariances do not depend on the realization of the data. 

\begin{figure}[!tb]
    \centering
    \subfigure[Nominal prior]{
		\includegraphics[width=0.45\columnwidth]{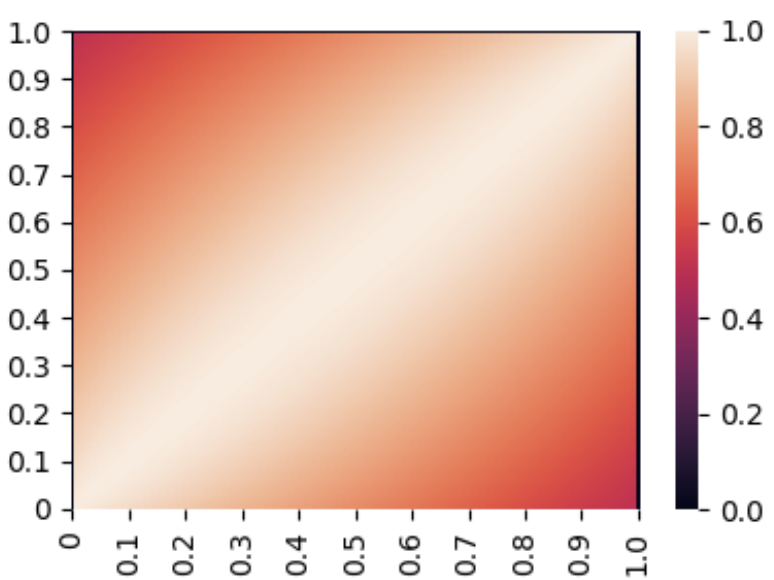}} \hspace{1mm}
	\subfigure[Worst-case prior]{
	\includegraphics[width=0.45\columnwidth]{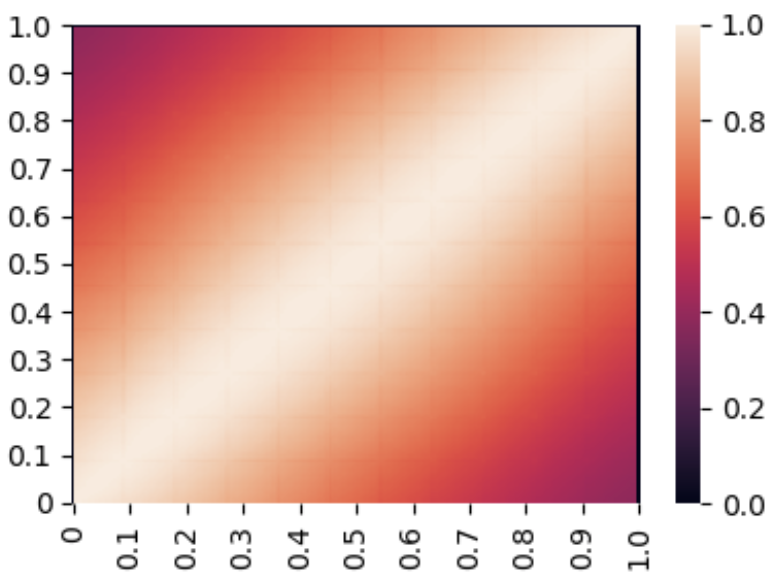}}
	    \hspace{1mm}
	\subfigure[Nominal posterior]{
		\includegraphics[width=0.45\columnwidth]{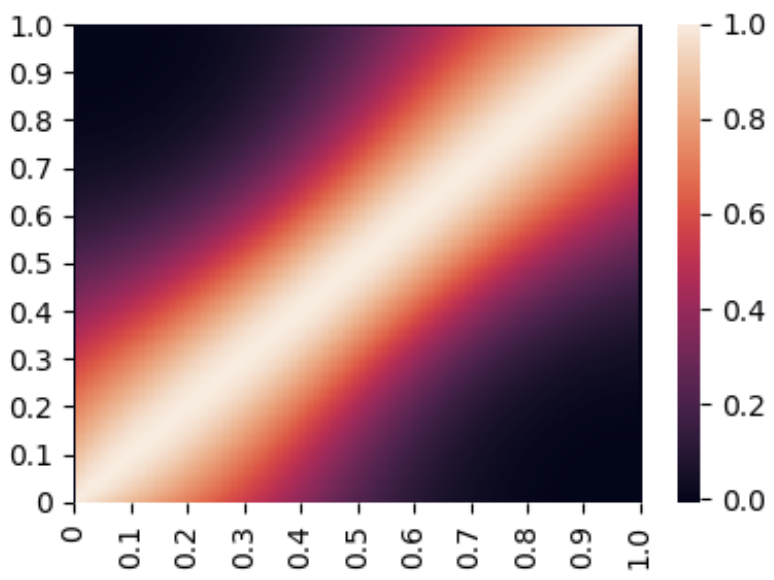}}
	\subfigure[Worst-case posterior]{
		\includegraphics[width=0.45\columnwidth]{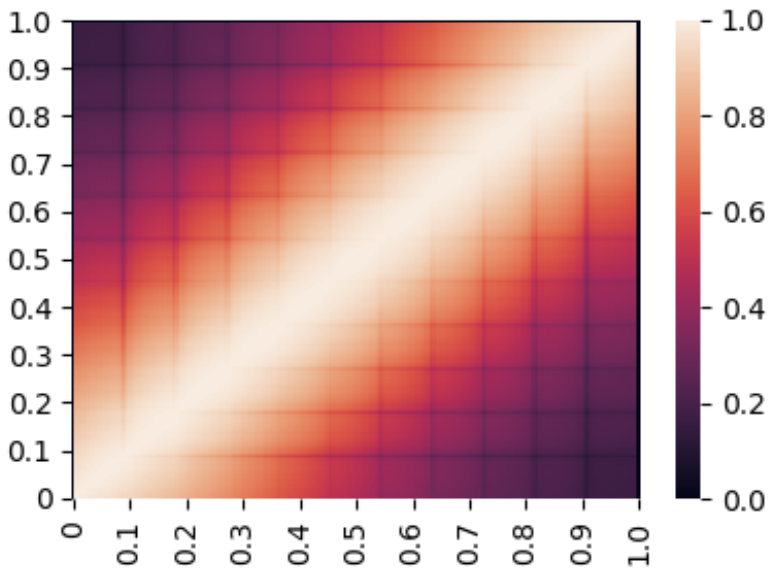}} 
	\caption{Correlation functions on $[0,1]^2$ with 10 design points equispaced on $(0,1)$.}
    \label{fig:heatmapwhole}
\end{figure}

\begin{figure}[h]
    \centering
    \subfigure[Nominal prior]{
		\includegraphics[width=0.45\columnwidth]{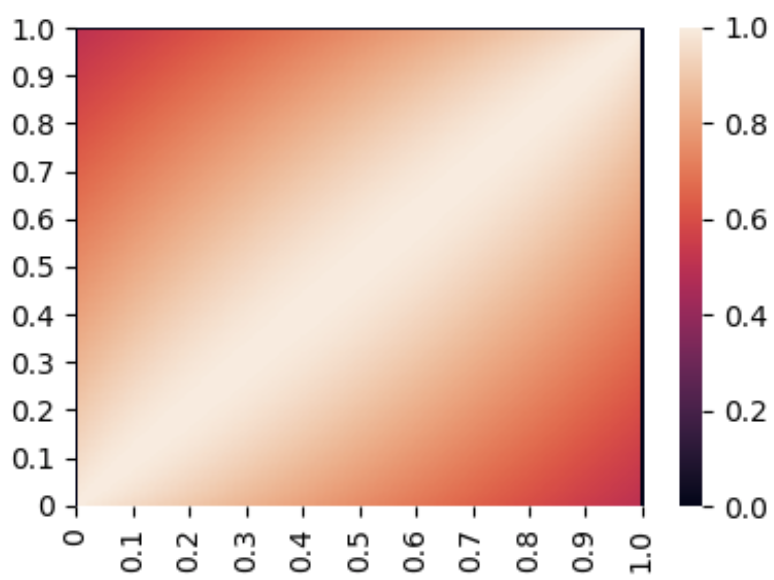}} \hspace{1mm}
	\subfigure[Worst-case prior]{
	\includegraphics[width=0.45\columnwidth]{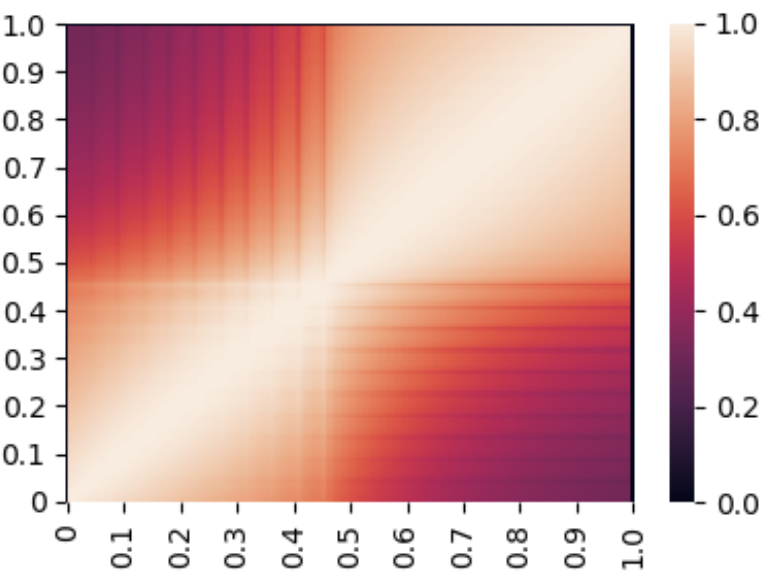}}
	    \hspace{1mm}
	\subfigure[Nominal posterior]{
		\includegraphics[width=0.45\columnwidth]{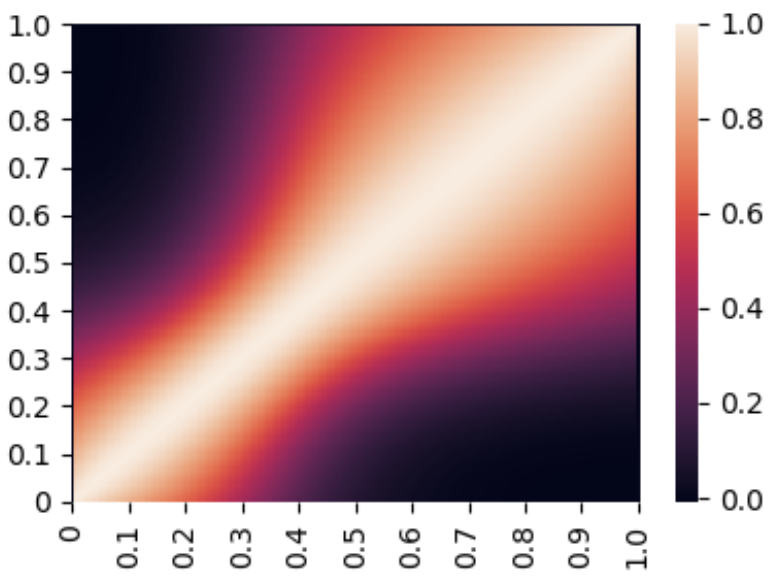}}
	\subfigure[Worst-case posterior]{
		\includegraphics[width=0.45\columnwidth]{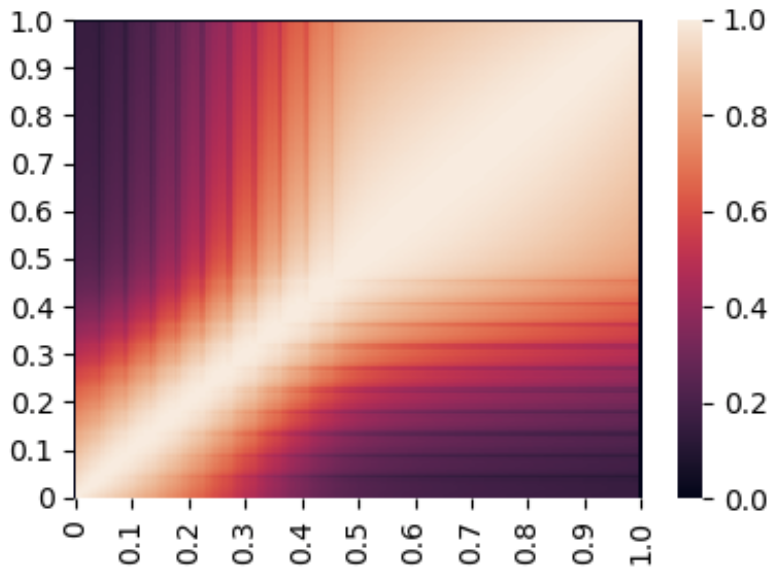}} 
	\caption{Correlation functions on $[0,1]^2$ with 10 design points equispaced on  $(0,0.5)$.}
    \label{fig:heatmaphalf}
\end{figure}

\begin{figure}[h]
    \centering
    \subfigure[Nominal prior]{
		\includegraphics[width=0.45\columnwidth]{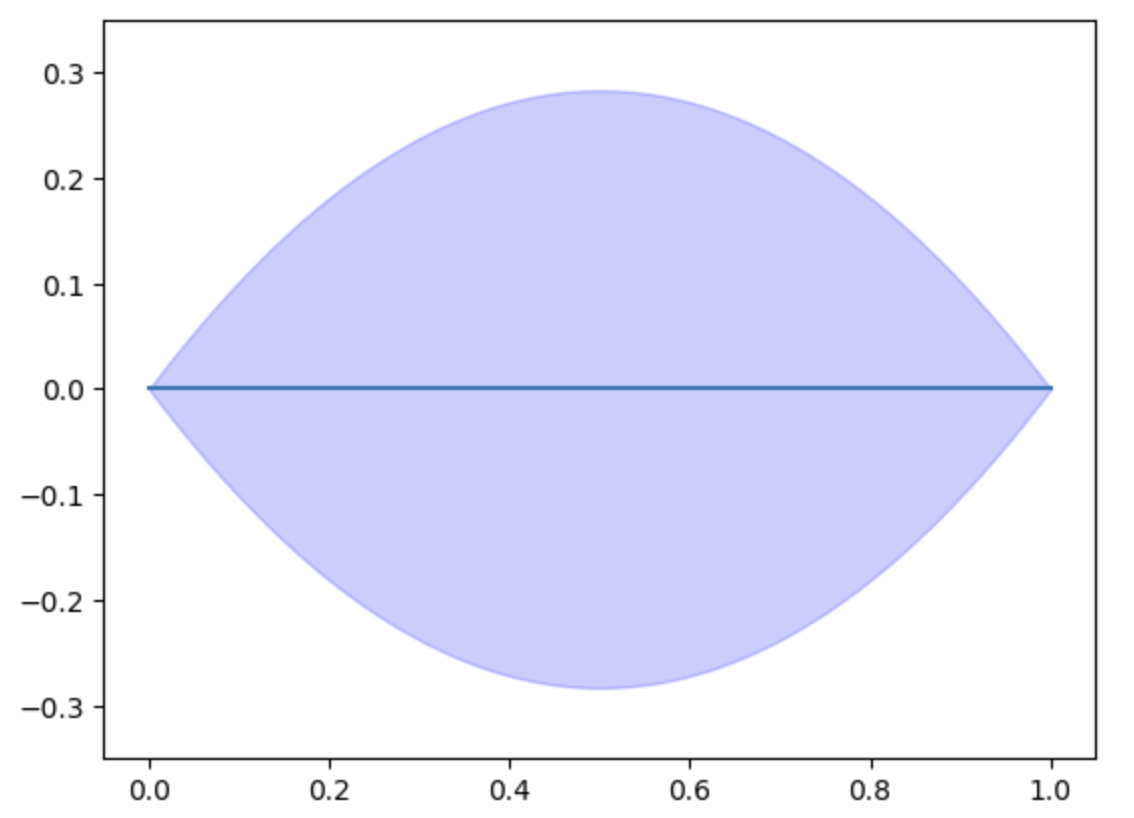}} \hspace{1mm}
	\subfigure[Worst-case prior]{
	\includegraphics[width=0.45\columnwidth]{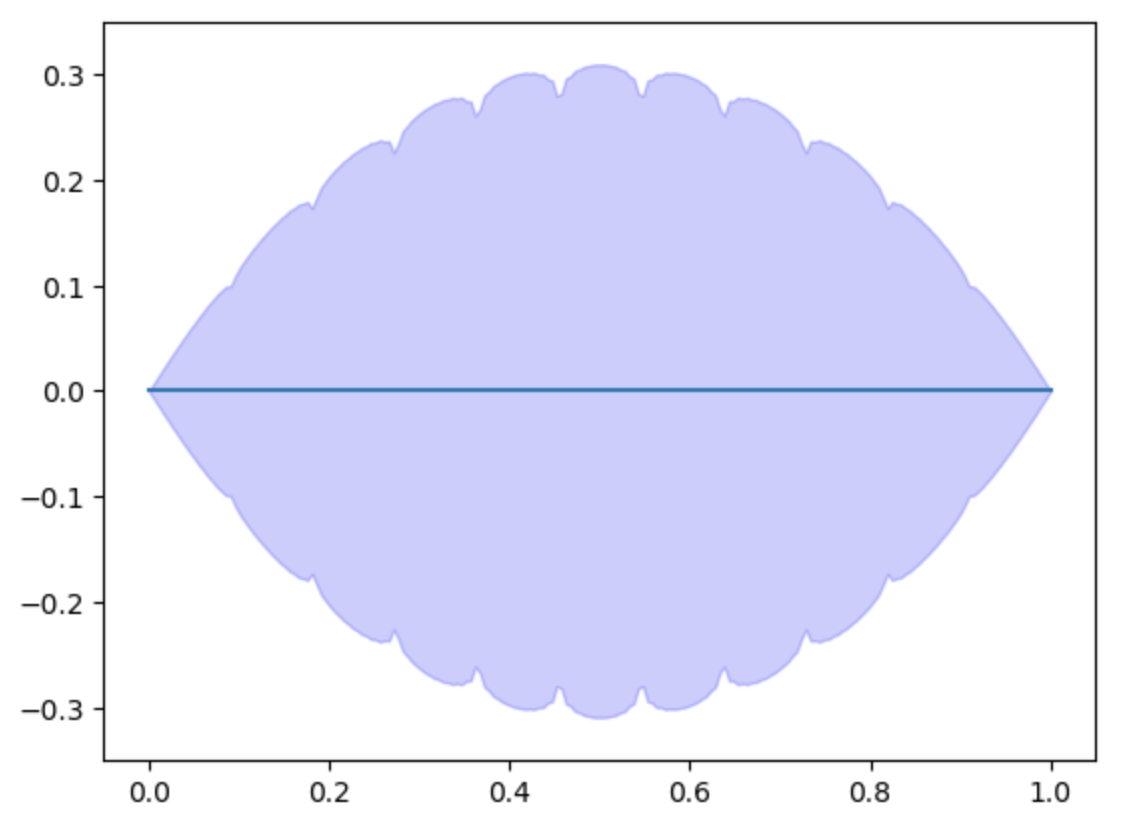}}
	    \hspace{1mm}
	\subfigure[Nominal posterior]{
		\includegraphics[width=0.45\columnwidth]{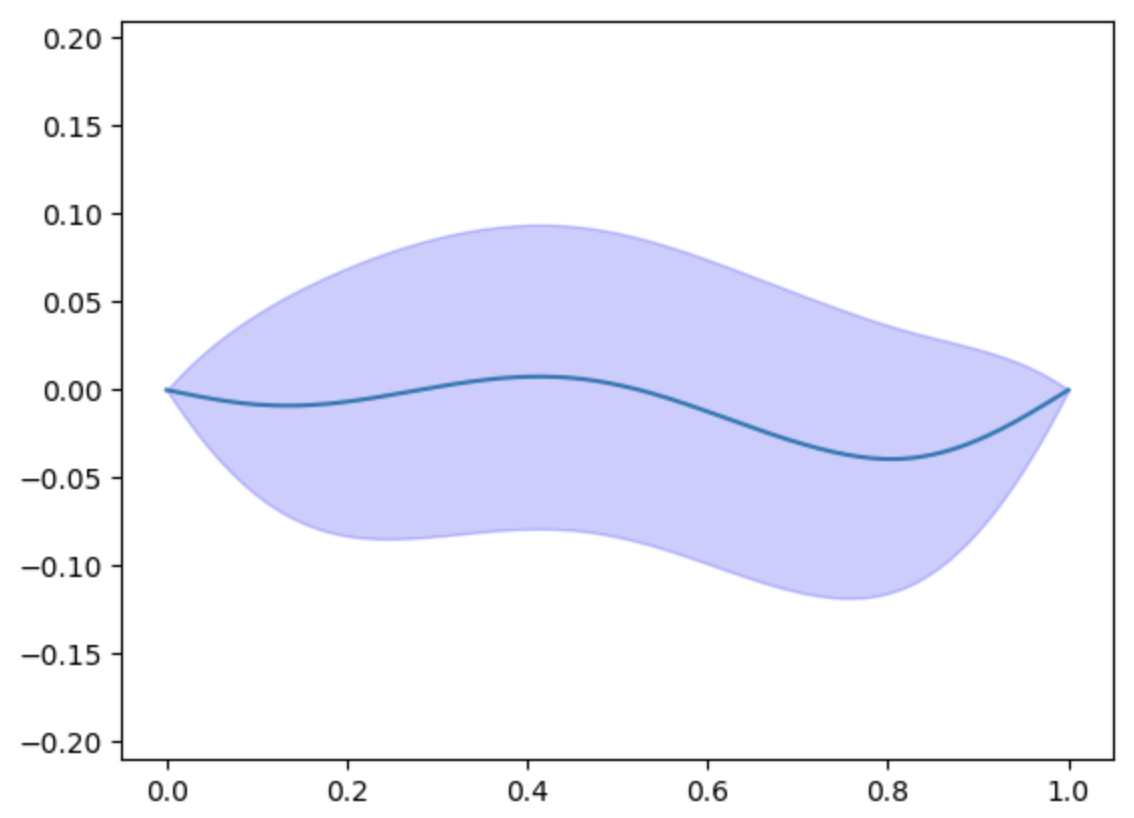}}
	\subfigure[Worst-case posterior]{
		\includegraphics[width=0.45\columnwidth]{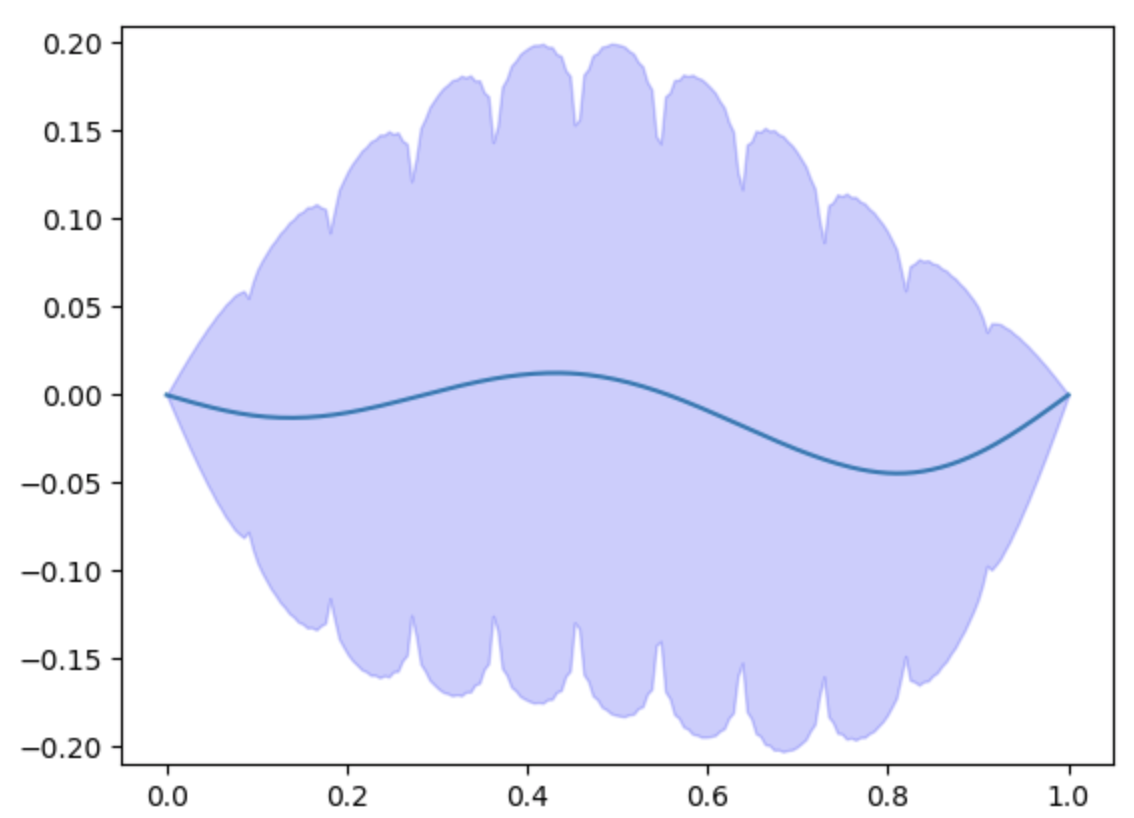}} 
	\caption{$95\%$ intervals of sample paths with 10 design points equispaced on  $(0,1)$.}
    \label{fig:pathswhole}
\end{figure}

\begin{figure}[h]
    \centering
    \subfigure[Nominal prior]{
		\includegraphics[width=0.45\columnwidth]{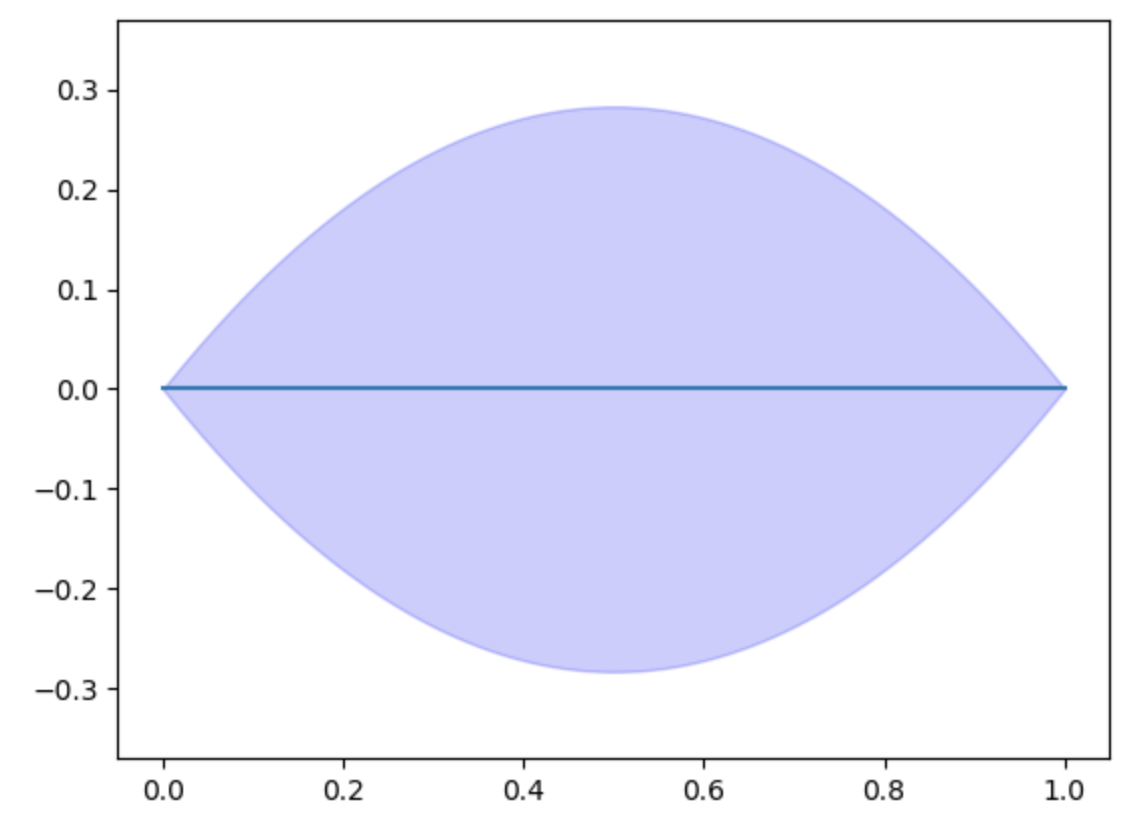}} \hspace{1mm}
	\subfigure[Worst-case prior]{
	\includegraphics[width=0.45\columnwidth]{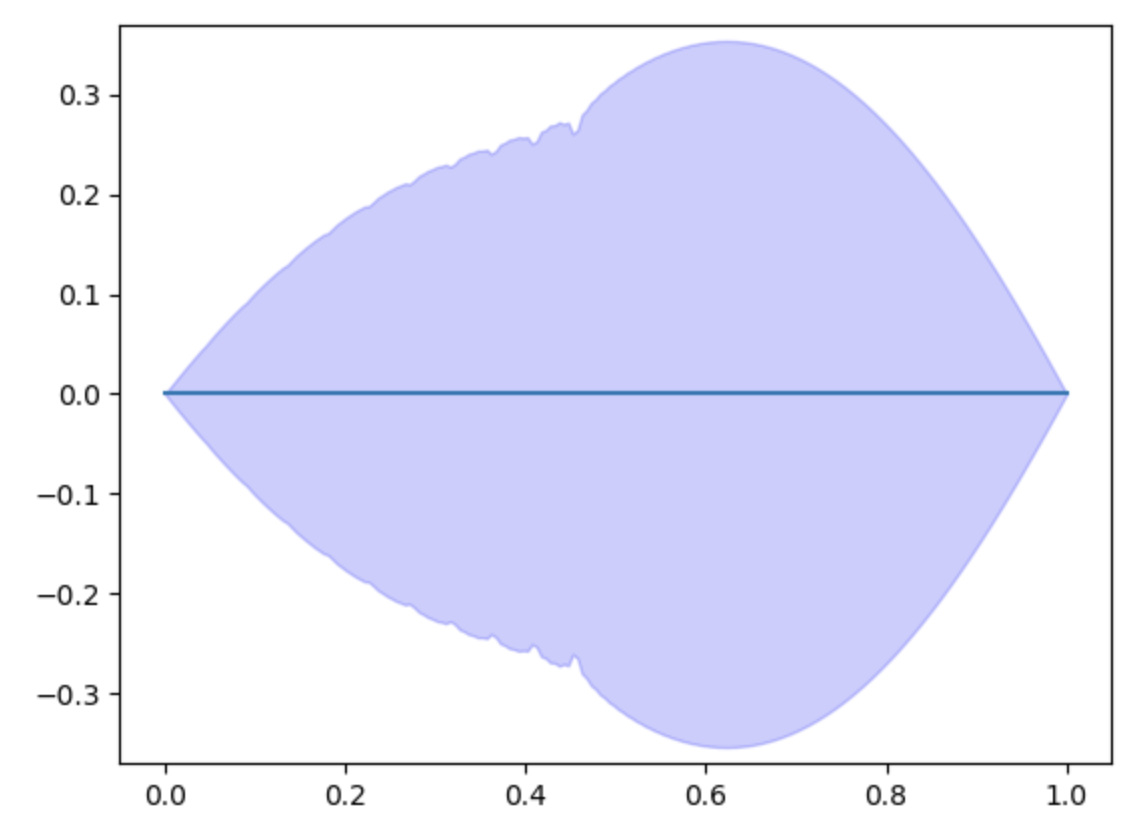}}
	    \hspace{1mm}
	\subfigure[Nominal posterior]{
		\includegraphics[width=0.45\columnwidth]{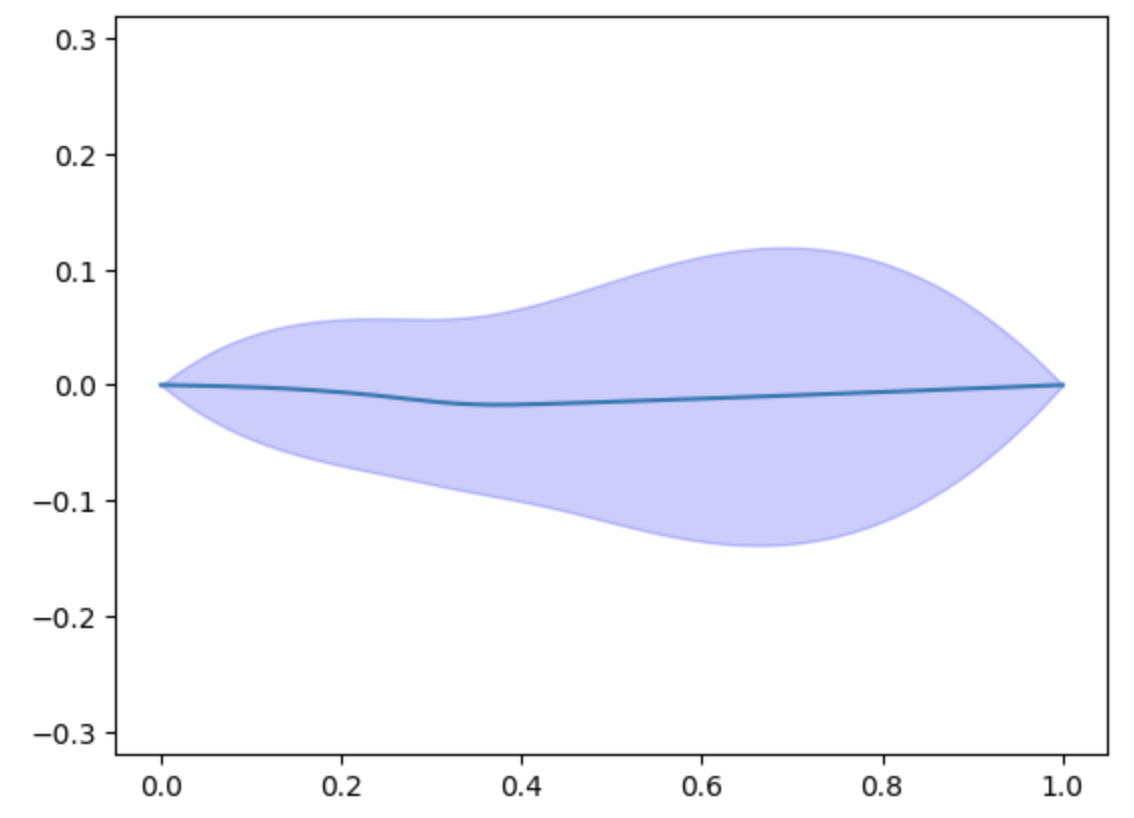}}
	\subfigure[Worst-case posterior]{
		\includegraphics[width=0.45\columnwidth]{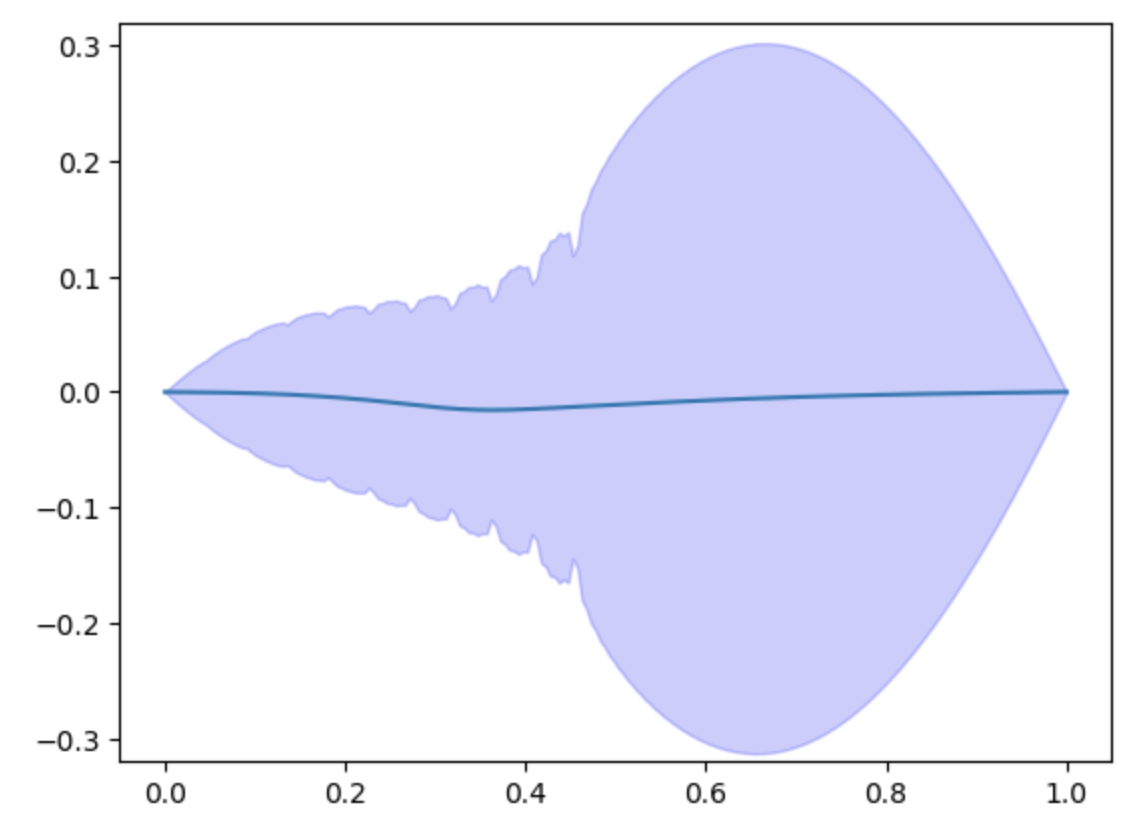}} 
	\caption{$95\%$ intervals of sample paths with 10 design points equispaced on  $(0,0.5)$.}
    \label{fig:pathshalf}
\end{figure}

We first visualize the four correlation functions corresponding to these covariances for our two designs in Figures~\ref{fig:heatmapwhole} and~\ref{fig:heatmaphalf}. The geometry of each design is evident in the worst-case measures. Compared to the nominal measures, we observe ``ripples'' in the worst-case measures corresponding to reductions of correlation between the observed locations. Next, in Figures~\ref{fig:pathswhole} and~\ref{fig:pathshalf}, we plot marginal intervals (at each $x$) containing 95\% of the nominal and the worst-case sample paths. Data to obtain the two posteriors shown here were drawn from the nominal prior measure. In particular, we used the vector of observations 
\begin{verbatim}
  (-0.17,-0.09,0.02,0.04,0.12,0.05,-0.03,0.03,-0.28,-0.15)
\end{verbatim}
for the $(0,1)$ design and
\begin{verbatim}
    (0.03,-0.05,0.08,-0.08,0.15,0.12,-0.25,-0.24,0.16,0.02)
\end{verbatim}
for the $(0,0.5)$ design. Comparing to the nominal prior measures, we observe that the worst-case prior measures have roughly the same overall variance magnitude, but a sharper contrast between the observed and unobserved locations, especially for the $(0,1)$-equispaced designs. On the other hand, comparing to the nominal posterior measures, we observe that the worst-case posterior measures have significantly higher marginal variances in regions away from the observed locations, while the variance increase moderately in regions surrounding the observations. The worst cases are thus perturbed to induce greater uncertainty in regions where information is limited; intuitively, this guarantees greater robustness of the estimates. 

\ymmtd{Figure 3: could we use the same range on the vertical axes for each row of figures? This will better illustrate how the variance changes. Same for Figure 4.}

Next, we vary the baseline parameters to see the effect of the worst-case perturbations on (prior and posterior) sample paths of $b$ compared to sample paths of $b$ under the nominal measures. We focus on the $(0,1)$-equispaced design. For each parameter setting, we draw and visualize five independent sample paths to gauge their qualitative behavior. The posterior sample paths are conditioned on the same observation values as before. As another way of quantifying the impact of the worst-case perturbation, we compute the distance between the prior and posterior covariance matrices on $\PD^{N}$, where these matrices are induced by either the nominal or the worst-case measure. The distance we employ is the natural geodesic distance on the manifold of symmetric positive-definite matrices, also known as the F\"{o}rstner distance~\cite{ref:Forstner2003metric} and (up to a constant) Rao's distance \cite{atkinson1981rao,rao1987differential}. This distance is invariant under affine transformations and under inversion and has been used extensively to compare covariance matrices in previous work \cite{ref:Forstner2003metric,spantini2015optimal,spantini2017goal}.

\begin{figure}[h]
    \centering
    \subfigure[Nominal prior]{
	 \includegraphics[width=0.45\columnwidth]{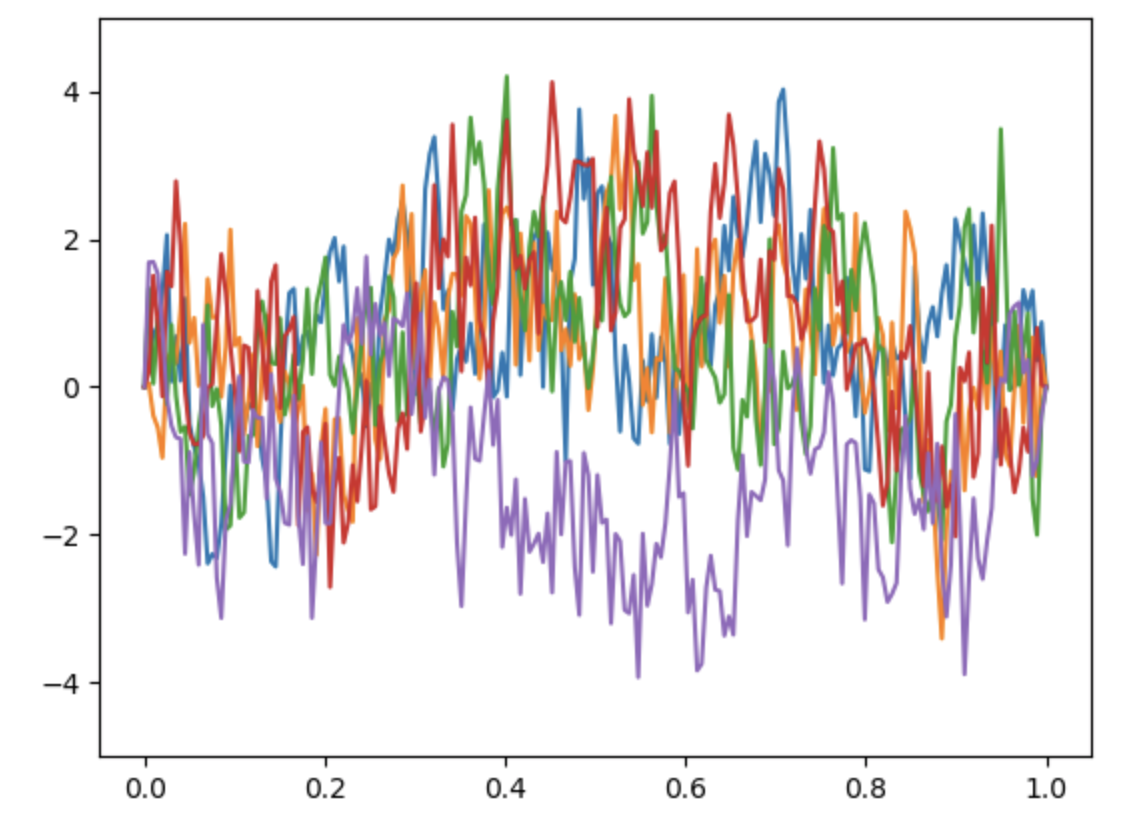}} \hspace{1mm}
	\subfigure[Worst-case prior]{
	\includegraphics[width=0.45\columnwidth]{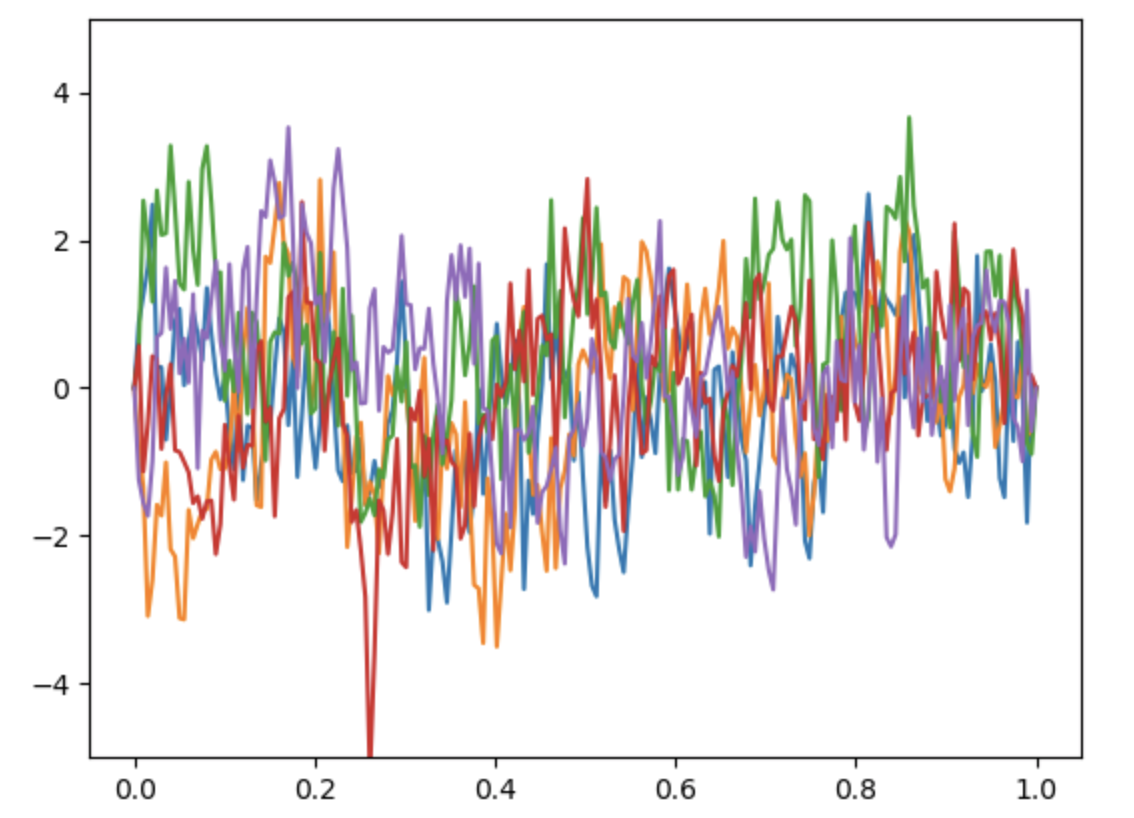}}
	    \hspace{1mm}
	\subfigure[Nominal posterior]{
		\includegraphics[width=0.45\columnwidth]{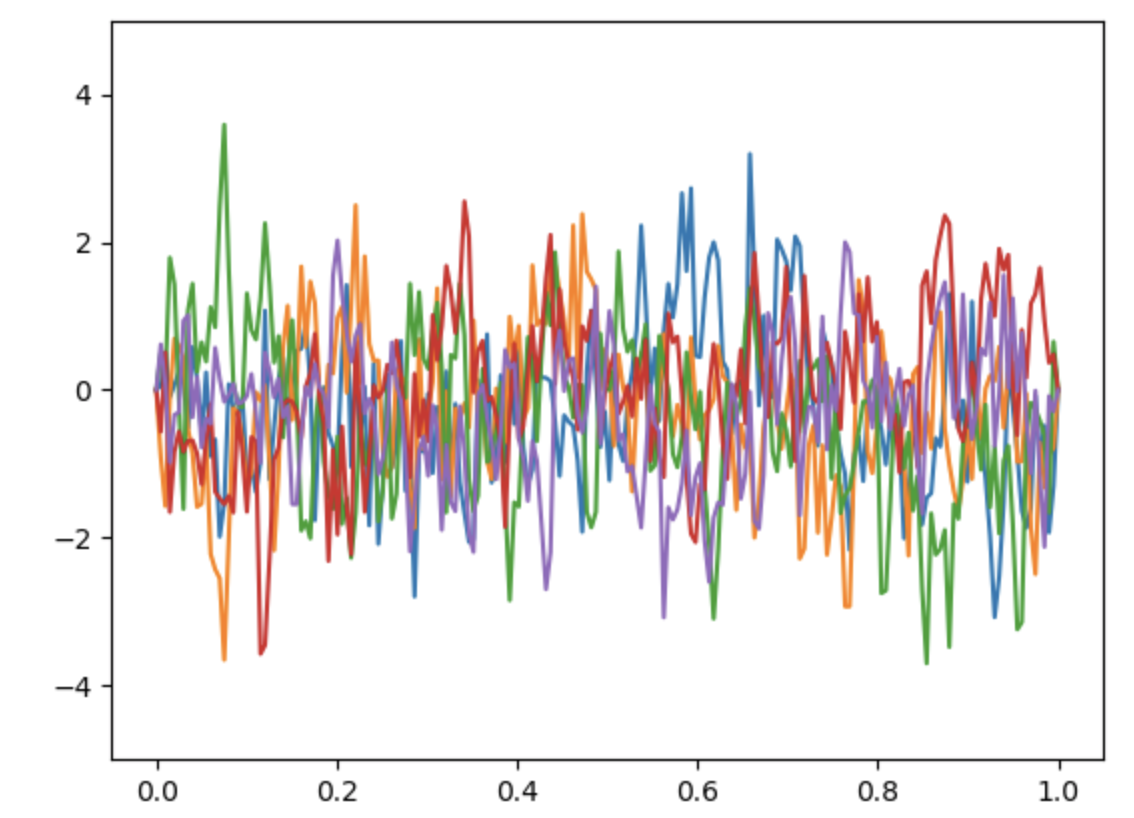}}
	\subfigure[Worst-case posterior]{
		\includegraphics[width=0.45\columnwidth]{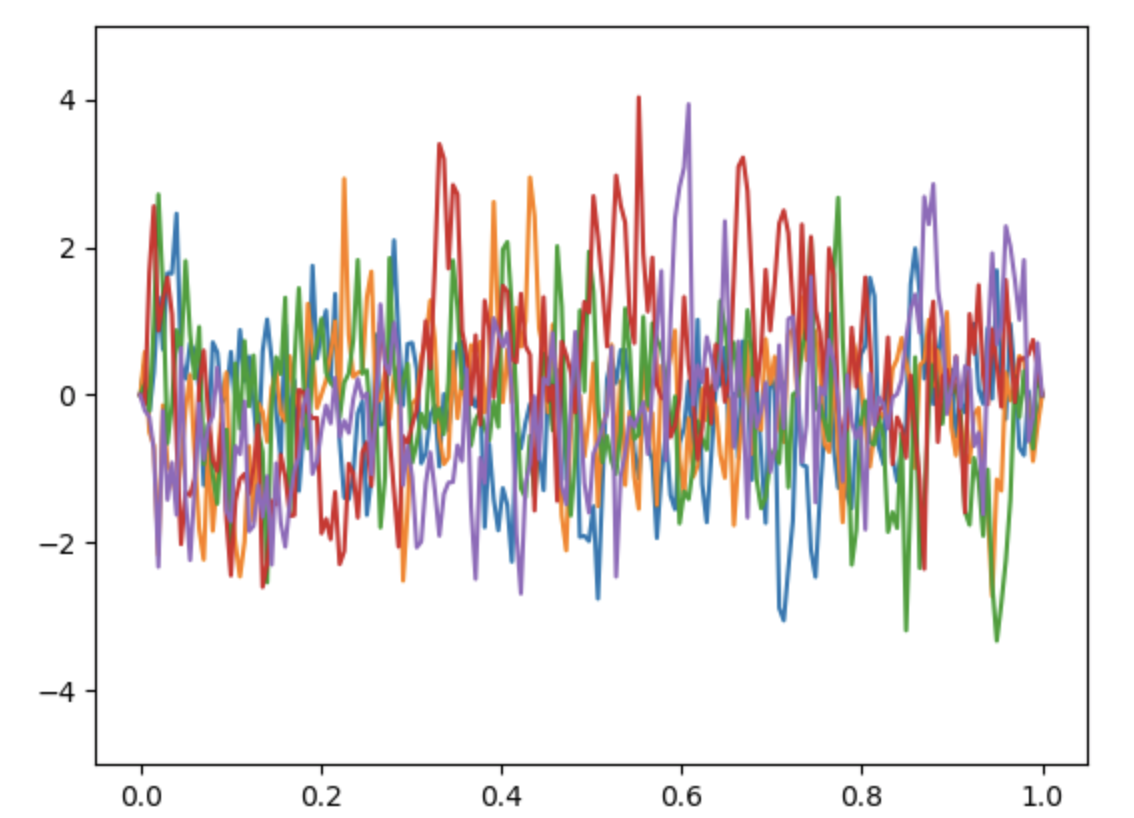}} 
	\caption{Sample paths with $\alpha=0.51$ and 10 designs equispaced in $(0,1)$.}
    \label{fig:alpha51pathswhole}
\end{figure}

\begin{figure}[h]
    \centering
    \subfigure[Nominal prior]{
		\includegraphics[width=0.45\columnwidth]{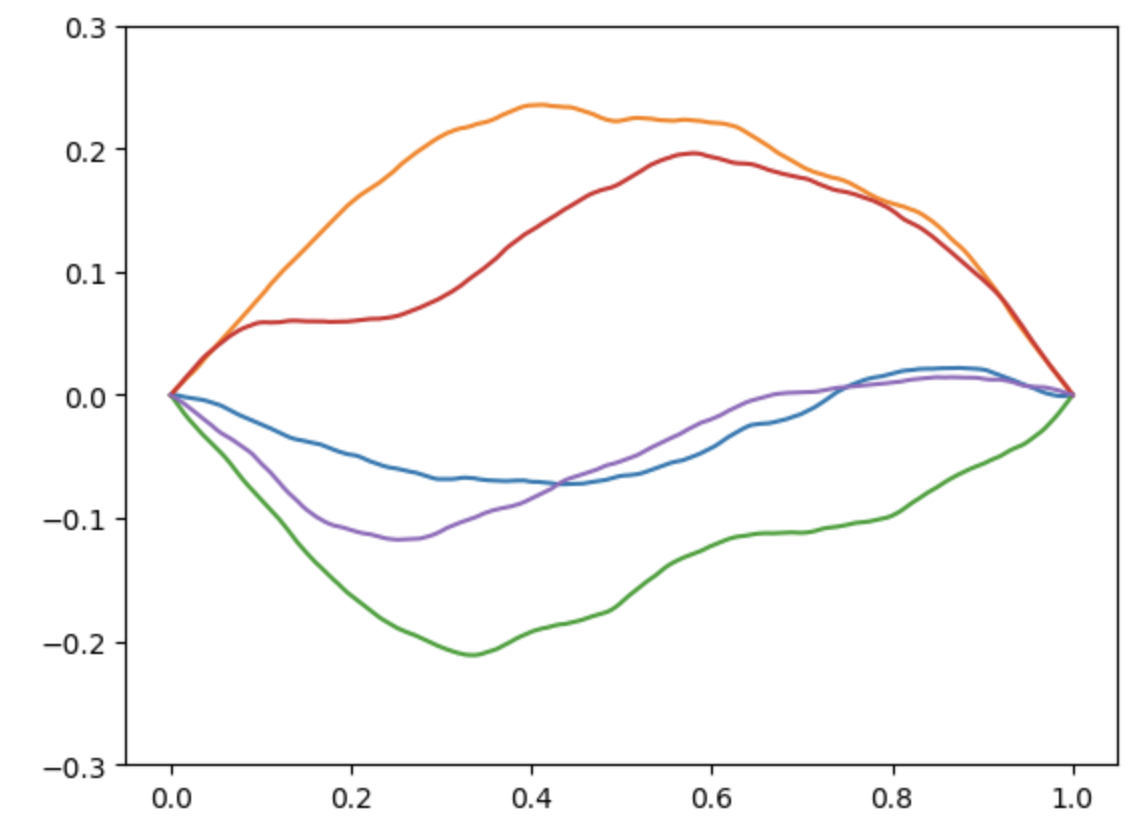}} \hspace{1mm}
	\subfigure[Worst-case prior]{
	\includegraphics[width=0.45\columnwidth]{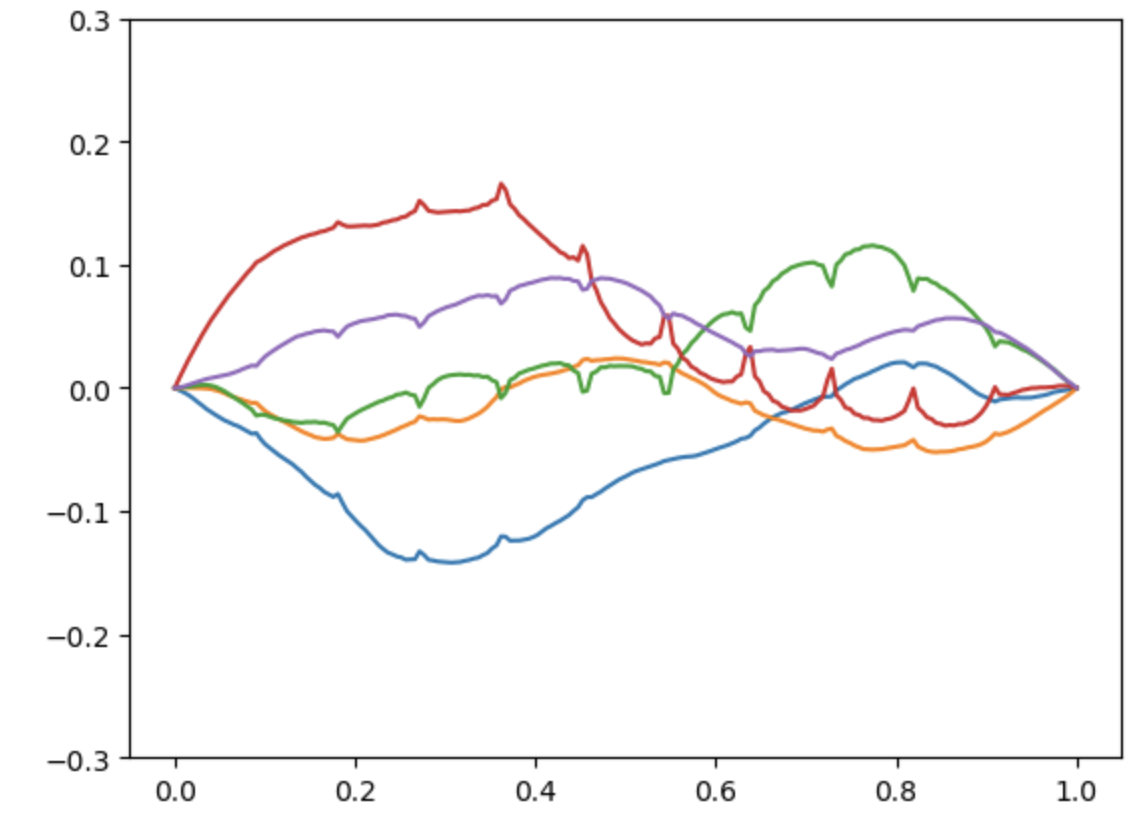}}
	    \hspace{1mm}
	\subfigure[Nominal posterior]{
		\includegraphics[width=0.45\columnwidth]{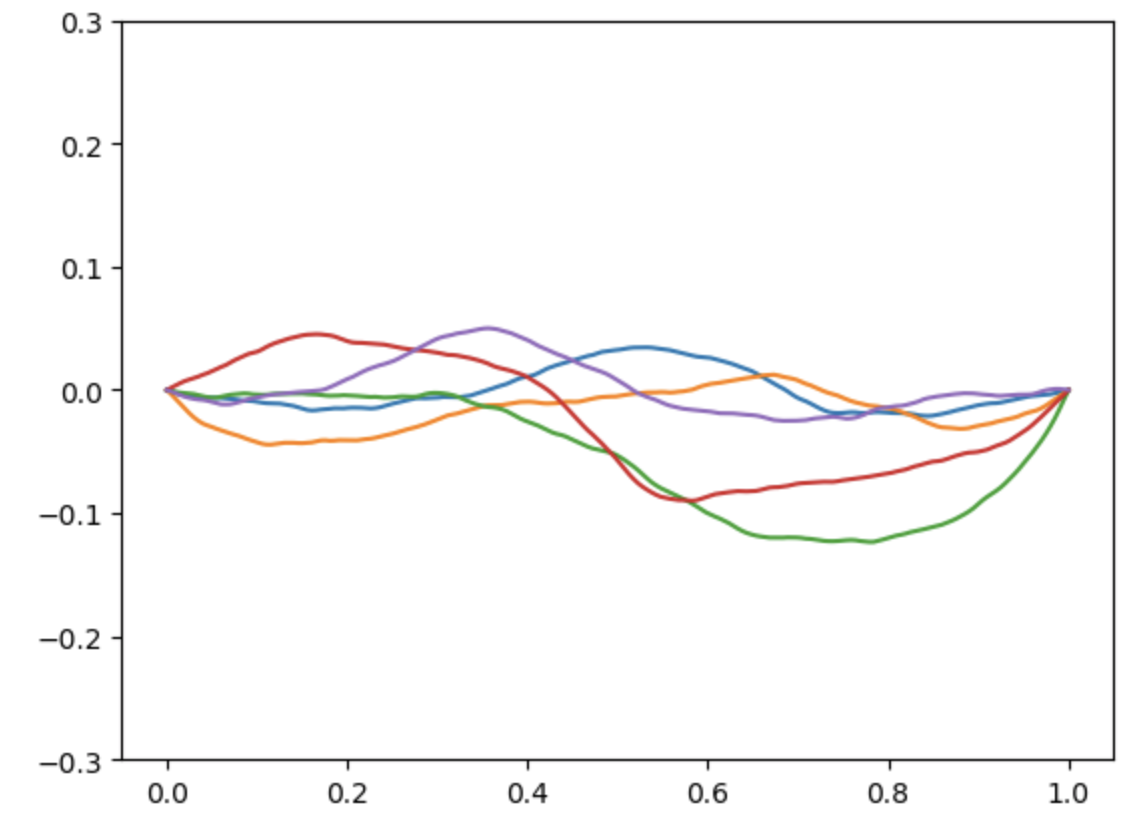}}
	\subfigure[Worst-case posterior]{
		\includegraphics[width=0.45\columnwidth]{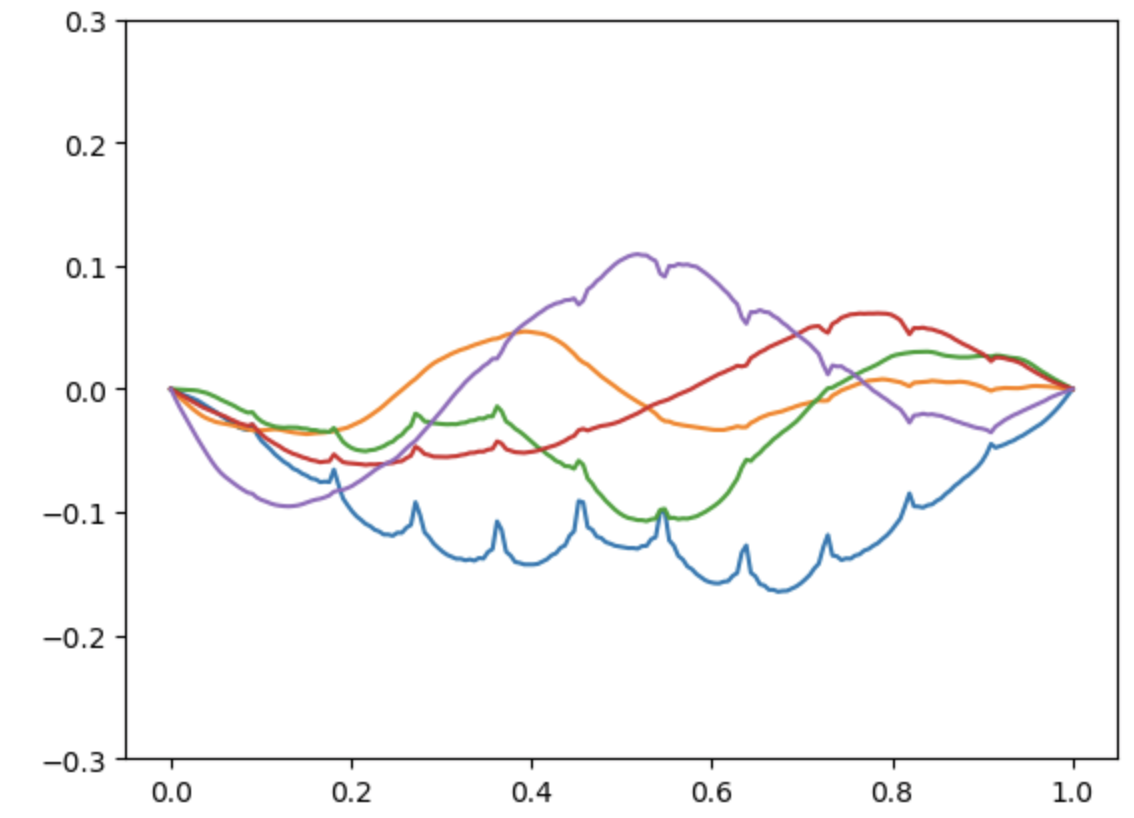}} 
	\caption{Sample paths with $\alpha=2$ and 10 designs equispaced in $(0,1)$.}
    \label{fig:alpha2pathswhole}
\end{figure}

\begin{figure}[h]
    \centering
    \subfigure[Nominal prior]{
		\includegraphics[width=0.45\columnwidth]{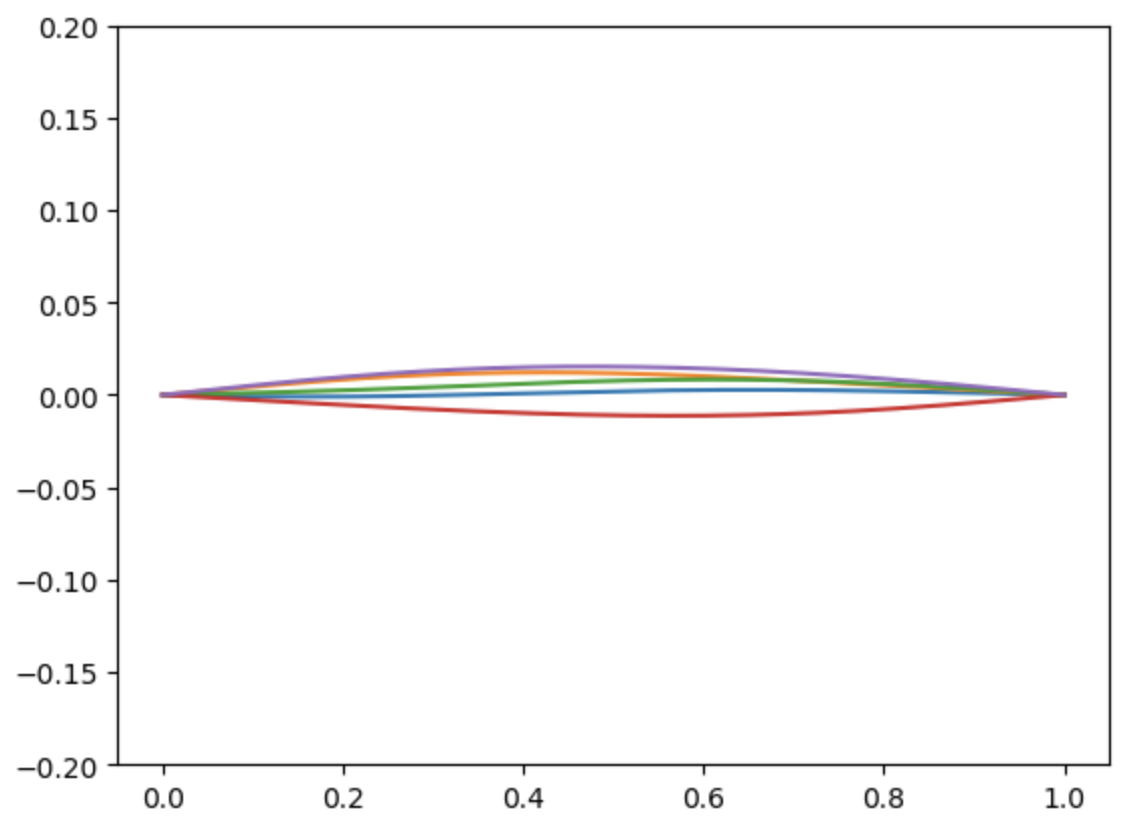}} \hspace{1mm}
	\subfigure[Worst-case prior]{
	\includegraphics[width=0.45\columnwidth]{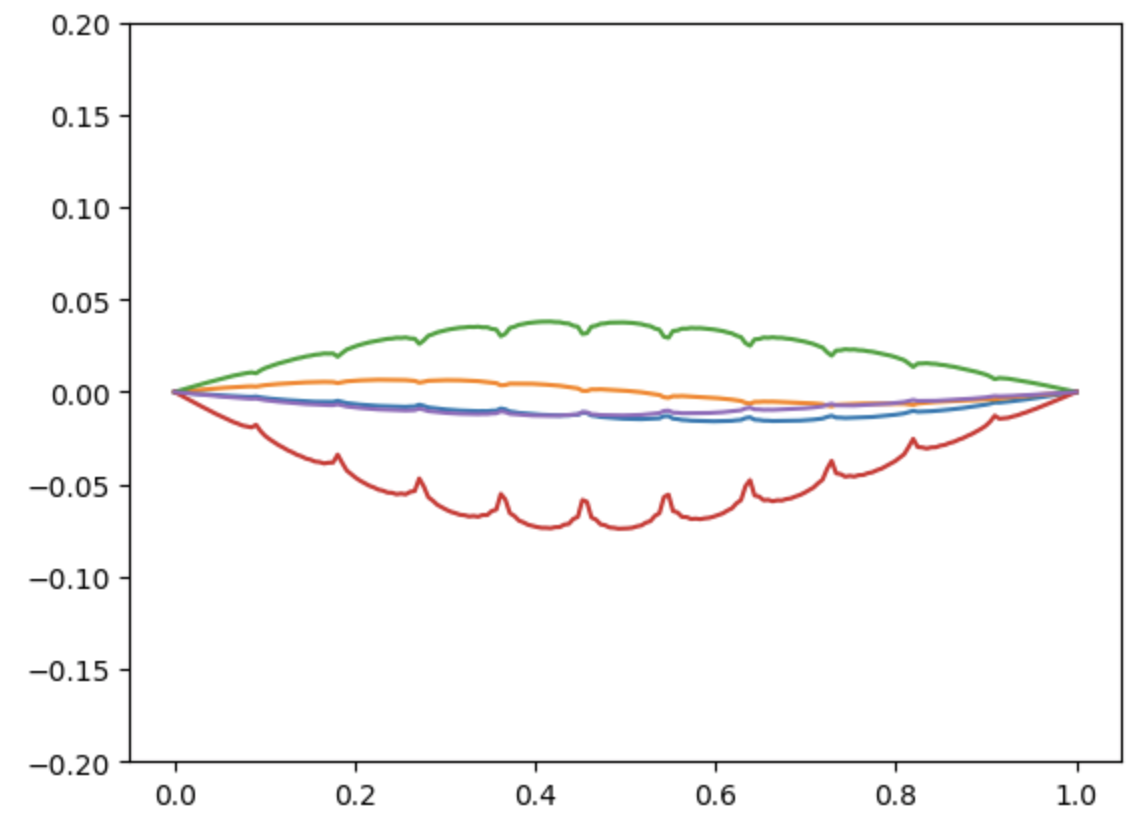}}
	    \hspace{1mm}
	\subfigure[Nominal posterior]{
		\includegraphics[width=0.45\columnwidth]{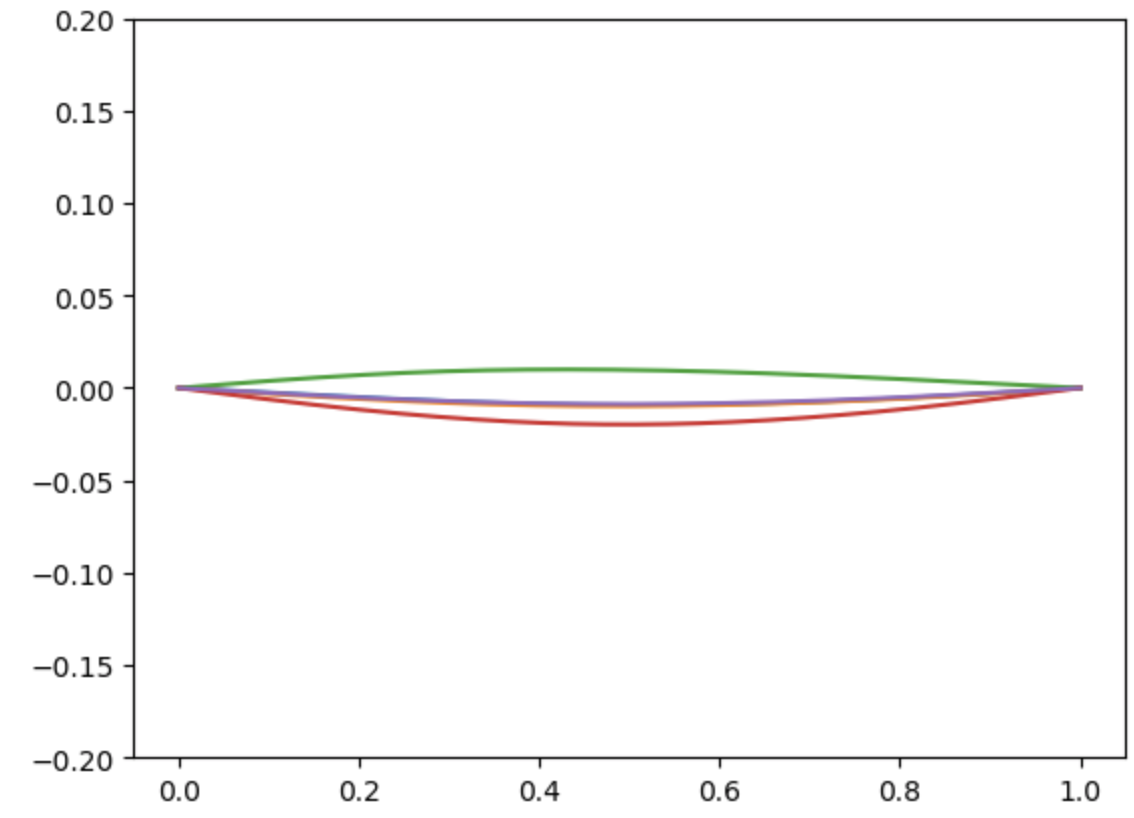}}
	\subfigure[Worst-case posterior]{
		\includegraphics[width=0.45\columnwidth]{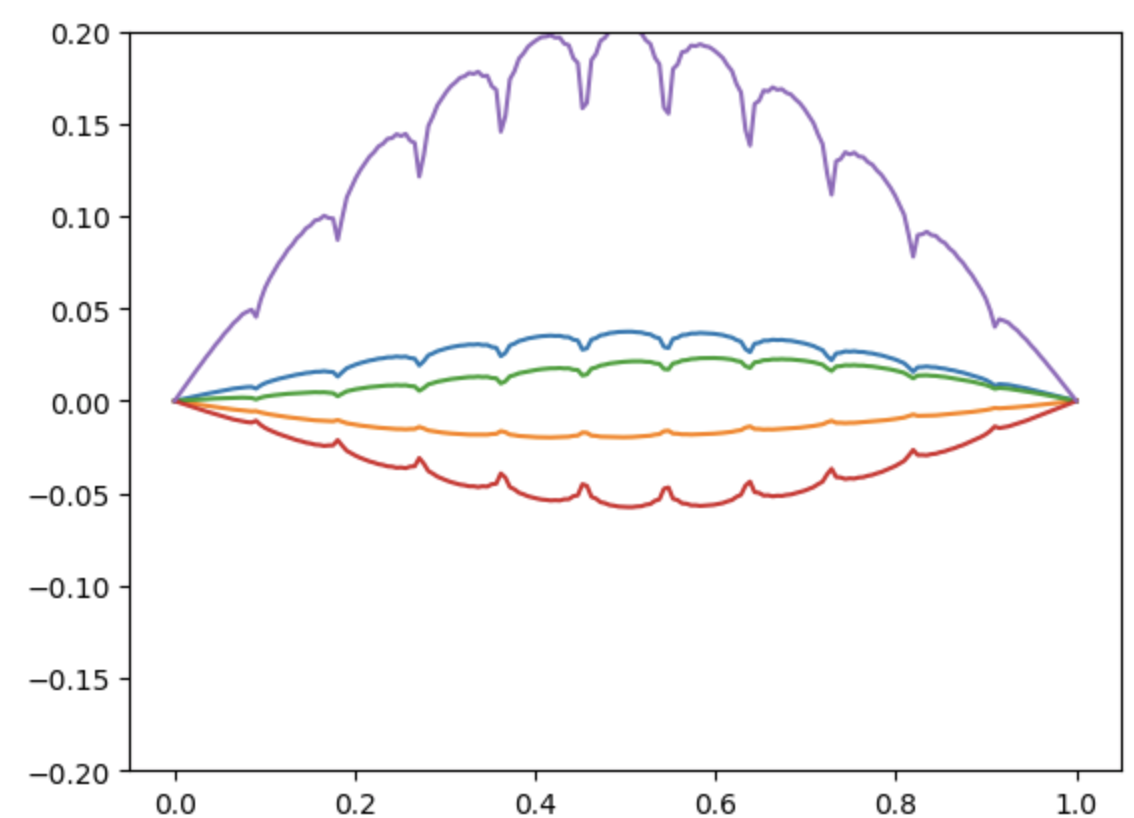}} 
	\caption{Sample paths with $\alpha=4$ and 10 designs equispaced in $(0,1)$.}
    \label{fig:alpha4pathswhole}
\end{figure}

\begin{figure}[h]
    \centering
    \subfigure[Prior, $\beta=0.7$]{
		\includegraphics[width=0.45\columnwidth]{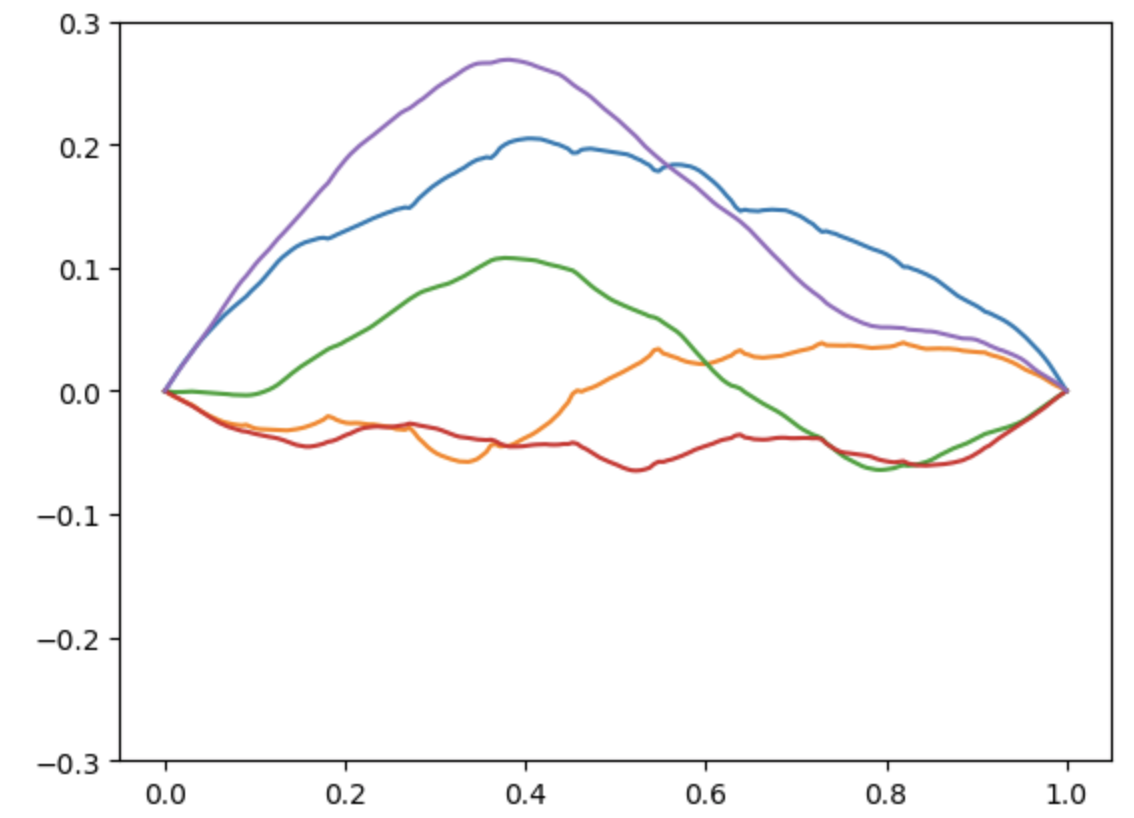}} \hspace{1mm}
	\subfigure[Posterior, $\beta=0.7$]{
	\includegraphics[width=0.45\columnwidth]{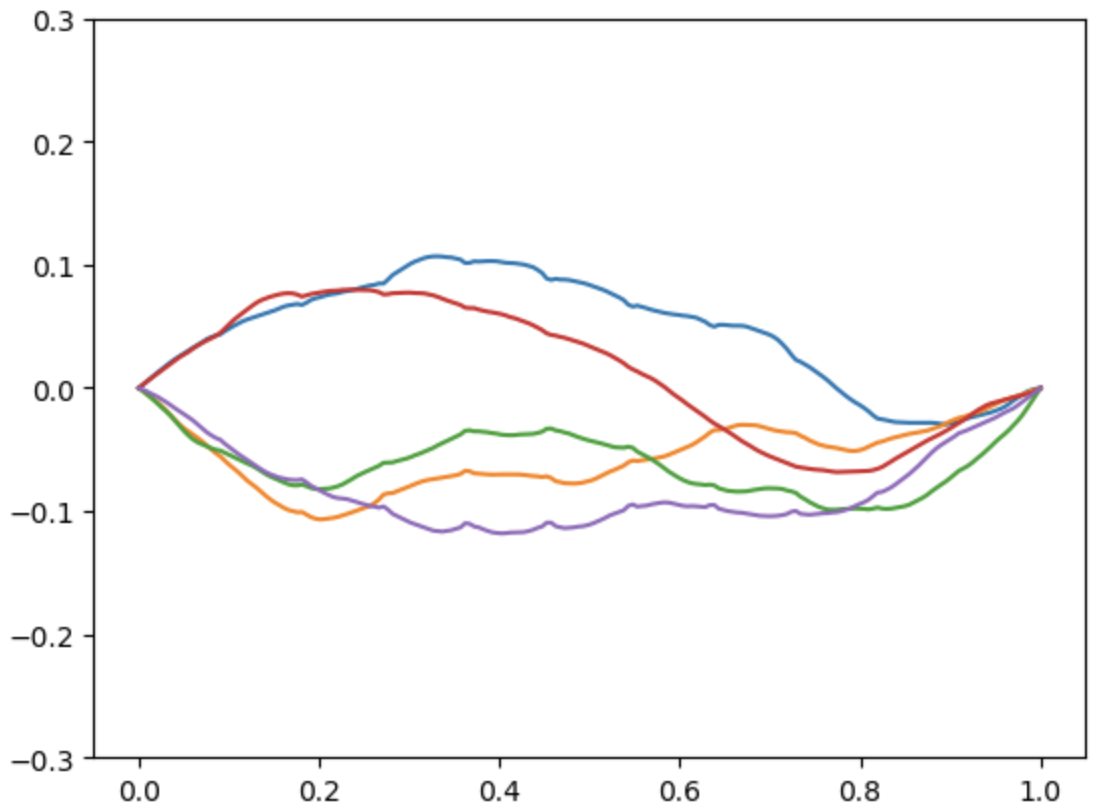}}
	    \hspace{1mm}
	\subfigure[Prior, $\beta=1$]{
		\includegraphics[width=0.45\columnwidth]{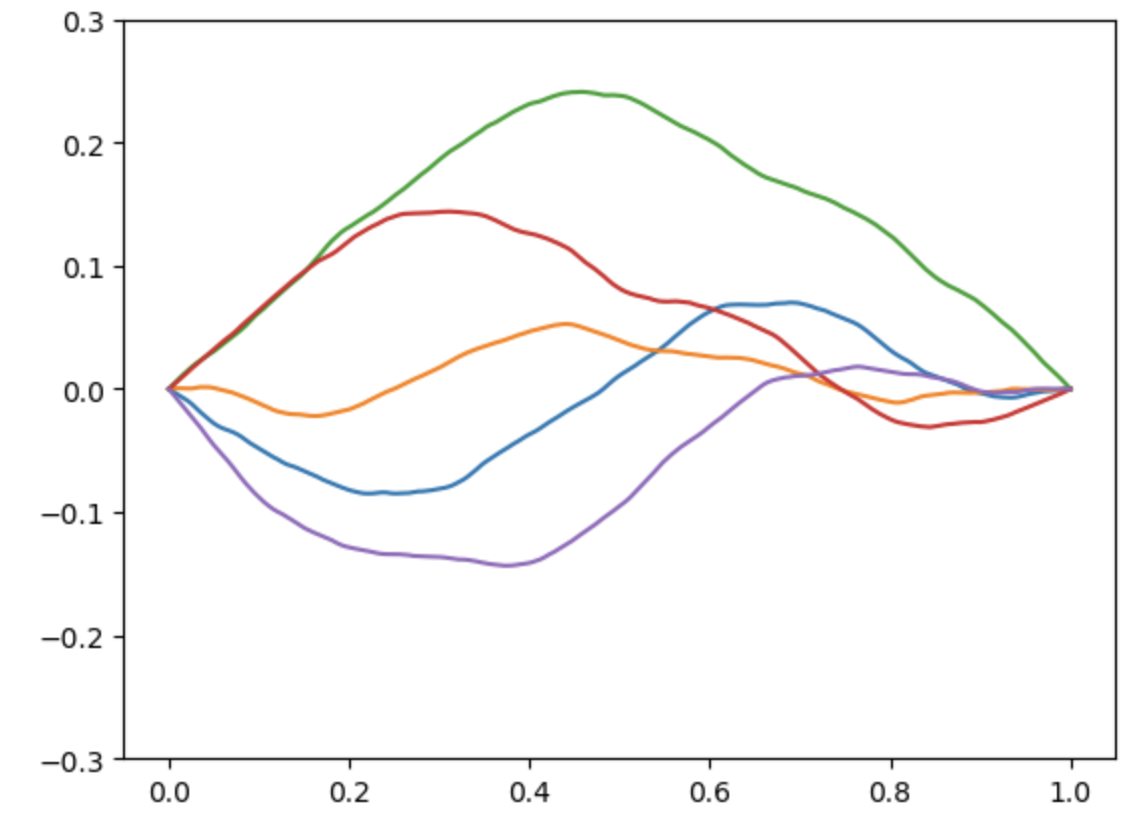}}
	\subfigure[Posterior, $\beta=1$]{
		\includegraphics[width=0.45\columnwidth]{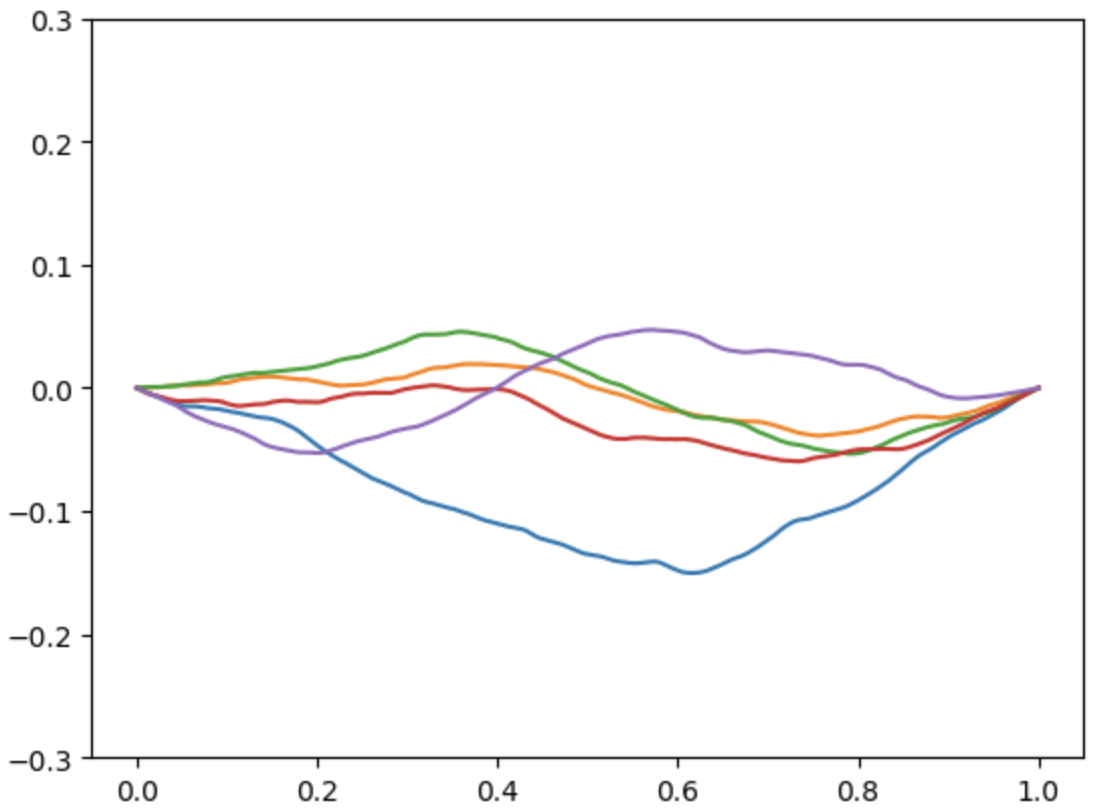}}
	\caption{Worst-case sample paths with varying $\beta$ and 10 designs equispaced in $(0,1)$.}
    \label{fig:betapathswhole}
\end{figure}

\begin{figure}[h]
    \centering
    \subfigure[Prior, $\delta=0.01$]{
		\includegraphics[width=0.45\columnwidth]{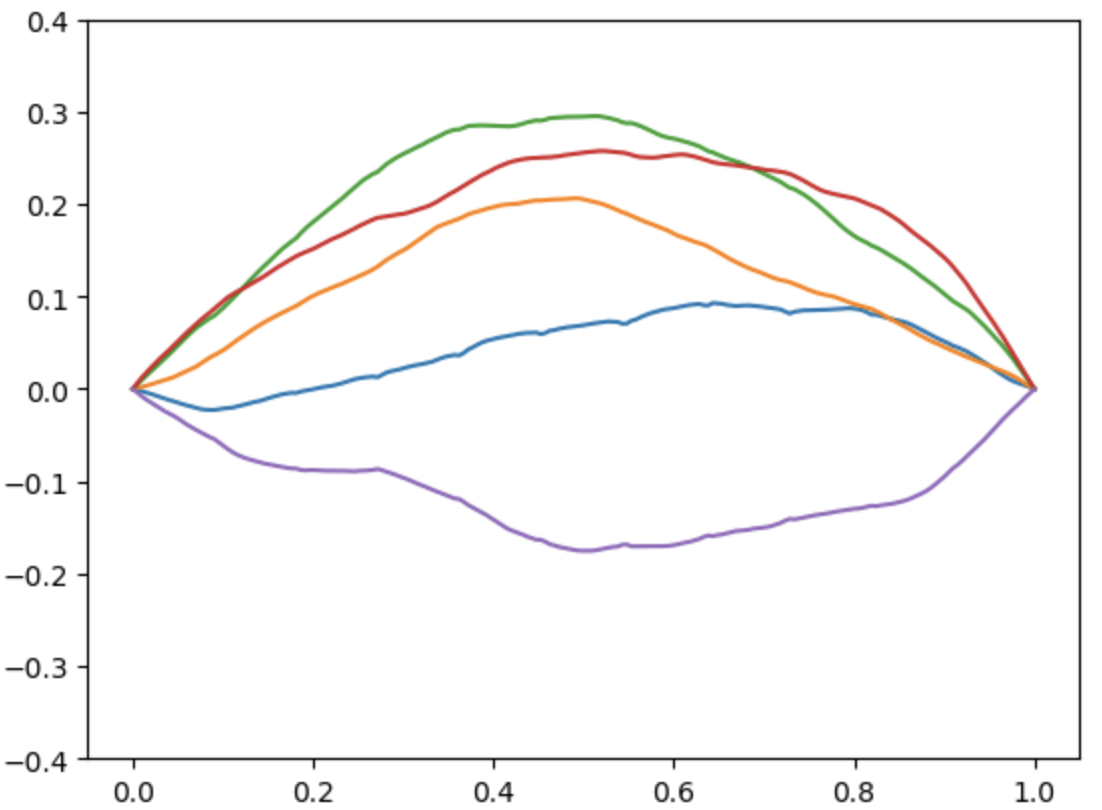}} \hspace{1mm}
	\subfigure[Posterior, $\delta=0.01$]{
	\includegraphics[width=0.45\columnwidth]{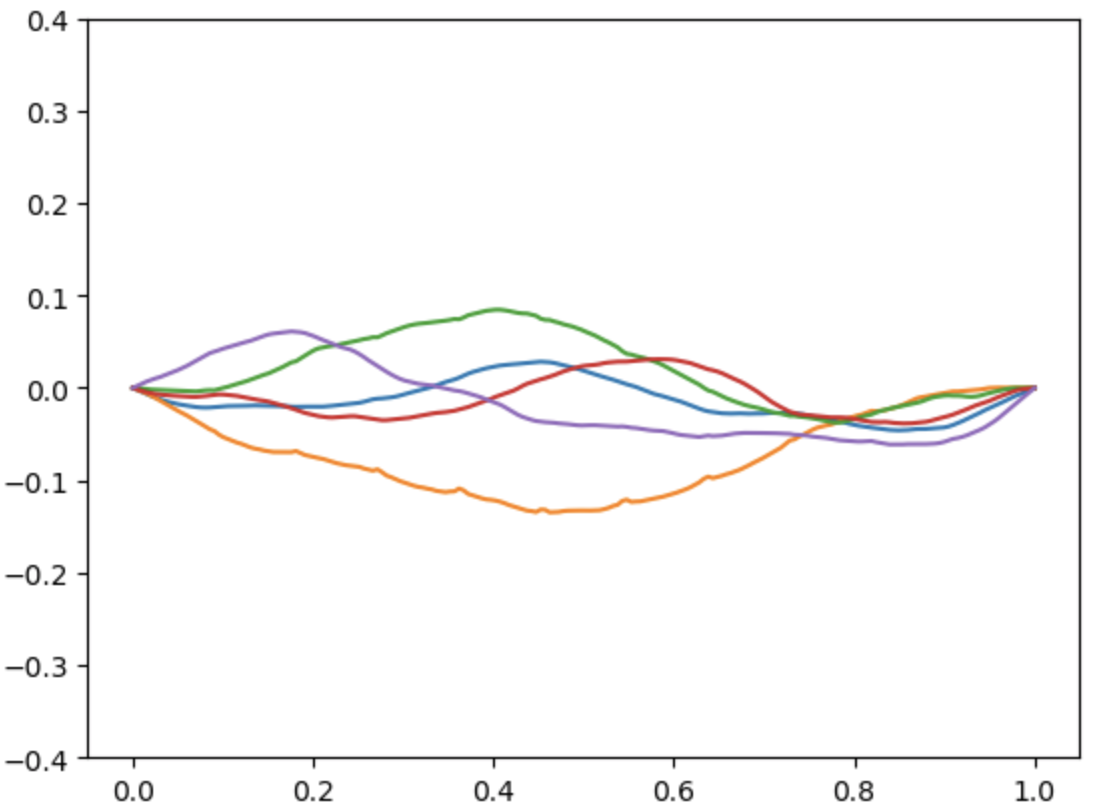}}
	    \hspace{1mm}
	\subfigure[Prior, $\delta=1$]{
		\includegraphics[width=0.45\columnwidth]{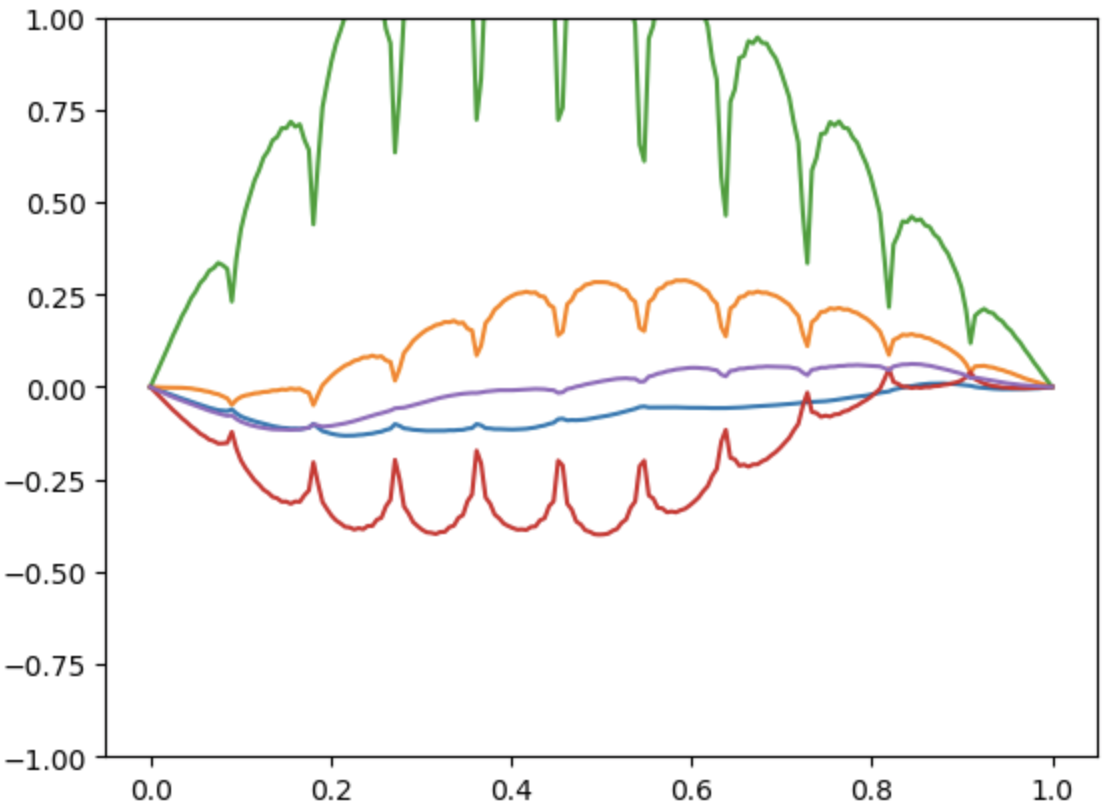}}
	\subfigure[Posterior, $\delta=1$]{
		\includegraphics[width=0.45\columnwidth]{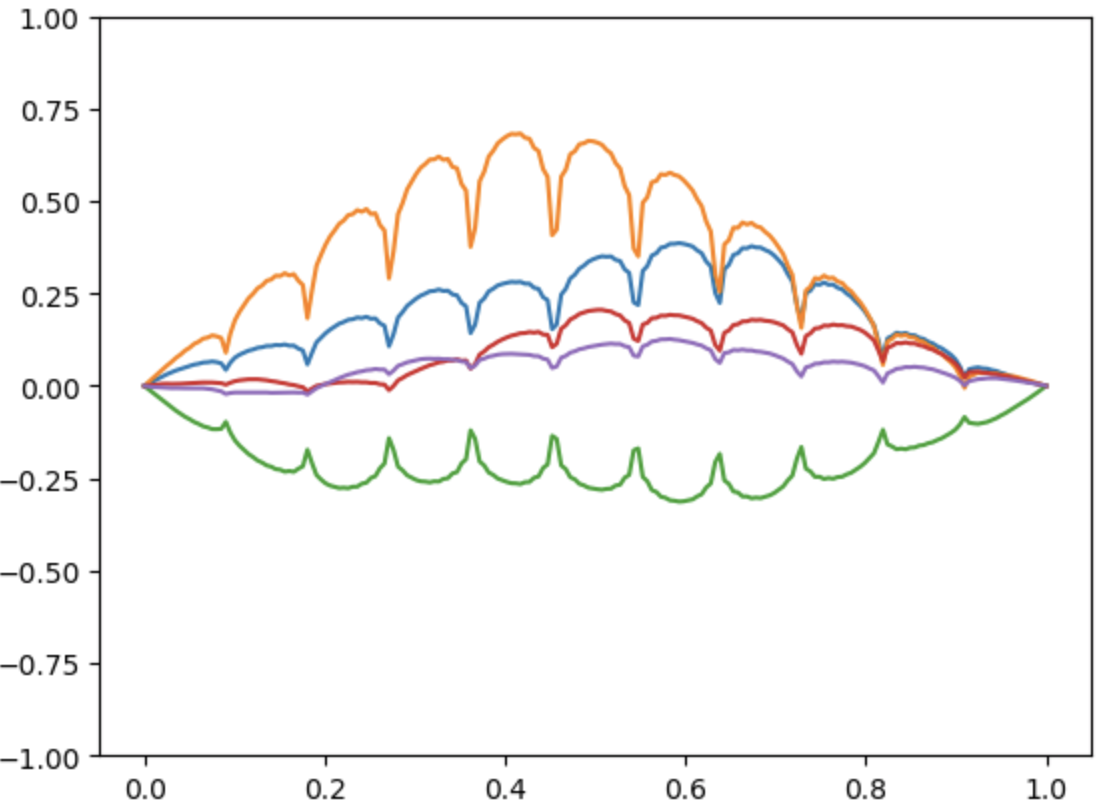}}
	\caption{Worst-case sample paths with varying $\delta$ and 10 designs equispaced in $(0,1)$.}
    \label{fig:deltapathswhole}
\end{figure}

\begin{figure}[h]
    \centering
    \subfigure[Nominal prior]{
		\includegraphics[width=0.45\columnwidth]{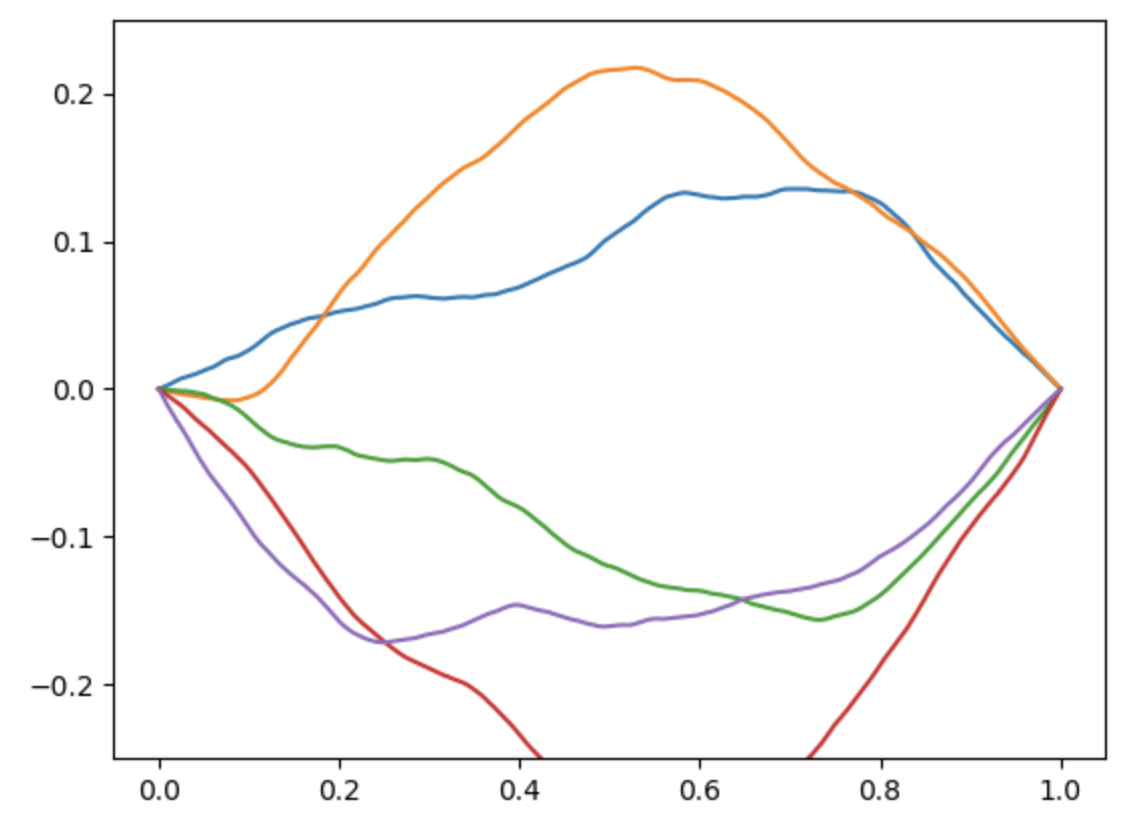}} \hspace{1mm}
	\subfigure[Worst-case prior]{
	\includegraphics[width=0.45\columnwidth]{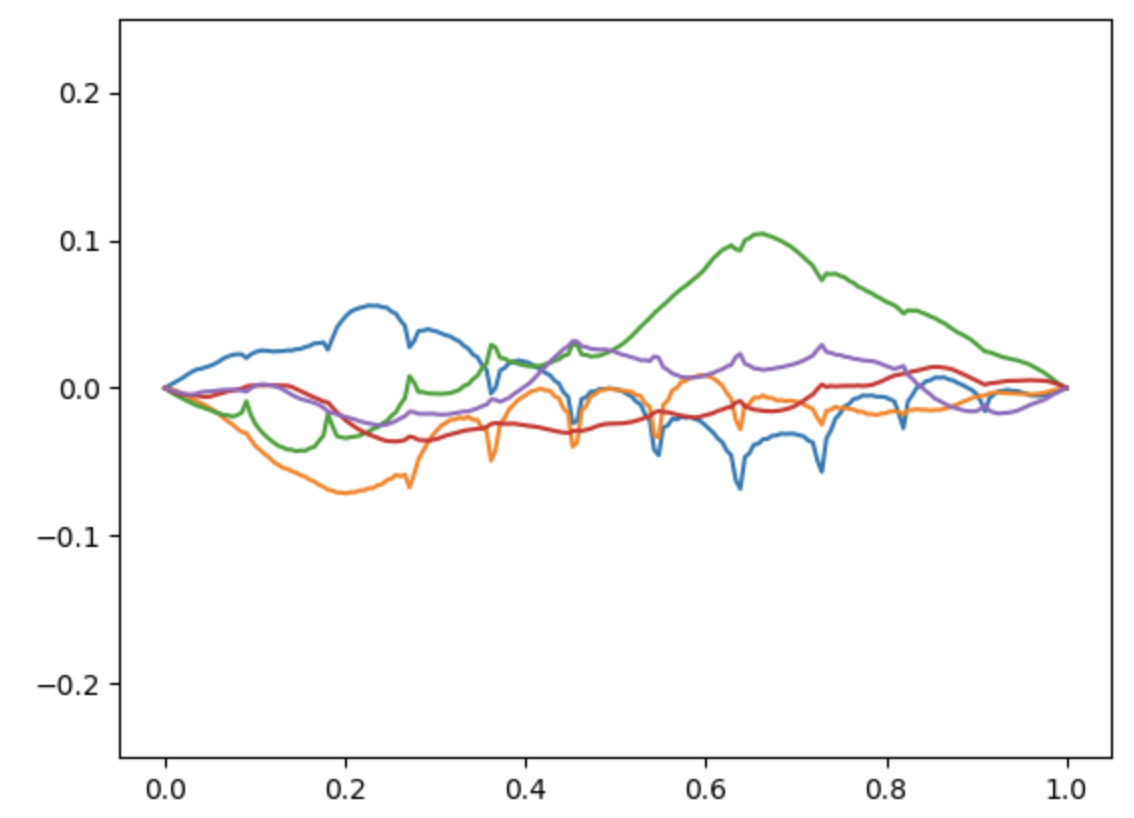}}
	    \hspace{1mm}
	\subfigure[Nominal posterior]{
		\includegraphics[width=0.45\columnwidth]{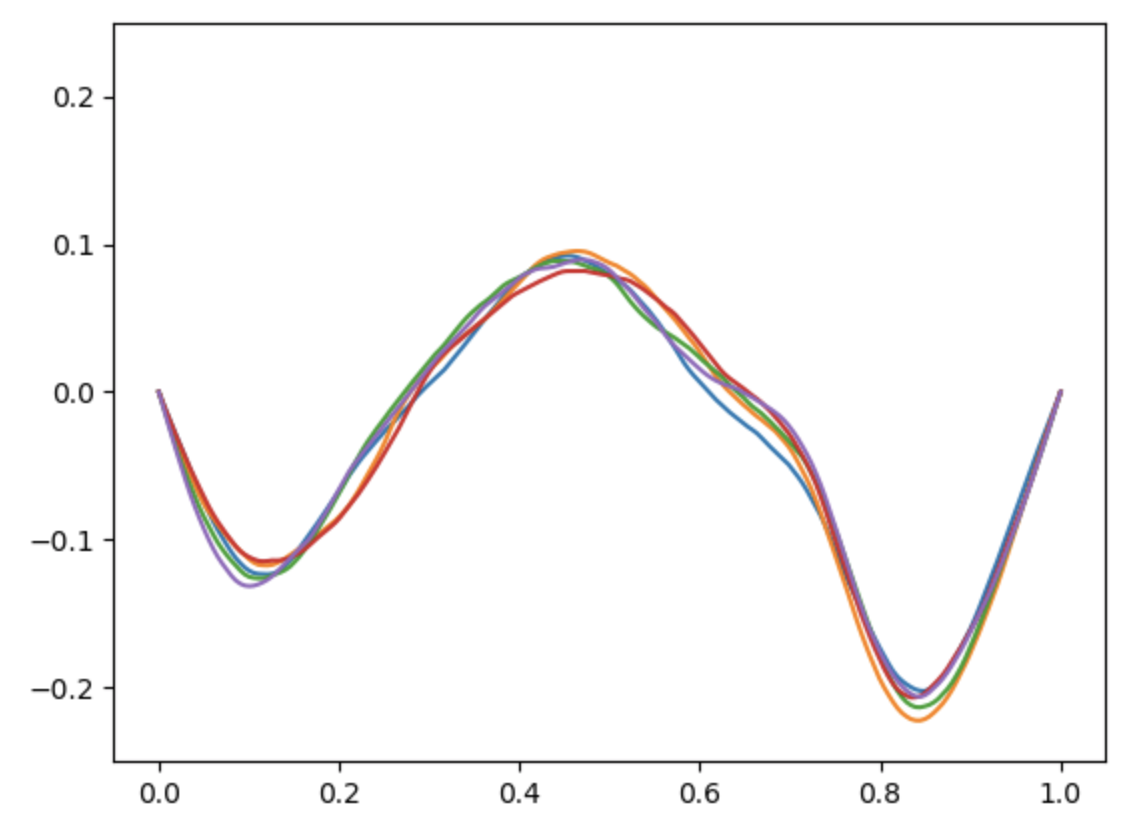}}
	\subfigure[Worst-case posterior]{
		\includegraphics[width=0.45\columnwidth]{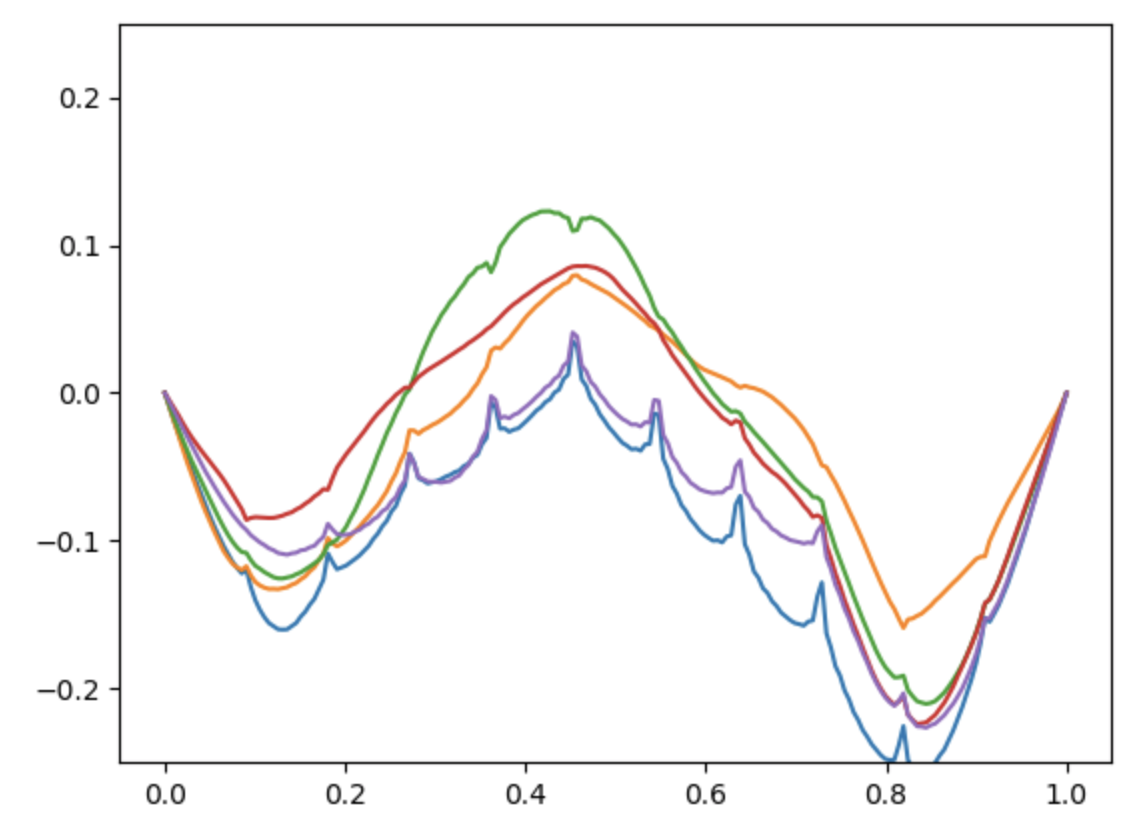}} 
	\caption{Sample paths with $\sigma=0.01$ and 10 designs equispaced in $(0,1)$.}
    \label{fig:sigma01pathswhole}
\end{figure}

\begin{figure}[h]
    \centering
    \subfigure[Nominal prior]{
		\includegraphics[width=0.45\columnwidth]{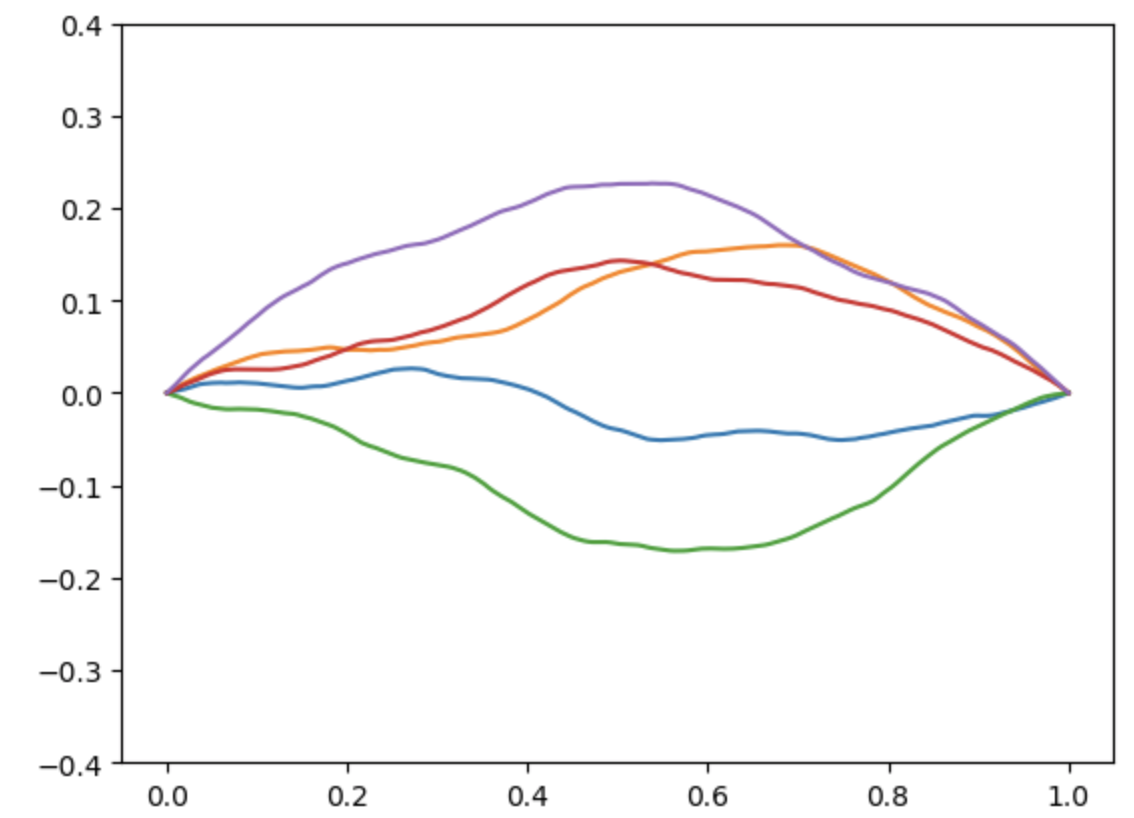}} \hspace{1mm}
	\subfigure[Worst-case prior]{
	\includegraphics[width=0.45\columnwidth]{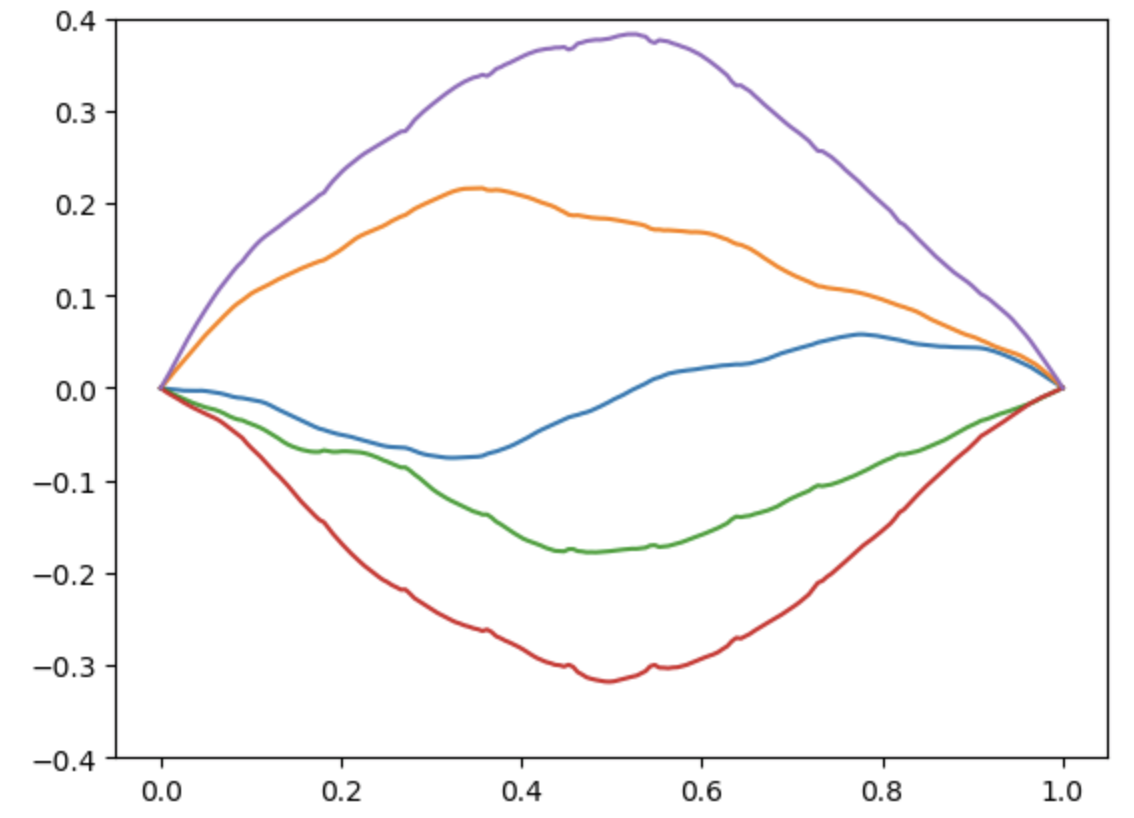}}
	    \hspace{1mm}
	\subfigure[Nominal posterior]{
		\includegraphics[width=0.45\columnwidth]{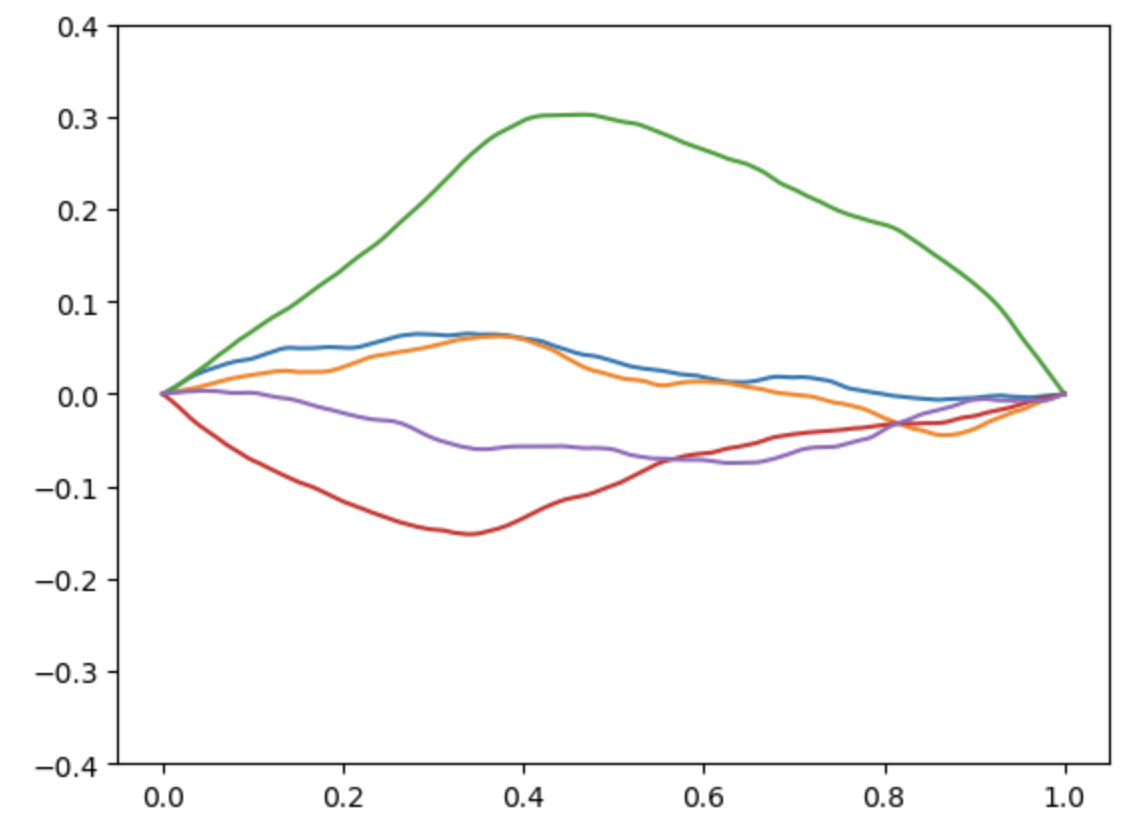}}
	\subfigure[Worst-case posterior]{
		\includegraphics[width=0.45\columnwidth]{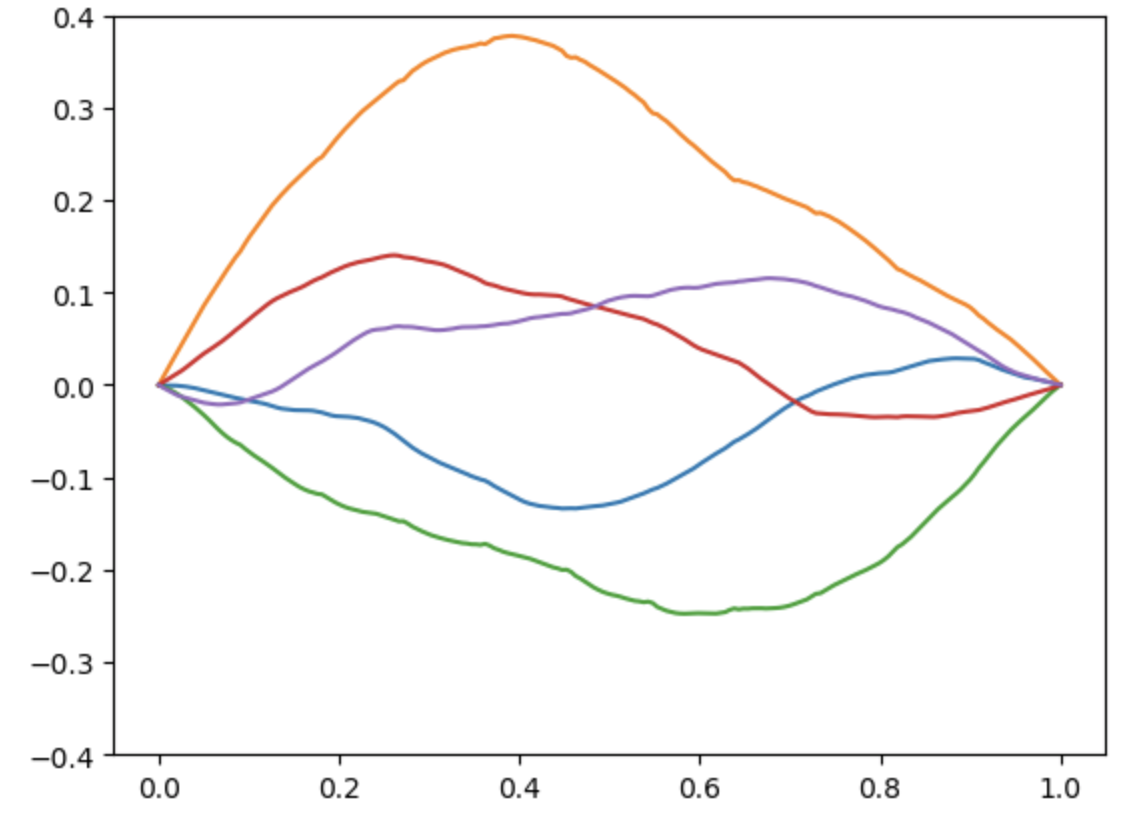}} 
	\caption{Sample paths with $\sigma=1$ and 10 designs equispaced in $(0,1)$.}
    \label{fig:sigma1pathswhole}
\end{figure}

With the remaining parameters fixed to the base case, we explore the following parameter variations.
\begin{enumerate}
    \item Prior smoothness: choose $\alpha\in\{0.51,2,4\}$. Results are shown in Figures~\ref{fig:alpha51pathswhole}--\ref{fig:alpha4pathswhole}; note that we include the baseline value $\alpha=2$ for comparison.
    \item Adversarial perturbation smoothness: choose $\beta\in\{0.7,1\}$. Results are shown in Figure~\ref{fig:betapathswhole}. Note that the nominal prior and posterior are the same as in the baseline setting.
    \item Size of the Wasserstein ambiguity set: choose $\delta^2\in\{0.01,1\}$. Results are shown in Figure~\ref{fig:deltapathswhole}. Note that the nominal prior and posterior are the same as in the baseline setting.
    \item Magnitude of the observation noise: choose $\sigma\in\{0.01,1\}$. Results are shown in Figures~\ref{fig:sigma01pathswhole} and~\ref{fig:sigma1pathswhole}.
\end{enumerate}

The corresponding nominal and worst-case prior-to-posterior distances are reported Table~\ref{tab:forstnerwhole}. Combining these qualitative and quantitative results, we observe that in cases where there is: a larger $\alpha$ (i.e., a smoother prior); a smaller $\beta$ (i.e., smaller penalty on modes that induce roughness); a larger $\delta$ (i.e., wider range of admissible perturbations); or a smaller $\sigma$ (i.e., smaller observation noise), the worst-case distributions induce sharper contrasts between the observed and unobserved locations in both the prior and posterior sample paths.
\begin{table}[h!]
\centering
\begin{tabular}{||c|c|c|c|c|c|c|c|c|c||} 
\hline
 &baseline&$\alpha=0.51$&$\alpha=4$&$\beta=0.7$&$\beta=1$&$\delta^2=0.01$&$\delta^2=1$&$\sigma=0.01$&$\sigma=1$\\
 \hline
Nominal&2.56&16.24&0.11&2.56&2.56&2.56&2.56&8.92&0.11 \\
Worst-case&12.82&16.15&5.33&8.97&2.13&9.88&10.90&19.74&3.64\\
 \hline
\end{tabular}
\caption{Natural distance between the prior and posterior covariance matrices under different problem settings, for 10 design points equispaced on $(0,1)$.}
\label{tab:forstnerwhole}
\end{table}

\begin{appendix}
\section{Proofs of Theorems and Additional Details}

\subsection{Proof of Theorem~\ref{thm:swap}}

\begin{proof}[Proof of the strong duality]
We start with a few definitions used in the proof. Recall that $\mathcal M= \{\phi: \R^m\to C(\mathcal D),~\phi \text{ is measurable}\}$ is the set of estimators for $b$.
\begin{definition}[Affine predictor]\label{def:affinepredictor_revision}
For $\phi(Y_1,\ldots,Y_m)$ a mapping from $\mathbb{R}^m$ (data space) to $C(\mathcal{D})$ (parameter space), we say that $\phi$ is affine if it is of the form
\[
\phi(Y_1,\ldots,Y_m) = \alpha_0 + \sum_{j=1}^m \alpha_j Y_j,
\]
for some functions $\alpha_j\in C(\mathcal{D})$\label{nomencl:alphaj}\nomenclature{$\alpha_j$}{coefficient functions in affine predictor}{nomencl:alphaj} , $0\leq j\leq m$. We write $\mathcal M_{\text{aff}}= \{\phi:\R^m\to C(\mathcal D),~ \phi \text{~is affine}\} \subset \mathcal M$.\label{nomencl:Maff}\nomenclature{$\mathcal M_{\text{aff}}$}{affine subspace of $\mathcal M$}{nomencl:Maff}
\end{definition}

Let $\Span\{e_n\}_{n=1}^N$ and $\Span\{T(e_n)\}_{n=1}^N$ denote the closed subspaces of $C(\mathcal D)$ spanned by the first $N$ basis functions from Assumption~\ref{assmp:fullrank}, and by the action of $T$ on these basis, respectively. For affine $\phi$, we say that $\coef(\phi)\in \Span\{e_n\}_{n=1}^N$ if for all $0\le j\le m$, 
$\alpha_j \in \Span\{e_n\}_{n=1}^N$, and we denote the set of such estimators by \[ \mathcal{M}_{\textrm{aff},N}= \{\phi \in \mathcal{M}_{\text{aff}}: \coef(\phi)\in \Span\{e_n\}_{n=1}^N \}.\]\label{nomencl:MaffN}\nomenclature{$\mathcal{M}_{\textrm{aff},N}$}{affine subspace of $\mathcal M$ with $\coef(\phi)\in \Span\{e_n\}_{n=1}^N$}{nomencl:MaffN}
% We say that $\phi$ is affine if it can be written as
% \[
% \phi(Y_1,\ldots,Y_m) = \alpha_0 + \sum_{j=1}^m \alpha_j Y_j,
% \]
%where $\alpha_j\in C(\mathcal{D})~\forall 0\leq j\leq m$. 
% We define $\supp(P)$ for a probability measure $P$ on the space $\mathcal{H}$ \viet{$P$ is on $\mathcal{H} \times \mathbb{R}^m$?} as the smallest closed subset $A$ such that $P(A^c)=0$. \viet{there's also a proper definition of support in Aliprantis and Border's book}
%Note that the supreme norm and $L^2$ norm 
For convenience, we write
\[
\Obj(\phi,P) = \mathbb{E}_{P}\left[\|b-\phi(Y_1,\ldots,Y_m)\|^2_{L^2(\mathcal{D})}\right]
\]
for the Bayes risk of an estimator $\phi\in \mathcal M$ under the distribution $P$. We denote the Wasserstein ball around $P_0$ as
\[
\mathcal{W}(\delta) = \{P\in\mathcal{P}: \W(P,P_0)\leq\delta\}.
\]
\label{nomencl:Wassdelta}\nomenclature{$\mathcal{W}(\delta)$}{Wassertein ball around $P_0$ of radius $\delta$}{nomencl:Wassdelta}Then, denoting the Wasserstein balls around $P_0$ arising from perturbations in the first $N$ coordinates by
\begin{equation}\label{eq:revision_Wassdelta}
\mathcal W_N(\delta) = \{P\in \mathcal{W}(\delta): \mathcal L^P \big(\{\langle b,e_{n}\rangle\}_{n=N+1}^\infty\big)=\mathcal L^{P_0} \big(\{\langle b,e_{n}\rangle\}_{n=N+1}^\infty\big)\}, 
\end{equation}
\label{nomencl:WassdeltaN}\nomenclature{$\mathcal W_N(\delta)$}{Wasserstein balls around $P_0$ arising from perturbations in the first $N$ coordinates}{nomencl:WassdeltaN}where $\mathcal{L}^P$ denotes the law under $P$.
We have for any $N\ge 1$ that
\begin{align}
 \inf_{\phi \in \mathcal M_{\text{aff},N}} \sup_{P\in\mathcal{W}(\delta)} \Obj(\phi,P) & \geq \inf_{\phi \in  \mathcal M_{\text{aff}}} \sup_{P\in\mathcal{W}(\delta)}\Obj(\phi,P)  \notag\\
 & \geq\inf_{\phi\in \mathcal M}\sup_{P\in\mathcal{W}(\delta)}\Obj(\phi,P) \notag\\
 &\geq \sup_{P\in\mathcal{W}(\delta)}\inf_{\phi\in\mathcal M}\Obj(\phi,P) \notag\\
 & \geq \sup_{P\in\mathcal{W}_N(\delta)} \inf_{\phi\in \mathcal{M}}\Obj(\phi,P).\label{eq:sandwich}
\end{align}
Above, the first two inequalities follow from the inclusion $\mathcal M_{\text{aff},N} \subset \mathcal M_{\text{aff}} \subset \mathcal M$, the third inequality follows from weak duality, and the last equality follows from the inclusion $\mathcal W_N(\delta) \subset \mathcal W(\delta)$.
Finally, for any $P\in\mathcal P$ and $N\geq1$, we denote the joint measure induced by projecting $b$ onto $\Span\{e_n\}_{n=1}^N$ whilst keeping $\epsilon$ intact by\label{nomencl:PuN}\nomenclature{$P^{(N)}$}{Projection of $P$ onto $\Span\{e_n\}_{n=1}^N$ whilst keeping $\epsilon$ intact}{nomencl:PuN}
\begin{equation}
    P^{(N)} \coloneqq \mathcal L^P \Big( \big(\sum_{1\le n\le N} e_n \langle b,e_n\rangle, \epsilon \big) \Big).\label{truncatedmeasure}
\end{equation}

Our proof consists of showing the following three claims:

\textbf{Claim 1}: The finite-dimensional version of strong duality holds, i.e.,
%\viet{something wrong in the below equation, i think an = sign missing}
\begin{align}
\inf_{\phi \in \mathcal M} \sup_{P\in \mathcal W_N(\delta)}~ \Obj(\phi,P^{(N)}) \notag &=\inf_{\phi\in \mathcal M_{\text{aff},N}} \sup_{P\in \mathcal W_N(\delta)} ~\Obj(\phi,P^{(N)}) \notag\\
&=\sup_{P \in \mathcal W_N(\delta)} \inf_{\phi \in \mathcal M}~\Obj(\phi,P^{(N)})\label{eq:finitedimduality} .
\end{align}

\textbf{Claim 2}: The truncation of $P$ to $P^{(N)}$ in the last term of~\eqref{eq:sandwich} preserves the chain of inequalities, i.e., 
\begin{equation}\label{eq:claimtoprove1}
\sup_{P\in \mathcal W_N(\delta)} \inf_{\phi\in \mathcal M}~\Obj(\phi,P) \geq \sup_{P\in \mathcal W_N(\delta)} \inf_{\phi\in \mathcal{M}}~\Obj(\phi,P^{(N)}).
\end{equation}

%We will later show that 
%\begin{equation}\label{eq:claimtoprove1}
%\sup_{P:\W(P,P_0)\leq \delta,\langle b,e_n\rangle \stackrel{d}{=}\langle b^0,e_n\rangle\forall n>N} \inf_{\phi}\Obj(\phi,P) -  \sup_{P:\W(P,P_0)\leq \delta,\langle b,e_n\rangle \stackrel{d}{=}\langle b^0,e_n\rangle\forall n>N} \inf_{\phi}\Obj(\phi,P^{(N)}) \geq 0,
%\end{equation}
%and
\textbf{Claim 3}: The truncation of $P$ to $P^{(N)}$ in the first term of~\eqref{eq:sandwich} has an error which is asymptotically negligible, i.e., as $N\to\infty$
\begin{equation}\label{eq:claimtoprove2}
\inf_{\phi\in \mathcal M_{\text{aff},N} } \sup_{P\in\mathcal W(\delta)}~\Obj(\phi,P)  - \inf_{\phi\in \mathcal M_{\text{aff},N} } \sup_{P\in \mathcal W_N(\delta)}~ \Obj(\phi,P^{(N)}) = o(1).
\end{equation}
%We now aim to show that the finite-dimensional version of the strong duality holds, i.e.,
%\begin{align}
%&\inf_{\phi,\text{affine},\coef(\phi)\in \Span\{e_n\}_{n=1}^N } \sup_{P\in \mathcal P(\delta, N)} \Obj(\phi,P^{(N)}) \notag\\
%&=\sup_{P:\W(P,P_0)\leq \delta,\langle b,e_n\rangle \stackrel{d}{=}\langle b^0,e_n\rangle\forall n>N} \inf_{\phi}\Obj(\phi,P^{(N)})\label{eq:finitedimduality} .
%\end{align}
Combining the above three claims and the chain of inequalities~\eqref{eq:sandwich}, we conclude that 
\begin{equation}\label{eq:thmstrongdualityregressappendix}
\inf_{\phi \in \mathcal M }\sup_{P\in\mathcal W(\delta)}~\Obj(\phi,P) = \sup_{P\in\mathcal W(\delta)}\inf_{\phi \in \mathcal M}~\Obj(\phi,P)
\end{equation}
by letting $N\to\infty$.
\end{proof}

We now provide the proof of the three claims.

\begin{proof}[Proof of Claim 1] \ssw{Got rid of the definition of marginal} 
For $P\in\mathcal{P}$, we denote the marginal distribution of the first $N$ basis coefficients $\langle b, e_n\rangle$ and $\epsilon$ by
\begin{equation}\label{QN}
    Q^{(N)}\coloneqq \mathcal L^P \big( (\langle b,e_1\rangle,\cdots,\langle b,e_N\rangle,\epsilon) \big).
\end{equation}
Note that $Q^{(N)}$ is a probability measure on $\mathbb{R}^{N+m}$, while the distribution $P^{(N)}$ is supported on $C(\mathcal D) \times \R^m$. In particular, we have $Q_0^{(N)} \coloneqq \mathcal L^{P_0} \big( (\langle b,e_1\rangle,\cdots,\langle b,e_N\rangle,\epsilon) \big)$ is a multivariate Gaussian. 
Now, for the weighted cost function
\[
c_N(r,s) = \sum_{n=1}^N (r_n-s_n)^2 w_n + \sum_{j=1}^m (r_{N+j} - s_{N+j})^2 \quad \forall  r,s\in\mathbb{R}^{N+m}
\]
(where the variable order of $r$ (or $s$) is the same as the variable order in equation~\eqref{QN}) on $\R^{N+m}$, we denote the corresponding optimal transport cost between probability measures $\tau$ and $\nu$ on $\R^{N+m}$ by
\[
\W_N^2(\tau,\nu)= \min_{\pi: \pi_r = \tau, \pi_s= \nu}\mathbb{E}_{(R,S)\sim \pi} [c_N(R,S)],
\]
where we implicitly imposed $\pi$ as a probability measure to avoid cluttered notations.
Thus, by writing the shorthand \[
g(r) \coloneqq\left\|\sum_{n=1}^Nr_n e_n - \phi\Big(\big(\sum_{n=1}^N r_nT( e_n)(x_j)+ r_{N+j}\big)_{j=1}^m\Big)\right\|^2_{L^2(\mathcal{D})},
\]
we obtain that
\begin{align*}
\sup_{P\in \mathcal W_N(\delta)} \Obj(\phi,P^{(N)}) = \sup_{Q^{(N)}:\W_{N}(Q^{(N)},Q_0^{(N)})\leq\delta}\mathbb{E}_{R\sim Q^{(N)}}[g(R)].
\end{align*}
Assume momentarily that $\phi\in \mathcal M_{\text{aff},N}$. In this case, $g$ is a convex (non-constant) quadratic function and \cite[Theorem~1]{ref:blanchet2019quantifying} implies the dual formulation
\begin{align}
& \sup_{Q^{(N)}:\W_N(Q^{(N)},Q_0^{(N)})\leq\delta}\mathbb{E}_{R\sim Q^{(N)}}[g(R)] = \inf_{\gamma\geq0}\left(\gamma\delta^2 + \mathbb{E}_{R\sim Q_0^{(N)}}\left[\sup_{s\in\mathbb{R}^{N+m}}\left\{g(s)-\gamma c_N(R,s)\right\}\right]\right)\label{eq:subprimal}.
\end{align}
The same theorem moreover implies that there exists a dual optimizer $\gamma^\star$ to the right-hand side minimization problem of~\eqref{eq:subprimal}. Since the distributions in the set
\[
\{Q^{(N)}:\W_N(Q^{(N)},Q_0^{(N)})\leq\delta\}
\]
are uniformly bounded in the second moment and $g$ is quadratic in $r$, we have that 
\[
\sup_{Q^{(N)}:\W_N(Q^{(N)},Q_0^{(N)})\leq\delta}\mathbb{E}_{R\sim Q^{(N)}}[g(R)]<\infty.
\] 
Since also $\sup_{s\in\mathbb{R}^{N+m}}g(s)=\infty$, we have necessarily that $\gamma^\star>0$.

Observe that $g(r)$ can be written in the form $g(r) = r^\top Gr + c^\top r + \|\alpha_0\|^2_{L^2(\mathcal{D})}$  for some $G\in\mathbb{S}_{+}^{N+m}$ and $c\in\mathbb{R}^{N+m}$, where we denote $\mathbb{S}_{+}^{N+m}$\label{nomencl:psdm}\nomenclature{$\mathbb{S}_{+}^{m}$}{positive semi-definite matrix of dimension $m$-by-$m$}{nomencl:psdm} as the set of positive semi-definite matrices with a dimension of $N+m$. Denoting \review{the matrix (note that not to be confused with $\W_N$, which denotes optimal transport distance)}\label{nomencl:WmatN}\nomenclature{$W_N$}{matrix with diagonals $w_n$}{nomencl:WmatN} 
\[W_N=\diag(w_1,\ldots,w_N,1,\ldots,1),\] 
then, almost surely for $R\sim Q_0^{(N)}$, we have
\begin{align*}
g(s)- \gamma^\star c_N(R,s) &= s^\top (G-\gamma^\star W_N) s + (2\gamma^\star R^\top W_N  +c ^\top) s +\|\alpha_0\|^2_{L^2(\mathcal{D})}  + \gamma^\star R^\top W_NR.
\end{align*}
Since $R$ follows a non-degenerate multivariate Gaussian distribution, we have that necessarily $\gamma^\star W_N - G \in\mathbb{S}^{N+m}_{++}$, otherwise $\sup_{s\in\mathbb{R}^{N+m}}g(s) - \gamma^\star c_N(R,s)=\infty$ almost surely. 
Note that~\cite[Theorem~1]{ref:blanchet2019quantifying} also implies that the optimal coupling between $Q_0^{(N)}$ and the maximizer to the left side of \eqref{eq:subprimal}, if it exists, must be given as the law of $(R,s^\star(R))$, for the affine push-forward map
\begin{equation}\label{eq:optimalplanfd}
s^\star(R) = -\frac{1}{2} (G-\gamma^\star W_N)^{-1}(2\gamma^\star W_N R + c).
\end{equation}
The existence of a solution $Q_
\star^{(N)}$ to the left hand side of \eqref{eq:subprimal} as well as of an optimal coupling between $Q_0^{(N)}$ and $Q_
\star^{(N)}$, which we denote by $\pi^\star$, can be verified by~\cite[Corollary 1(i)]{ref:gao2016distributionally}. Indeed, since $\gamma^\star W_N - G\in\mathbb{S}^{N+m}_{++}$, for the growth rate $\kappa$ defined to be 
\[
\kappa = \limsup_{r\to\infty}\frac{r^\top G r}{r^\top W_N r},
\]
we have that $\kappa<\gamma^\star$ \review{so that Condition (a) in \cite[Corollary 1(i)]{ref:gao2016distributionally} holds}.
%{\color{red}in which case}
%or $\kappa =\gamma^\star$ and equation~(20) in CITE if $\gamma^\star W_N -  G \in\mathbb{S}^{N+m}_+$. 
%{\color{red}in which case}
%Moreover $\pi^\star$ is unique and 
%concentrates on the graph of a non-degenerate linear transform of $R$. 
By the affine push-forward map~\eqref{eq:optimalplanfd}, $Q_\star^{(N)}$ is also Gaussian, whence we may write~(\ref{eq:subprimal}) as
\begin{align*}
\sup_{Q^{(N)}:\W_N(Q^{(N)},Q_0^{(N)})\leq\delta}\mathbb{E}_{R\sim Q^{(N)}}[g(R)]
& = \sup_{Q^{(N)}:\W_N(Q^{(N)},Q_0^{(N)})\leq\delta,Q^{(N)}\textrm{  normal}}\mathbb{E}_{R\sim Q^{(N)}}[g(R)]\\
%& = \sup_{Q^{(N)}:\W_N(Q^{(N)},Q_0^{(N)})\leq\delta, Q^{(N)} = (A,a)^{\#}Q_0^{(N)}}\mathbb{E}_Q[g(R)]\\
& = \sup_{(\mu,\Sigma)\in \mathcal{S}_N} \langle G, \Sigma+\mu\mu^\top \rangle + c^\top \mu + \|\alpha_0\|^2_{L^2(\mathcal{D})}\\
& \geq  \sup_{(\mu,\Sigma)\in \mathcal{S}_N} \langle G, \Sigma+\mu\mu^\top \rangle,
\end{align*}
where the last inequality is because $c^\top\mu\geq0$ after a possible sign change in $\mu$. Here $\mathcal{S}_N$ is a compact and convex set\label{nomencl:SSN}\nomenclature{$\mathcal{S}_N$}{modified Gelbrich distance}{nomencl:SSN} 
\[
\mathcal{S}_N = \left\{(\mu,\Sigma):\|\mu\|_2^2 + tr(W_N\Sigma)+tr(W_N\Sigma_0) - 2tr\left[\left(\sqrt{\Sigma_0}W_N\Sigma W_N\sqrt{\Sigma_0}\right)^{1/2}\right]\leq\delta^2\right\},
\]
where $\Sigma_0 \in \PD^{N+m}$\label{nomencl:Sigma0}\nomenclature{$\Sigma_0$}{covariance matrix of $Q_0^{(N)}$}{nomencl:Sigma0} is the covariance matrix of $Q_0^{(N)}$. To see how the set $\mathcal{S}_N$ arises, note that the usual squared-Euclidean cost function is $ \|r-s\|_2^2$ between $r$ and $s\in\mathbb{R}^{N+m}$, which gives rise to \review{the usual Gelbrich distance (see~\cite[Proposition 2.2]{ref:shafieezadeh2018wasserstein})
\[
\sqrt{\|\mu\|_2^2 + tr(\Sigma)+tr(\Sigma_0) - 2tr\left[\left(\sqrt{\Sigma_0}\Sigma \sqrt{\Sigma_0}\right)^{1/2}\right]}
\]
that coincides with the type-2 Wasserstein distance between multivariate Gaussians}. Herein, our new cost is $(r-s)^\top W_N(r-s)$, thus the optimal coupling $\pi^\star$ for the new cost solves
\[
\min_{\pi: \pi_r = \tau, \pi_s= \nu} \left(-\int r^\top W_N s \mathrm{d}\pi(r,s)\right).
\]
Using substitution of variables $\tilde r = \sqrt{W_N}r$ and $\tilde s = \sqrt{W_N}s$, we find that the above problem reduces to the optimal transport problem with a squared-Euclidean cost.
% there is a one-to-one correspondence between $\pi(r,s)$ and $\tilde \pi(\tilde r,\tilde s)$, and $\tilde \pi^\star$ (corresponding to $\pi^\star$) solves
% \[
% \min_{\tilde\pi(\tilde r,\tilde s)} \left(-\int \tilde r^\top \tilde s \mathrm{d}\tilde\pi(\tilde r,\tilde s)\right).
% \]

In the sequel we denote the Wasserstein balls around $P_0$ arising from perturbations in the first $N$ coordinates and restricting to the family of Gaussian distributions as\label{nomencl:wnorN}\nomenclature{$\mathcal{W}_{\textrm{nor},N}(\delta)$}{Wasserstein balls around $P_0$ arising from perturbations in the first $N$ coordinates and restricting to the family of Gaussian distributions}{nomencl:wnorN} 
\begin{equation}\label{eq:wnorN_revision}
\mathcal{W}_{\textrm{nor},N}(\delta) = \{P\in \mathcal W_N(\delta): Q^{(N)}\textrm{ is centered full-rank normal}\}.
\end{equation}
Therefore, we obtain \review{(recall $\alpha_0$ in Definition~\ref{def:affinepredictor_revision})}\ssw{Fix this $\alpha_0$ business.{\color{red}done}}
\begin{align*}
&\inf_{\phi\in\mathcal{M}_{\textrm{aff},N},\alpha_0=0}  \sup_{Q^{(N)}:\W_N(Q^{(N)},Q_0^{(N)})\leq\delta}\mathbb{E}_{R\sim Q^{(N)}}[g(R)] \\
& \ge\inf_{\phi\in\mathcal{M}_{\textrm{aff},N}  }  \sup_{Q^{(N)}:\W_N(Q^{(N)},Q_0^{(N)})\leq\delta}\mathbb{E}_{R\sim Q^{(N)}}[g(R)]\\
%&  \inf_{\phi(\cdot),\text{affine},\coef(\phi)\in \Span\{T(e_n)\}_{n=1}^N } \sup_{(\mu,\Sigma)\in \mathcal{S}} \langle G, \Sigma+\mu\mu^\top \rangle + \sum_{n=1}^N\langle \alpha_0,e_n\rangle \langle g_n,\mu\rangle + \sum_{n=1}^N(\langle \alpha_0,e_n\rangle)^2\\
&\ge \inf_{\phi\in\mathcal{M}_{\textrm{aff},N}
}  \sup_{(\mu,\Sigma)\in \mathcal{S}_N} \langle G, \Sigma + \mu\mu^\top\rangle\\
& \ge \inf_{\phi\in\mathcal{M}_{\textrm{aff},N}
}  \sup_{(\mu,\Sigma)\in \mathcal{S}_N} \langle G, \Sigma\rangle\\
& =\inf_{\phi\in\mathcal{M}_{\textrm{aff},N},\alpha_0=0}  \sup_{(0,\Sigma)\in \mathcal{S}_N} \langle G, \Sigma \rangle\\
& = \inf_{\phi\in\mathcal{M}_{\textrm{aff},N},\alpha_0=0}  \sup_{Q^{(N)}:\W_N(Q^{(N)},Q_0^{(N)})\leq\delta}\mathbb{E}_{R\sim Q^{(N)}}[g(R)] ,
\end{align*}
where the last equality is because $s^\star(R)$ is zero-mean whenever $c=0$, which in turn follows from $\alpha_0=0$. Hence all inequalities become equalities, and we conclude that
\[
\inf_{\phi\in\mathcal{M}_{\textrm{aff},N}  }  \sup_{Q^{(N)}:\W_N(Q^{(N)},Q_0^{(N)})\leq\delta}\mathbb{E}_{R\sim Q^{(N)}}[g(R)] = \inf_{\phi\in\mathcal{M}_{\textrm{aff},N},\alpha_0=0}  \sup_{(0,\Sigma)\in \mathcal{S}_N} \langle G, \Sigma \rangle.
\]
Remind that the matrix $G$ is dependent on the coefficients of the affine estimator $\phi$. By Sion's minimax theorem~\cite{ref:sion1958minimax}, we have
\begin{align}
%&  \min_{\phi(\cdot),\text{affine},\coef(\phi)\in \Span\{T(e_n)\}_{n=1}^N } \sup_{Q\in\mathbb{R}^N:\W(Q,Q_0^{(N)})\leq\delta_N}\\
&   \inf_{\phi\in\mathcal{M}_{\textrm{aff},N},\alpha_0=0}  \sup_{(0,\Sigma)\in \mathcal{S}_N} \langle G, \Sigma\rangle\label{eq:minmaxfd}
\\
%& =  \sup_{(\mu,\Sigma)\in \mathcal{S}_N} \inf_{\phi(\cdot),\text{affine},\coef(\phi)\in \Span\{T(e_n)\}_{n=1}^N,\alpha_0=0}\langle G, \Sigma\rangle \notag\\
& =   \sup_{(0,\Sigma)\in \mathcal{S}_N} \inf_{\phi\in\mathcal{M}_{\textrm{aff},N},\alpha_0=0}\langle G, \Sigma\rangle\label{eq:maxminfd} \\
%& = \sup_{Q\in\mathbb{R}^{N+d}:\W(Q,Q_0^{(N)})\leq\delta_N, Q\textrm{ is centered normal}}\inf_{\phi(\cdot),\text{affine},\coef(\phi)\in \Span\{T(e_n)\}_{n=1}^N,\alpha_0=0}\mathbb{E}_{Q}[g(R)]\notag\\
& = \sup_{Q^{(N)}:\W_N(Q^{(N)},Q_0^{(N)})\leq\delta, Q^{(N)}\textrm{ is centered full-rank normal}}\inf_{\phi\in\mathcal{M}_{\textrm{aff},N},\alpha_0=0}\mathbb{E}_{R\sim Q^{(N)}}[g(R)].\notag
\end{align}
The last ``full-rank'' assertion comes from the fact that $s^\star(R)$ is a non-degenerate linear transformation of $R$, and that a linear transformation of a multivariate Gaussian is also Gaussian.  
% Moreover, by Shafieezadeh-Abadeh et al. 2018 (Lemma A.2), we may restrict $\Sigma\succeq \Sigma_0$ in~(\ref{eq:maxminfd}).{\color{red}:not clear} 
Thus
\begin{align}
\inf_{\phi\in\mathcal{M}_{\textrm{aff},N} }\sup_{P\in\mathcal{W}_N(\delta)} \Obj(\phi,P^{(N)})\notag & =  \inf_{\phi\in\mathcal{M}_{\textrm{aff},N}} \sup_{Q^{(N)}:\W_N(Q^{(N)},Q_0^{(N)})\leq\delta}\mathbb{E}_{R\sim Q^{(N)}}[g(R)]\notag \\
& =  \inf_{\phi\in\mathcal{M}_{\textrm{aff},N},\alpha_0=0} \sup_{Q^{(N)}:\W_N(Q^{(N)},Q_0^{(N)})\leq\delta}\mathbb{E}_{R\sim Q^{(N)}}[g(R)]\notag \\
%& = \inf_{\phi(\cdot),\text{affine},\coef(\phi)\in \Span\{T(e_n)\}_{n=1}^N }\sup_{P:\W(P,P_0^{(N)})\leq \delta_N,\supp(P)\subseteq \Span\{T(e_n)\}_{n=1}^N,h^{\#}(P)\in\mathcal{N}_N}\notag\\
&= \sup_{P\in \mathcal{W}_{\textrm{nor},N}(\delta)}\inf_{\phi\in\mathcal{M}_{\textrm{aff},N},\alpha_0=0}\Obj(\phi,P^{(N)})\label{eq:psubnreduction},
%& \geq \sup_{P:\W(P,P_0^{(N)})\leq \delta_N,\supp(P)\subseteq \Span\{T(e_n)\}_{n=1}^N,\coef(P)\in\mathcal{N}}\inf_{\phi(\cdot),\text{affine},\coef(\phi)\in \Span\{T(e_n)\}_{n=1}^N }\\
\end{align}
\review{where recall that $\mathcal{W}_{\textrm{nor},N}(\delta)$ was defined in equation~\eqref{eq:wnorN_revision}.}
For $\W(P,P_0)\leq \delta$, note that $b$ is a process with continuous sample paths, whence we can interchange the integration
\begin{align*}
\mathbb{E}_P\left[\| b-\phi(Y_1,\ldots,Y_m)\|^2_{L^2(\mathcal{D})}\right] &= \mathbb{E}_P\left[\int_\mathcal{D} |b(x) - \phi(Y_1,\ldots,Y_m)(x)|^2 \dd x\right] \\
& = \int_\mathcal{D} \mathbb{E}_P\left[|b(x) -\phi(Y_1,\ldots,Y_m)(x)|^2\right] \dd x.
\end{align*}
Thus the optimal solution to the estimation problem
\begin{equation}\label{eq:potsol}
\inf_{\phi\in\mathcal{M}} \mathbb{E}_P\left[\| b-\phi(Y_1,\ldots,Y_m)\|^2_{L^2(\mathcal{D})}\right]
\end{equation}
is given by the conditional expectation function
\begin{equation}\label{eq:condimean}
\phi(Y_1,\ldots,Y_m)(x)=\mathbb{E}_{P}\left[b(x)|Y_1,\ldots,Y_m\right].
\end{equation}
For any distribution $P$ that is a centered Gaussian random variable, it is easy to see that for $x\in\mathcal{D}$, the random vector
\[
(b(x),Y_1 ,\ldots,Y_m) = (b(x), b(x_1) +\epsilon_1,\ldots,b(x_m)+\epsilon_m)
\]
is jointly Gaussian and the conditional expectation satisfies
\begin{equation}\label{eq:jointgaussian}
\mathbb{E}_{P}\left[b(x)|Y_1 ,\ldots,Y_m\right]= \left(k_{\epsilon,P}(x,x_1),\ldots,k_{\epsilon,P}(x,x_m)\right)\cdot (K_{\epsilon,P})^{-1} \cdot \left(Y_1 ,\ldots,Y_m\right)^\top
\end{equation}
where $k_{\epsilon,P}(x,x_j) = \mathbb{E}_P[b(x)(b(x_j)+\epsilon_j)]$, and 
\review{
\[
K_{\epsilon,P} = \left(\mathbb{E}_P[(b(x_i)+\epsilon_i)(b(x_j)+\epsilon_j)]\right)_{ij} \in \PD^{m}
\]}
provided the matrix $K_{\epsilon,P}$ is invertible. Moreover, 
the optimal value of the Bayes risk~\eqref{eq:potsol} is given by
\begin{align*}
&\int_\mathcal{D} \mathbb{E}_P[\mathbb{V}\mathrm{ar}_P[b(x)|Y_1 ,\ldots,Y_m]]\dd x\\
& = \int_\mathcal{D}  k_P(x,x) -  \left(k_{\epsilon,P}(x,x_1),\ldots,k_{\epsilon,P}(x,x_m)\right)\cdot (K_{\epsilon,P})^{-1} \cdot\left(k_{\epsilon,P}(x,x_1),\ldots,k_{\epsilon,P}(x,x_m)\right)^\top \dd x
\end{align*}
where $k_P(x,x) = \mathbb{E}_P[b(x)b(x)]$.\label{nomencl:varPb}\nomenclature{$k_P(x,x)$}{variance function of $b(x)$ under $P$}{nomencl:varPb} 
% For any $P^{(N)}$ that satisfies the constraint of equation~\eqref{eq:psubnreduction}, the optimal solution $\phi$ defined by~(\ref{eq:condimean}) and~(\ref{eq:jointgaussian}) under the measure $P^{(N)}$ is affine ($\phi$ is well-defined since $K_\epsilon$ is invertible), $\coef(\phi)\in \Span\{T(e_n)\}_{n=1}^N$ and $\alpha_0=0$. 
Thus
\begin{align*}
\sup_{P\in \mathcal{W}_{\textrm{nor},N}(\delta)}\inf_{\phi\in\mathcal{M}_{\textrm{aff},N},\alpha_0=0}\Obj(\phi,P^{(N)})\notag
%&\qquad \mathbb{E}_P\left[\| u(\cdot)-\phi(\cdot;u(x_1),\ldots,u(x_d))\|^2\right]\\
&=\sup_{P\in\mathcal{W}_{\textrm{nor},N}(\delta)}\inf_{\phi \in \mathcal M}\Obj(\phi,P^{(N)})\notag\\
%&\qquad \mathbb{E}_P\left[\| u(\cdot)-\phi(\cdot;u(x_1),\ldots,u(x_d))\|^2\right]\\
& \leq \sup_{P\in\mathcal{W}_N(\delta)}\inf_{\phi \in \mathcal M}\Obj(\phi,P^{(N)})\notag,
\end{align*}
On the other hand,
\begin{align*}
\sup_{P\in\mathcal{W}_N(\delta)}\inf_{\phi \in \mathcal M}\Obj(\phi,P^{(N)}) \leq \inf_{\phi \in \mathcal M}\sup_{P\in\mathcal{W}_N(\delta)}\Obj(\phi,P^{(N)})
%&\leq \sup_{P:\W(P,P_0)\leq \delta,\supp(P)\subseteq \Span\{T(e_n)\}_{n=1}^N}\inf_{\phi(\cdot),\text{affine},\coef(\phi)\in \Span\{T(e_n)\}_{n=1}^N }\\
\leq\inf_{\phi\in\mathcal{M}_{\textrm{aff},N} }\sup_{P\in\mathcal{W}_N(\delta)}\Obj(\phi,P^{(N)}).
\end{align*}
Thus combining the above chains of inequalities and the equality~\eqref{eq:psubnreduction}, we have
\begin{align}
\inf_{\phi\in\mathcal{M}_{\textrm{aff},N} }\sup_{P\in\mathcal{W}_N(\delta)}\Obj(\phi,P^{(N)})\notag
& = \inf_{\phi\in\mathcal{M}_{\textrm{aff},N},\alpha_0=0}\sup_{P\in\mathcal{W}_N(\delta)}\Obj(\phi,P^{(N)})\notag \\
& = \inf_{\phi \in \mathcal M} \sup_{P\in\mathcal{W}_N(\delta)}\Obj(\phi,P^{(N)})\notag\\
& = \sup_{P\in\mathcal{W}_N(\delta)}\inf_{\phi \in \mathcal M}\Obj(\phi,P^{(N)})\label{eq:sandwichpart}\\
& = \sup_{P\in \mathcal{W}_{\textrm{nor},N}(\delta)}\inf_{\phi \in \mathcal M}\Obj(\phi,P^{(N)})\notag.
\end{align}
%On the other hand,
%\begin{align*}
%&\inf_{\phi(\cdot),\text{affine},\coef(\phi)\in \Span\{T(e_n)\}_{n=1}^N }\sup_{P:\W(P,P_0)\leq \delta,\supp(P)\subseteq \Span\{T(e_n)\}_{n=1}^N} \\
%& \geq\sup_{P:\W(P,P_0)\leq \delta,\supp(P)\subseteq \Span\{T(e_n)\}_{n=1}^N}\inf_{\phi(\cdot),\text{affine},\coef(\phi)\in \Span\{T(e_n)\}_{n=1}^N }\\
%& \geq \sup_{P:\W(P,P_0^{(N)})\leq \delta_N,\supp(P)\subseteq \Span\{T(e_n)\}_{n=1}^N,P\in h^{-1}(\mathcal{N}_N)}\inf_{\phi(\cdot),\text{affine},\coef(\phi)\in \Span\{T(e_n)\}_{n=1}^N }
%\end{align*}
Thus we have established~\eqref{eq:finitedimduality}. 
\end{proof}
\begin{proof}[Proof of Claim 2]
By the proof of Claim~1, there exists a Nash equilibrium of~\eqref{eq:finitedimduality}, which we denote by $(\phi^\star_N,P^\star_N)$. In particular, by the previous Claim~1,\ssw{Specify 'by previous claim 1'} we can choose $\phi^\star_N$ as affine and $Q_N^{\star} \coloneqq \mathcal L^{P^\star_N} \big( \langle b,e_1\rangle,\cdots,\langle b,e_N\rangle,\epsilon\big)$ as some centered full-rank normal distribution. Also, the pair $(\phi^\star_N,Q_N^{\star})$ is the unique Nash equilibrium corresponding to the game~\eqref{eq:minmaxfd} and~\eqref{eq:maxminfd}.\ssw{This sentence is not written well and grammar is wrong.} Moreover, under the law $P^\star_N$, the random elements 
\begin{equation}\label{eq:independentassertionclaim2_app}
(\{\langle b,e_n\rangle\}_{n=1}^N,\epsilon)\quad\textrm{ and }\quad\{\langle b,e_{n}\rangle\}_{n=N+1}^\infty
\end{equation}
are independent. Finally, we have that
\[
 \mathcal L^{P^\star_N} \big(\{\langle b,e_{n}\rangle\}_{n=N+1}^\infty\big) = \mathcal L^{P_0} \big(\{\langle b,e_{n}\rangle\}_{n=N+1}^\infty\big).
\]
Note that $P^\star_N$ is a Gaussian measure supported on $C(\mathcal{D})$. We now show our claim~\eqref{eq:claimtoprove1}. Note that 
\[
Y_j = u(x_j) +\epsilon_j = \sum_{k=1}^N\langle b,e_k\rangle T(e_k)(x_j)+\epsilon_j + \left(u(x_j)-\sum_{k=1}^N\langle b,e_k\rangle T(e_k)(x_j)\right),
\]
where the last term (denoted by $R_{j,N}$ as a shorthand) is of zero-mean under $P_N^\star$. Thus any affine estimator $\phi$ can be written as (recall  Definition~\ref{def:affinepredictor_revision})
\[
\phi(Y_1,\ldots,Y_m) = \sum_{j=1}^m\alpha_j\left(\sum_{k=1}^N\langle b,e_k\rangle T(e_k)(x_j)+\epsilon_j\right) + \sum_{j=1}^m\alpha_jR_{j,N},
\]
where w.l.o.g. $\alpha_0=0$.
\review{Let $\{\tilde e_k\}_{k=1}^{\infty}$ = $\{ e_k\}_{k=1}^{\infty}$. (In the proof of Theorem~\ref{thm:swapinverse} we modify $\{\tilde e_k\}_{k=1}^{\infty}$ to  be the Gram-Schmidt orthonormalization of $\{T(e_k)\}_{k=1}^\infty$\ssw{Is this the finite-dimensional space or the infinite dimensional one??}. See discussion in Appendix~\ref{sec:regressionproof_revision})}. 
% with $\tilde{N} = \tilde{N}(N)\to\infty$\ssw{What is $N$? That has to be specified!} as $N\to\infty$, where if $\{T(e_n)\}_{n=1}^\infty$ is finite-dimensional, then we formally add $\tilde e_n$ \ssw{Sorry this has to be written precisely.}as the zero function if necessary to avoid notation burdens.
We have
\begin{align*}
&\sup_{P\in\mathcal{W}_N(\delta)} \inf_{\phi \in \mathcal M}\Obj(\phi,P) \\
& \geq \inf_{\phi \in \mathcal M}\Obj(\phi,P^\star_N)\\
& =\inf_{\phi \in \mathcal M}\mathbb{E}_{P^\star_N}[\|b-\phi(Y_1,\ldots,Y_m)\|^2_{L^2(\mathcal{D})}]\\
& = \inf_{\phi\in\mathcal{M}_{\textrm{aff}}}\mathbb{E}_{P^\star_N}[\|b-\phi(Y_1,\ldots,Y_m)\|^2_{L^2(\mathcal{D})}]\quad\textrm{ due to Gaussianity}\\
%& \geq \inf_{\phi}\sum_{n=1}^N\mathbb{E}_{P^\star_N}(\langle u-\phi,e_n\rangle^2)\\
% & \geq \inf_{\phi\in\mathcal{M}_{\textrm{aff}}}\mathbb{E}_{P^\star_N}[\|u|_{\Span\{T(e_n)\}_{n=1}^N}-\phi(Y_1,\ldots,Y_m)|_{\Span\{T(e_n)\}_{n=1}^N}\|^2_{L^2(\mathcal{D})}]\\
& = \inf_{\phi\in\mathcal{M}_{\textrm{aff}}} \mathbb{E}_{P^\star_N}\left[\sum_{n=1}^\infty \left(\left\langle b - \phi(Y_1,\ldots,Y_m),\tilde e_n\right\rangle\right)^2\right] \\
& \geq \inf_{\phi\in\mathcal{M}_{\textrm{aff}}} \mathbb{E}_{P^\star_N}\left[\sum_{n=1}^N \left(\left\langle b - \phi(Y_1,\ldots,Y_m),\tilde e_n\right\rangle\right)^2\right]\\
& = \inf_{\phi\in\mathcal{M}_{\textrm{aff}}}\sum_{n=1}^{N } \mathbb{E}_{P^\star_N}\left(\left\langle b - \sum_{j=1}^m\alpha_j\left(\sum_{k=1}^N\langle b,e_k\rangle T(e_k)(x_j)+\epsilon_j\right),\tilde e_n\right\rangle^2\right) \\
&\quad  + 2\sum_{n=1}^N\mathbb{E}_{P^\star_N}\left(\left\langle b - \sum_{j=1}^m\alpha_j\left(\sum_{k=1}^N\langle b,e_k\rangle T(e_k)(x_j)+\epsilon_j\right),\tilde e_n\right\rangle\left\langle\sum_{j=1}^m\alpha_jR_{j,N},\tilde e_n\right\rangle\right)\\
&\quad + \sum_{n=1}^{N} \mathbb{E}_{P^\star_N}\left(\left\langle\sum_{j=1}^m\alpha_jR_{j,N},\tilde e_n\right\rangle^2\right)\\
& = \inf_{\phi\in\mathcal{M}_{\textrm{aff}}}\sum_{n=1}^{N } \mathbb{E}_{P^\star_N}\left(\left\langle b - \sum_{j=1}^m\alpha_j\left(\sum_{k=1}^N\langle b,e_k\rangle T(e_k)(x_j)+\epsilon_j\right),\tilde e_n\right\rangle^2\right) \\
&\quad + \sum_{n=1}^{N} \mathbb{E}_{P^\star_N}\left(\left\langle\sum_{j=1}^m\alpha_jR_{j,N},\tilde e_n\right\rangle^2\right)\\
&\geq \inf_{\phi\in\mathcal{M}_{\textrm{aff}}}\sum_{n=1}^{N}  \mathbb{E}_{P^\star_N}\left(\left\langle b - \sum_{j=1}^m\alpha_j\left(\sum_{k=1}^N\langle b,e_k\rangle T(e_k)(x_j)+\epsilon_j\right),\tilde e_n\right\rangle^2\right) \\
& =\sup_{P\in\mathcal{W}_N(\delta)} \inf_{\phi \in \mathcal M}\Obj(\phi,P^{(N)}).
\end{align*}\ssw{In the fifth line, I propose a different notation for the restriction. In the seventh line, 'zero-meanness' is informal. In the sixth line, we have a problem: $T(e_n)$ might not be orthogonal if $T$ is not diagonal?? }
% where $f|_{\Span\{T(e_n)\}_{n=1}^N}$ denotes the orthogonal projection of a function $f\in C(\mathcal{D})$ to the subspace $\Span\{T(e_n)\}_{n=1}^N$. 
The penultimate equality is because that
\begin{align*}
    &\,\,\mathbb{E}_{P^\star_N}\left(\left\langle b - \sum_{j=1}^m\alpha_j\left(\sum_{k=1}^N\langle b,e_k\rangle T(e_k)(x_j)+\epsilon_j\right),\tilde e_n\right\rangle\left\langle\sum_{j=1}^m\alpha_jR_{j,N},\tilde e_n\right\rangle\right)\\
    = &\,\,\mathbb{E}_{P^\star_N}\left(\left\langle b - \sum_{j=1}^m\alpha_j\left(\sum_{k=1}^N\langle b,e_k\rangle T(e_k)(x_j)+\epsilon_j\right),\tilde e_n\right\rangle\right)\mathbb{E}_{P^\star_N}\left(\left\langle\sum_{j=1}^m\alpha_jR_{j,N},\tilde e_n\right\rangle\right)\\
     = &\,\,0,
\end{align*}
where we have used the independence in~\eqref{eq:independentassertionclaim2_app} and the fact that $R_{j,N}$ have zero mean.
Therefore, we have established Claim 2.\ssw{Change notation.}
\end{proof}
\begin{proof}[Proof of Claim 3]
% \begin{align}
% &\lim_{N\to\infty} \inf_{\phi(\cdot),\text{affine},\coef(\phi)\in \Span\{T(e_n)\}_{n=1}^N }\sup_{P:\W(P,P_0)\leq \delta,\supp(P|_{\mathcal{H}})\subseteq \Span\{T(e_n)\}_{n=1}^N}=\notag\\
% &\lim_{N\to\infty}\inf_{\phi(\cdot),\text{affine},\coef(\phi)\in \Span\{T(e_n)\}_{n=1}^N } \sup_{P:\W(P,P_0)\leq\delta}.\label{eq:errorfd}
% \end{align}
% Denote $\phi_N^\star$ to be the linear predictor that solves 
%  \[
%  \inf_{\phi(\cdot),\text{affine},\coef(\phi)\in \Span\{T(e_n)\}_{n=1}^N }\sup_{P:\W(P,P_0)\leq \delta,\supp(P|_{\mathcal{H}})\subseteq \Span\{T(e_n)\}_{n=1}^N}\Obj(\phi,P).
%  \]
Since $\phi^\star_N$ is affine, we can write $\phi^\star_N$ in the form $\phi^\star_N(Y_1,\ldots,Y_m) = \sum_{j=1}^m \alpha_j^\star Y_j$, where $\alpha_j^\star\in\Span\{e_n\}_{n=1}^N$ and w.l.o.g. $\alpha_0^\star=0$. Here, we have suppressed the dependence of $\alpha_j^\star$ on $N$ to simplify the notations. Also note that there exists a constant $C$ independent of $N$, such that for any $1\leq j\leq m$,
\[
\|\alpha_j^\star\|^2_{L^2(\mathcal{D})} = \sum_{n=1}^{N} |\langle \alpha_j^\star,\tilde e_n\rangle|^2\leq C<\infty.
\]\ssw{This is not true. Or, we need to properly define $\tilde e_n$!}
This is because we can choose $(0,\Sigma)\in \mathcal{S}_N$ in the inner constraint of~\eqref{eq:minmaxfd} to be the nominal measure, and thereby conclude that the optimal objective value in~\eqref{eq:minmaxfd} is lower bounded by $C^{'}\sum_{j=1}^m \|\alpha_j\|^2_{L^2(\mathcal{D})}$, for some constant $C^{'}$ independent of $N$. For any $N$, we have
\begin{align}
&\sup_{P\in\mathcal{W}(\delta)} \Obj(\phi^\star_N,P) -\sup_{P\in\mathcal{W}_N(\delta)}\Obj(\phi_N^\star,P^{(N)})\notag\\
& =\sup_{P\in\mathcal{W}(\delta)}\left(\sum_{n=1}^{N} \mathbb{E}_{P}(\langle b-\phi_N^\star,\tilde e_n\rangle)^2+\sum_{n=N+1}^\infty\mathbb{E}_{P}(\langle b,\tilde e_n\rangle)^2 \right)\notag\\
&\qquad-\sup_{P\in\mathcal{W}_N(\delta)}\sum_{n=1}^{N} \mathbb{E}_{P^{(N)}}(\langle b-\phi_N^\star,\tilde e_n\rangle)^2\notag\\
 &\leq \sup_{P\in\mathcal{W}(\delta)} \sum_{n=1}^{N} \mathbb{E}_{P}(\langle b-\phi_N^\star,\tilde e_n\rangle)^2 - \sup_{P\in\mathcal{W}_N(\delta)} \sum_{n=1}^{N}\mathbb{E}_{P^{(N)}}(\langle b-\phi_N^\star,\tilde e_n\rangle)^2\notag\\
 & \qquad + \sup_{P\in\mathcal{W}(\delta) }\sum_{n=N+1}^\infty\mathbb{E}_{P}(\langle b,\tilde e_n\rangle)^2 \label{eq:nogapsplit}.
\end{align}\ssw{the infinite series of $\tilde e_n$ is not even defined.}
For the last term in~\eqref{eq:nogapsplit}, note that
\begin{align*}
&\sup_{P\in\mathcal{W}(\delta)}\mathbb{E}_{P}\left[\sum_{n=N+1}^\infty(\langle b,\tilde e_n\rangle)^2\right]\\
&\leq 2\mathbb{E}_{P_0}\left[\sum_{n=N+1}^\infty(\langle b^0,\tilde e_n\rangle)^2\right]+2\sup_{P\in\mathcal{W}(\delta)}\mathbb{E}_{\pi}\left[\sum_{n=N+1}^\infty(\langle b-b^0,\tilde e_n\rangle)^2\right]\\
& \leq o(1) + O(1)\sup_{P\in\mathcal{W}(\delta)}\mathbb{E}_\pi\left[\|b-b^0-\sum_{n=1}^N\langle b-b^0,e_n\rangle e_n\|^2_{L^2(\mathcal{D})}\right],
\end{align*}
where $\pi$ is the optimal coupling between the marginals $P_0$ and $P$, and \ssw{Need to assume that $\mathcal H_{\tilde w}$ embeds continuously into $\mathcal C(\mathcal D)$, make precise in Assumption 2.3. Say we use assumption 2.3}
\begin{align*}
&  \sup_{P\in\mathcal{W}(\delta)}\mathbb{E}_\pi\left[\|b-b^0-\sum_{n=1}^N\langle b-b^0,e_n\rangle e_n\|^2_{L^2(\mathcal{D})}\right]\\
& \leq O(1) \sup_{P\in\mathcal{W}(\delta)}\mathbb{E}_\pi\left[\|b-b^0-\sum_{n=1}^N\langle b-b^0,e_n\rangle e_n\|^2_{\mathcal{H}_{\tilde w}}\right]\\
&\leq o(1)\sup_{P\in\mathcal{W}(\delta)}\mathbb{E}_\pi\left[\|b-b^0-\sum_{n=1}^N\langle b-b^0,e_n\rangle e_n\|^2_{\mathcal{H}_w}\right]\\
&\to 0\textrm{ as } N\to\infty,
\end{align*}
% \begin{align*}
% &  \sup_{P\in\mathcal{W}(\delta)}\mathbb{E}_\pi\left[\|T(b-b^0-\sum_{n=1}^N\langle b-b^0,e_n\rangle e_n)\|^2_{L^2(\mathcal{D})}\right]\\
% & \leq O(1) \sup_{P\in\mathcal{W}(\delta)}\mathbb{E}_\pi\left[\|T(b-b^0-\sum_{n=1}^N\langle b-b^0,e_n\rangle e_n)\|^2_{\mathcal{H}_{\tilde w}}\right]\\
% & \leq O(1) \sup_{P\in\mathcal{W}(\delta)}\mathbb{E}_\pi\left[\|b-b^0-\sum_{n=1}^N\langle b-b^0,e_n\rangle e_n\|^2_{\mathcal{H}_{\tilde w}}\right]\\
% &\leq o(1)\sup_{P\in\mathcal{W}(\delta)}\mathbb{E}_\pi\left[\|b-b^0-\sum_{n=1}^N\langle b-b^0,e_n\rangle e_n\|^2_{\mathcal{H}_w}\right]\\
% &\to 0\textrm{ as } N\to\infty,
% \end{align*}
where the second inequality comes from  Assumption~\ref{assmp:operator2} \ssw{Just quote Assumption 2.4!!!!} and the fact that\ssw{'the fact that'. We need to use correct grammar!!} $b-b^0-\sum_{n=1}^N\langle b-b^0,e_n\rangle e_n$ is in the space spanned by $\{e_n\}_{n=N+1}^\infty$. 
For the first two terms in~\eqref{eq:nogapsplit}, we write
\[
\phi^\star_N(Y_1,\ldots,Y_m) = \sum_{j=1}^m\alpha_j^\star\left(\sum_{k=1}^N\langle b,e_k\rangle T(e_k)(x_j)+\epsilon_j\right) + \sum_{j=1}^m\alpha_j^\star R_{j,N}.
\]
Thus for any feasible $P$ in $\mathcal{W}(\delta)$,
we have
\begin{align*}
&\mathbb{E}_P(\langle b-\phi^\star_N,\tilde e_n\rangle)^2\\
& = 
  \mathbb{E}_P\left(\left\langle b - \sum_{j=1}^m\alpha_j^\star(x)\left(\sum_{k=1}^N\langle b,e_k\rangle T(e_k)(x_j)+\epsilon_j\right),\tilde e_n\right\rangle^2\right)\\
  & -2\mathbb{E}_P\left(\left\langle b - \sum_{j=1}^m\alpha_j^\star(x)\left(\sum_{k=1}^N\langle b,e_k\rangle T(e_k)(x_j)+\epsilon_j\right),\tilde e_n\right\rangle\left\langle \sum_{j=1}^m\alpha_j^\star(x)R_{j,N},\tilde e_n\right\rangle\right)\\
  & + \mathbb{E}_P\left(\left\langle\sum_{j=1}^m\alpha_j^\star(x)R_{j,N},\tilde e_n\right\rangle^2\right).
\end{align*}
Note that
\begin{align*}
 & \sum_{n=1}^{N}\mathbb{E}_P\left(\left\langle b - \sum_{j=1}^m\alpha_j^\star(x)\left(\sum_{k=1}^N\langle b,e_k\rangle T(e_k)(x_j)+\epsilon_j\right),\tilde e_n\right\rangle^2\right)\\
 & \leq \sup_{P\in\mathcal{W}_N(\delta)} \sum_{n=1}^{N} \mathbb{E}_{P^{(N)}}(\langle b-\phi_N^\star,\tilde e_n\rangle)^2.
%  &\leq \sup_{P:\W(P,P_0)\leq \delta,\supp(P|_{\mathcal{H}})\subseteq \Span\{T(e_n)\}_{n=1}^N} \sum_{n=1}^N \mathbb{E}_{P}(\langle u-\phi_N^\star,e_n\rangle)^2 + \lambda_N^\star \mathbb{E}_{P_0}\left[\sum_{n=N+1}^\infty(\langle b,e_n\rangle)^2w_n\right]\\
%  &\leq \sup_{P:\W(P,P_0)\leq \delta,\supp(P|_{\mathcal{H}})\subseteq \Span\{T(e_n)\}_{n=1}^N} \sum_{n=1}^N \mathbb{E}_{P}(\langle u-\phi_N^\star,e_n\rangle)^2 + O(1) \mathbb{E}_{P_0}\left[\sum_{n=N+1}^\infty(\langle b,e_n\rangle)^2w_n\right]
 \end{align*}
% where $\lambda_N^\star$ is the dual optimal in~\eqref{eq:subprimal} under $\phi^\star_N$. 
Also, by Cauchy-Schwarz, we obtain
\begin{align*}
& \sum_{n=1}^{N}  \mathbb{E}_P\left(\left\langle b - \sum_{j=1}^m\alpha_j^\star(x)\left(\sum_{k=1}^N\langle b,e_k\rangle T(e_k)(x_j)+\epsilon_i\right),\tilde e_n\right\rangle\left\langle \sum_{j=1}^m\alpha_j^\star(x)R_{j,N},\tilde e_n\right\rangle\right)\\
&\leq \sqrt{\mathbb{E}_P\left(\sum_{n=1}^{N}\left\langle b - \sum_{j=1}^m\alpha_j^\star(x)\left(\sum_{k=1}^N\langle b,e_k\rangle T(e_k)(x_j)+\epsilon_i\right),\tilde e_n\right\rangle^2\right)}\cdot\\
&\quad\sqrt{\mathbb{E}_P\left(\sum_{n=1}^{N}\left\langle \sum_{j=1}^m\alpha_j^\star(x)R_{j,N},\tilde e_n\right\rangle^2\right)}\\
&\leq\sqrt{ \sup_{P\in\mathcal{W}_N(\delta)} \sum_{n=1}^{N} \mathbb{E}_{P^{(N)}}(\langle b-\phi_N^\star,\tilde e_n \rangle)^2}\sqrt{mC\left(\sum_{j=1}^m\mathbb{E}_P\left[R_{j,N}^2\right]\right)}.
\end{align*}
Therefore, denoting $\pi$ as the optimal coupling between $P$ and $P_0$, we have
\begin{align*}
& \sup_{P\in\mathcal{W}(\delta)}\mathbb{E}_P\left(R_{j,N}\right)^2\notag\\
& = \sup_{P\in\mathcal{W}(\delta)}\mathbb{E}_P\left (u(x_j) - \sum_{k=1}^N \langle b,e_k\rangle T(e_k)(x_j)\right)^2\notag\\
& \leq 2\mathbb{E}_{P_0}\left (u^0(x_j) - \sum_{k=1}^N \langle b^0,e_k\rangle T(e_k)(x_j)\right)^2 \\
&\qquad + 2\sup_{P\in\mathcal{W}(\delta)}\mathbb{E}_{\pi}\left(u(x_j)-u^0(x_j) - \sum_{k=1}^N \langle b-b^0,e_k\rangle T(e_k)(x_j)\right)^2\\
&\leq o(1) + 2\sup_{P\in\mathcal{W}(\delta)}\mathbb{E}_{\pi}\left(u(x_j)-u^0(x_j) - \sum_{k=1}^N \langle b-b^0,e_k\rangle T(e_k)(x_j)\right)^2.
\end{align*}
Moreover, we also find
\begin{align*}
    & \sup_{P\in\mathcal{W}(\delta)}\mathbb{E}_{\pi}\left(u(x_j)-u^0(x_j) - \sum_{k=1}^N \langle b-b^0,e_k\rangle T(e_k)(x_j)\right)^2\\
    &=   \sup_{P\in\mathcal{W}(\delta)}\mathbb{E}_{\pi}\left[\left(T(b-b^0-\sum_{k=1}^N\langle b-b^0,e_k\rangle e_k)(x_j)\right)^2\right]\\
    & \leq O(1) \sup_{P\in\mathcal{W}(\delta)}\mathbb{E}_{\pi}\left[\|T(b-b^0-\sum_{k=1}^N\langle b-b^0,e_k\rangle e_k)\|_{\mathcal{H}_{\tilde w}}^2\right]\\
    & \leq O(1) \sup_{P\in\mathcal{W}(\delta)}\mathbb{E}_{\pi}\left[\|b-b^0-\sum_{k=1}^N\langle b-b^0,e_k\rangle e_k\|_{\mathcal{H}_{\tilde w}}^2\right]\\
    & \leq o(1) \sup_{P\in\mathcal{W}(\delta)}\mathbb{E}_{\pi}\left[\|b-b^0-\sum_{k=1}^N\langle b-b^0,e_k\rangle e_k\|_{\mathcal{H}_w}^2\right]\\
    &\leq o(1) \sup_{P\in\mathcal{W}(\delta)}\mathbb{E}_{\pi}\left[\|b-b^0\|_{\mathcal{H}_w}^2\right]\\
    &\to 0\textrm{ as }N\to\infty,
\end{align*}
where the first inequality follows because $\mathcal{H}_{\tilde w}$ is an RKHS consisting of continuous functions, so that point evaluations at the $x_j$'s are bounded, according to Assumption~\ref{assmp:operator}(ii). The second inequality also follows from Assumption~\ref{assmp:operator}(ii), while the third equality follows from Assumption~\ref{assmp:operator2}. Finally, the last inequality follows because $\{e_n\}_{n=1}^\infty$ is an orthogonal system.
Therefore, we have that
\[
\sup_{P\in\mathcal{W}(\delta)}\mathbb{E}_P\left(R_{j,N}\right)^2 \to0 \textrm{ as }N\to\infty.
\]
It then follows that
\[
\sup_{P\in\mathcal{W}(\delta)} \Obj(\phi^\star_N,P) -\sup_{P\in\mathcal{W}_N(\delta)}\Obj(\phi_N^\star,P^{(N)})= o(1)\textrm{ as }N\to\infty,
\]
establishing our Claim 3.
\end{proof}

\subsection{Proof of Proposition~\ref{prop:fdvalueconverges}}
\begin{proof}[Proof of Proposition~\ref{prop:fdvalueconverges}]
This is an immediate consequence of the inequality~\eqref{eq:sandwich} and Claims 1-3 in the proof of Theorem~\ref{thm:swap}.
\end{proof}

\subsection{Proof of Proposition~\ref{prop:weakcompact}}

\begin{proof}[Proof of Proposition~\ref{prop:weakcompact}]
Since the space $C(\mathcal{D})\times\mathbb{R}^m$ is Polish, if suffices to show that the sequence  $P_N,N=1,\ldots,\infty$ is tight. First note that $P_0$ is a finite Radon measure on $C(\mathcal{D})\times\mathbb{R}^m$, and hence $P_0$ is tight. Therefore, for any $0<\eta<1$, there exists a compact set $\mathcal{C}_0$ in $C(\mathcal{D})\times\mathbb{R}^m$, such that $P_0(\mathcal{C}_0)\geq 1-\frac{\eta}{2}$. Now consider the set
\[
\mathcal{C}_1 = \{(b,\epsilon)\in C(\mathcal{D})\times\mathbb{R}^m: \exists~ (b^0,\epsilon^0)\in \mathcal{C}_0,~\|b-b^0\|_{\mathcal{H}_w}^2+\|\epsilon-\epsilon_0\|_2^2\leq 2\delta^2/\eta\},
\]
which is compact since the space $\mathcal{H}_w$ is compactly embedded in $C(\mathcal{D})$ by  Assumption~\ref{assmp:perturbation}. Denote $Q_N$ as the law of the difference $(b-b^0,\epsilon-\epsilon^0)$ under the optimal coupling between $(b,\epsilon)\sim P_N$ and $(b^0,\epsilon^0)\sim P_0$. Then, 
\begin{align*}
1-P_N(\mathcal{C}_1) & \leq 1-P_0(\mathcal{C}_0) + Q_N(\{(b,\epsilon)\in\mathcal{H}_w\times\mathbb{R}^m: \|b\|^2_{\mathcal{H}_w}+\|\epsilon\|_2^2 > 2\delta^2/\eta\})\\
& \leq \frac{\eta}{2} + \frac{\W^2(P_N,P_0)}{2\delta^2/\eta} \leq \frac{\eta}{2} + \frac{\eta}{2}=\leq\eta,
\end{align*}
where in the second inequality we used the Cauchy-Schwarz inequality. Therefore, we conclude that
\[
P_N(\mathcal{C}_1)\geq 1-\eta~~~\forall N,
\]
and hence that the sequence $P_N$ is tight. It follows that there exists a weakly convergent subsequence $P_{N_l},l\geq1$ with $P_{N_l}\Rightarrow P_\infty$. Since the Wasserstein distance is lower semi-continuous, the limit $P_\infty$ necessarily satisfies $\W(P_\infty,P_0)\leq\delta$.

\ssw{I think we also use here that the Wasserstein ball is closed for the space of probability measures equipped with the weak topology, right?}
\end{proof}

\subsection{Proof of Theorem~\ref{thm:existnash}}
We start with a simple but useful result, which is well-known~\cite{ref:von2004theory}. \ssw{If it is well-known, cite where it is to be found the first time! }
\begin{lemma}[Minimax theorem and Nash equilibrium]\label{thm:nash}
Let $\phi^\star$ be an optimal solution to the outer infimum of the problem 
\[
\inf_{\phi\in\mathcal{M}}\sup_{P\in\mathcal{W}(\delta)}\Obj(\phi,P),
\]
and let $P^\star$ be an optimal solution to the outer supremum of the problem
\[
 \sup_{P\in\mathcal{W}(\delta)}\inf_{\phi\in\mathcal{M}}\Obj(\phi,P).
\]
Then $(\phi^\star,P^\star)$ is a Nash equilibrium of~\eqref{eq:thmstrongdualityregressappendix}.
\end{lemma}

\begin{proof}[Proof of Lemma~\ref{thm:nash}]
We have
\begin{align*}
 \sup_{P\in\mathcal{W}(\delta)}\inf_{\phi \in \mathcal M}\Obj(\phi,P)&=\inf_{\phi \in \mathcal M}\Obj(\phi,P^\star) \leq \Obj(\phi^\star,P^\star)\\
& \leq \sup_{P\in\mathcal{W}(\delta)}\Obj(\phi^\star,P) = \inf_{\phi \in \mathcal M}\sup_{P\in\mathcal{W}(\delta)}\Obj(\phi,P).
\end{align*}
By strong duality, all inequalities become equalities. Thus, we have
\[
\Obj(\phi^\star,P^\star) = \min_{\phi \in \mathcal M}\Obj(\phi,P^\star) = \max_{P\in\mathcal{W}(\delta)}\Obj(\phi^\star,P).
\]
This completes the proof.
\end{proof} 

\begin{proof}[Proof of Theorem~\ref{thm:existnash}]
Recall that we denote
 \[
 Q_N^{\star} \coloneqq \mathcal L^{P^\star_N} \big( \langle b,e_1\rangle,\cdots,\langle b,e_N\rangle,\epsilon\big)
 \]
in the proof to Claim 2 of Theorem~\ref{thm:swap}. We also denote additionally
 \[P_\star^{(N)} \coloneqq \mathcal L^{P^\star_N} \Big( \big(\sum_{1\le n\le N} e_n \langle b,e_n\rangle, \epsilon \big) \Big).\]
\label{nomencl:PlsuN}\nomenclature{$P_\star^{(N)}$}{truncation of $P^\star_N$ without the tails}{nomencl:PlsuN}Letting $k^{(N)}(x,x)=\mathbb{E}_{P_\star^{(N)}}[b(x)^2]$, the objective value $\Obj(\phi^\star_N,P_\star^{(N)})$ can be written as
\begin{align*}
& \Obj(\phi^\star_N,P_\star^{(N)})\\
&=\int_\mathcal{D} \mathbb{E}_{P_\star^{(N)}}[\mathbb{V}\mathrm{ar}_{P_\star^{(N)}}[b(x)|Y_1,\ldots,Y_m]] \dd x\\
& = \int_\mathcal{D}  k^{(N)}(x,x)\\
&\qquad\qquad-  \left(k^{(N)}_\epsilon(x,x_1),\ldots,k^{(N)}_\epsilon(x,x_m)\right) (K^{(N)}_\epsilon)^{-1} \left(k^{(N)}_\epsilon(x,x_1),\ldots,k^{(N)}_\epsilon(x,x_m)\right)^\top \dd x.
\end{align*}
 While the limiting objective value $\Obj(\phi^\star_\infty,P^\star_\infty)$ is
\begin{align*}
 &\Obj(\phi^\star_{\infty},P^\star_{\infty})\\ &=\int_\mathcal{D}\mathbb{E}_{P^\star_\infty}[\mathbb{V}\mathrm{ar}_{P^\star_\infty}[b(x)|Y_1,\ldots,Y_m]] \dd x\\
& = \int_\mathcal{D}  k_{P^\star_\infty}(x,x) -  \left(k_{\epsilon,P^\star_\infty}(x,x_1),\ldots,k_{\epsilon,P^\star_\infty}(x,x_m)\right)\cdot (K_{\epsilon,P^\star_\infty})^{-1} \cdot\left(k_{\epsilon,P^\star_\infty}(x,x_1),\ldots,k_{\epsilon,P^\star_\infty}(x,x_m)\right)^\top \dd x,
\end{align*}
where $k_{P^\star_\infty}(x,x)=\mathbb{E}_{P^\star_\infty}[b(x)^2]$. Recall that $P^\star_{N_l}$ is a (Gaussian) subsequence of $P^\star_N$, which converges weakly to $P^\star_\infty$. Note that $T$ is a bounded linear operator and that the tail-difference in $P_\star^{(N_l)}$ and  $P^\star_{N_l}$ is negligible. By~\cite[Exercise 2.1.4]{ref:gine2015mathematical}, convergence in distribution implies convergence of the second moment, thus
uniformly for $x\in\mathcal{D}$, $k^{(N_l)}(x,x)\to k_{P^\star_\infty}(x,x)$ and $k_\epsilon^{(N_l)}(x,x_i)\to k_{\epsilon,P^\star_\infty}(x,x_i)$. Also, it holds that $K_\epsilon^{(N_l)}\to K_{\epsilon,P^\star_\infty}$.
Thus, we have
\[
\Obj(\phi^\star_{N_l},P_\star^{(N_l)}) = \min_{\phi \in \mathcal M}\Obj(\phi,P_\star^{(N_l)}) \to \Obj(\phi^\star_\infty,P^\star_\infty) = \min_{\phi \in \mathcal M} \Obj(\phi,P^\star_\infty), 
\]
and hence that $P^\star_\infty$ solves the right-hand side of~\eqref{eq:thmstrongdualityregressappendix}.

%$P^\star_\infty$ is Gaussian by pointwise-convergence of $\mathcal{H}$, since the a.s. limit (actually everywhere)
%\[
%\sum_{...}\langle u,  e_n\rangle e_n(x)
%\]
%s Gaussian (a.s convergence implies weak convergence, weak limit of Gaussians is Gaussian). $\langle u,e_n\rangle$ for finite number of $n$ is jointly Gaussian by considering the marginals of $P^\star_N$.

We next show that  $\phi^\star_\infty$ solves the left-hand side of~\eqref{eq:thmstrongdualityregressappendix}. We first note that  due to the convergence of moments of $P^\star_{N_l}$, it holds that
\[\lim_{l\to\infty}\phi^\star_{N_l}(Y_1,\ldots,Y_m)(x) = \phi^\star_\infty(Y_1,\ldots,Y_m)(x),~~~\forall~ Y_1,\ldots,Y_m\in \R,~~~ x\in\mathcal{D},\]
% \[
% \phi^\star_{N_l}(\cdot,\ldots,\cdot)(\cdot)\to\phi^\star_\infty(\cdot,\ldots,\cdot)(\cdot)\textrm{ everywhere point-wise as }l\to\infty.
% \]
Next, by Fatou's lemma, we have that for any $P\in\mathcal{W}(\delta)$ ,
\[
\mathbb{E}_P\left[\|b- \phi^\star_{\infty}(Y_1,\ldots,Y_m)\|^2_{L^2(\mathcal{D})}\right]\leq\liminf_{l\to\infty} \mathbb{E}_P\left[\|b - \phi^\star_{N_l}(Y_1,\ldots,Y_m)\|^2_{L^2(\mathcal{D})}\right].
\]
Therefore,
\begin{align*}
\sup_{P\in\mathcal{W}(\delta)}\mathbb{E}_P\left[\|b - \phi^\star_{\infty}(Y_1,\ldots,Y_m)\|^2_{L^2(\mathcal{D})}\right]\leq\liminf_{l\to\infty} \sup_{P\in\mathcal{W}(\delta)}\mathbb{E}_P\left[\|b- \phi^\star_{N_l}(Y_1,\ldots,Y_m)\|^2_{L^2(\mathcal{D})}\right].
\end{align*}
Note that due to the finite-dimensional strong duality~\eqref{eq:finitedimduality}, $\phi^\star_N$ also solves
\[
\min_{\phi\in\mathcal{M}_{\textrm{aff},N} }\max_{P\in\mathcal{W}_N(\delta)} \Obj(\phi,P^{(N)}).
\]
By~\eqref{eq:claimtoprove2} and strong duality, we have
\[
\liminf_{N\to\infty} \sup_{P\in\mathcal{W}(\delta)}\mathbb{E}_P\left[\|b - \phi^\star_N(Y_1,\ldots,Y_m)\|^2_{L^2(\mathcal{D})}\right]  =\min_{\phi \in \mathcal M} \sup_{P\in\mathcal{W}(\delta)} \Obj(\phi,P)
\]
Thus the estimator $\phi^\star_\infty$ solves the min-max problem on the left-hand side of~\eqref{eq:thmstrongdualityregressappendix}.

% Since for any $\phi^\star$ that is a point-wise limit of $\phi^\star_N$ has to solve $\min_{\phi \in \mathcal M} \Obj(\phi,P^\star_\infty)$ by
By Lemma~\ref{thm:nash}, since $\phi^\star_\infty$ has to solve $\min_{\phi \in \mathcal M}\Obj(\phi,P^\star_\infty)$ , we have that $\phi^\star_\infty$ is uniquely determined (independent of the choice of the subsequence $\phi^\star_{N_l}$), thus\ssw{Replace `everywhere pointwise' by a correct phrase}
\[\lim_{N\to\infty}\phi^\star_{N}(Y_1,\ldots,Y_m)(x) = \phi^\star_\infty(Y_1,\ldots,Y_m)(x)~~~\forall~ Y_1,\ldots,Y_m\in \R,~~~ x\in\mathcal{D}.\]
% \[
% \phi^\star_{N}(\cdot,\ldots,\cdot)(\cdot)\to\phi^\star_\infty(\cdot,\ldots,\cdot)(\cdot)\textrm{ everywhere point-wise as }N\to\infty.
% \]
Moreover, by Lemma~\ref{thm:nash}, $P^\star_\infty$ has to solve $\max_{P\in\mathcal{W}(\delta)}\Obj(\phi^\star_\infty,P)$. We show in the sequel that $P^\star_\infty$ is uniquely determined (independent of the choice of the subsequence $P^\star_{N_l}$),\ssw{Clarify this sentence?? I.e. the sequence $P_N^\star$ converges as a whole?} from which the sequence $P_N^\star$ converges to $P^\star_\infty$ in the weak topology as a whole. 
% By Lemma~\ref{thm:nash}, any $P^\star$ has to solve $\max_{P\in\mathcal{W}(\delta)}\Obj(\phi^\star_\infty,P)$. 
By Theorem 1 in~\cite{ref:blanchet2019quantifying}, we have the reformulation
\begin{align*}
    & \max_{P\in\mathcal{W}(\delta)}\Obj(\phi^\star_\infty,P)\\
    = & \inf_{\gamma\geq0}\Big(\gamma\delta^2 + \\
    & \mathbb{E}_{P_0}\Big[\sup_{(b,\epsilon)\in C(\mathcal{D})\times\mathbb{R}^m}\left(\|b -\phi^\star_\infty((T(b)(x_i)+\epsilon_i)_i)\|^2_{L^2(\mathcal{D})} - \gamma (\|b-b^0\|^2_{\mathcal{H}_w}+\|\epsilon-\epsilon^0\|^2_2)\right)\Big]\Big).
\end{align*}
We claim that the optimal dual $\gamma^\star$ is sufficiently large so that for all $( b,\epsilon)\in\mathcal{H}_w\times\mathbb{R}^m$ such that $(b,\epsilon)\neq0$:
\begin{equation}\label{eq:claimlargegammainfinite}
\|b -\phi^\star_\infty((T(b)(x_i)+\epsilon_i)_i)\|^2_{L^2(\mathcal{D})} - \gamma^\star (\|b\|^2_{\mathcal{H}_w}+\|\epsilon\|^2_2)<0.
\end{equation}
It is easy to see that $\gamma^\star>0$. Suppose that~\eqref{eq:claimlargegammainfinite} does not hold, then
\begin{equation}\label{eq:rev_contra}
\|b^\star -\phi^\star_\infty( (T(b^\star)(x_i)+\epsilon^\star_i)_i)\|^2_{L^2(\mathcal{D})} - \gamma^\star (\|b^\star\|^2_{\mathcal{H}_w}+\|\epsilon^\star\|^2_2)\geq0,
\end{equation}
for some \review{$(b^\star,\epsilon^\star)\in\mathcal{H}_w\times\mathbb{R}^m$ and $(b^\star,\epsilon^\star)\neq0$}. Then for all $(b^0,\epsilon^0)\in C(\mathcal{D})\times\mathbb{R}^m$ satisfying 
\[
\langle b^\star -\phi^\star_\infty((T(b^\star)(x_i)+\epsilon^\star_i)_i),b^0 -\phi^\star_\infty( (T(b^0)(x_i)+\epsilon^0_i)_i)\rangle\neq0,
\]
we have that 
\begin{align*}
&\sup_{(b,\epsilon)\in C(\mathcal{D})\times\mathbb{R}^m}\left(\|b -\phi^\star_\infty((T(b)(x_i)+\epsilon_i)_i)\|^2_{L^2(\mathcal{D})} - \gamma^\star (\|b-b^0\|^2_{\mathcal{H}_w}+\|\epsilon-\epsilon^0\|^2_2)\right) \\
&\geq \sup_t \left(\|t\left(b^\star -\phi^\star_\infty((T(b^\star)(x_i)+\epsilon^\star_i)_i)\right)+\left(b^0 -\phi^\star_\infty( (T(b^0)(x_i)+\epsilon^0_i)_i)\right)\|^2_{L^2(\mathcal{D})} \right)\\
& \qquad - t^2\gamma^\star(\|b^\star\|_{\mathcal{H}_w}^2 + \|\epsilon^\star\|^2_2)\\
& \geq\sup_t t^2\left(\|b^\star -\phi^\star_\infty( (T(b^\star)(x_i)+\epsilon^\star_i)_i)\|^2_{L^2(\mathcal{D})} - \gamma^\star (\|b^\star\|^2_{\mathcal{H}_w}+\|\epsilon^\star\|^2_2)\right)\\
& \qquad+ 2t\langle b^\star -\phi^\star_\infty((T(b^\star)(x_i)+\epsilon^\star_i)_i),b^0 -\phi^\star_\infty( (T(b^0)(x_i)+\epsilon^0_i)_i)\rangle \\
& \qquad + \|b^0 -\phi^\star_\infty( (T(b^0)(x_i)+\epsilon^0_i)_i)\|^2_{L^2(\mathcal{D})}\\
& =\infty,
\end{align*}
where the infinite optimal value is due to the nonnegativity of the coefficient associated with the quadratic term $t^2$.
Moreover, observe that the constraint
\begin{equation}\label{eq:linearconstraintappendix}
\langle b^\star -\phi^\star_\infty((T(b^\star)(x_i)+\epsilon^\star_i)_i),b^0 -\phi^\star_\infty( (T(b^0)(x_i)+\epsilon^0_i)_i)\rangle=0
\end{equation}
is a linear constraint in $(b^0,\epsilon^0)$,  thus the collection of $(b^0,\epsilon^0)$ satisfying~\eqref{eq:linearconstraintappendix}, denoted by $\mathcal{C}_{(b^\star,\epsilon^\star)}$, is a linear subspace of $C(\mathcal{D})\times \mathbb{R}^m$. This collection is closed (with respect to the product norm $\|\cdot\|_{C(\mathcal{D})}\times \|\cdot\|_{\ell_2}$) because $T:C(\mathcal{D})\to C(\mathcal{D})$ is bounded. Note that the support of $P_0$, denoted by $\textrm{supp}(P_0)$, contains basis vectors of the form
\[
\bar{B} = \{(e_k, 0,\ldots, 0)\}_{k=1}^\infty  \cup \{(\mathbf{0}, \bar{e}_j)\}_{j=1}^m
\]
where $\{ e_k\}_{k=1}^{\infty}$ is the orthonormal basis of $C(\mathcal{D})$, $\mathbf{0}\in C(\mathcal{D})$ is the zero-function, and $\bar{e}_j\in \mathbb{R}^m$ denotes the j-th standard basis on $\mathbb{R}^m$. Suppose that $\bar{B}\subseteq \mathcal{C}_{(b^\star,\epsilon^\star)}$. Since $b^\star\in\mathcal{H}_{w}\subseteq \mathcal{H}_{\tilde w}$, by expanding $(b^\star,\epsilon^\star)$ in the basis vectors from $\bar{B}$ (which converges in the product norm $\|\cdot\|_{\mathcal{H}_{\tilde w}}\times \|\cdot\|_{\ell_2}$), and considering that $T: \mathcal{H}_{\tilde w}\to \mathcal{H}_{\tilde w}$ is bounded from Assumption~\ref{assmp:operator} (ii) as well as the fact that $\|\cdot\|_{\mathcal{H}_{\tilde w}}$ is a stronger norm than the $L^2(\mathcal{D})$-norm, we conclude that $(b^\star,\epsilon^\star) \in \mathcal{C}_{(b^\star,\epsilon^\star)}$, which is a contradiction to~\eqref{eq:rev_contra}. Hence there is at least one basis vector of $\bar{B}$ not in $\mathcal{C}_{(b^\star,\epsilon^\star)}$. Since $\textrm{supp}(P_0)$ is the smallest closed subset of $C(\mathcal{D})\times \mathbb{R}^m$ on which the probability measure evaluates to one, we have that $P_0(\mathcal{C}_{(b^\star,\epsilon^\star)}) < 1$. This concludes the proof of claim~\eqref{eq:claimlargegammainfinite}.

Denoting the shorthand
\[
J(b,\epsilon) = \|b -\phi^\star_\infty((T(b)(x_i)+\epsilon_i)_i)\|^2_{L^2(\mathcal{D})} - \gamma^\star (\|b-b^0\|^2_{\mathcal{H}_w}+\|\epsilon-\epsilon^0\|^2_2),
\]
% By calculus of variation\ssw{What does this mean? Let's just leave out this phrase}
\ssw{The notation $\Delta b,\Delta \eps$ doesn't make sense in light of the Laplacian. Could use e.g. $h_b$, $h_\eps$ as notation of the displacement}any maximizer of $J$ has to satisfy the first-order optimality condition
\[
\frac{\partial}{\partial t}J(b+th_ b,\epsilon+th_\epsilon)\Bigg|_{t=0} = 0\quad \forall (h_b,h_\epsilon)\in  C(\mathcal{D})\times\mathbb{R}^m.
\]
It is easy to compute that
\begin{align*}
&\frac{1}{2}\frac{\partial}{\partial t}J(b+th_ b,\epsilon+th_\epsilon)\Bigg|_{t=0}\\
= & \left\langle b-\phi^\star_\infty( (T(b)(x_i)+\epsilon_i)_i), h_ b-\phi^\star_\infty((T(h_ b)(x_i)+(h_\epsilon)_i)_i)\right\rangle \\
& -\gamma^\star \left(\langle b-b^0,h_ b\rangle_{\mathcal{H}_w}+\langle \epsilon-\epsilon^0,h_\epsilon\rangle\right).
\end{align*}
Suppose that $(b,\epsilon)$ and $(\tilde b,\tilde\epsilon)$ are two maximizers of $J$, by choosing $(h_ b,h_\epsilon) = (b-\tilde b,\epsilon-\tilde\epsilon)$, we have
\[
\|h_ b -\phi^\star_\infty((T(h_ b)(x_i)+(h_\epsilon)_i)_i)\|^2_{L^2(\mathcal{D})} - \gamma^\star (\|h_ b\|^2_{H_w}+\|h_\epsilon\|^2_2) = 0.
\]
From our earlier claim~\eqref{eq:claimlargegammainfinite}, we conclude that $(b,\epsilon)=(\tilde b,\tilde\epsilon)$. By Theorem 1 in~\cite{ref:blanchet2019quantifying}, this shows that the solution to  $\max_{P\in\mathcal{W}(\delta)}\Obj(\phi^\star_\infty,P)$ is unique.
\end{proof}
%\ssw{we haven't spoken of the primal optimal transport plan, please explain to the reader what it is!! You can't assume they will understand} 

\subsection{Proof of Theorem~\ref{thm:nashproperty_revision}}\review{

\begin{proof}[Proof of Theorem~\ref{prop:nasheqparity}]

Since $P^\star_N \Rightarrow P^\star_\infty$ as $N\to\infty$, and each $P^\star_N$ is Gaussian, it follows directly that $P^\star_\infty$ is also a Gaussian process. Since $\phi^\star_\infty$ has to solve $\Obj(\phi^\star_\infty,P^\star_\infty) = \min_{\phi \in \mathcal M} \Obj(\phi,P^\star_\infty)$, $\phi^\star_\infty$ is the conditional mean of the worst-case distribution $P^\star_\infty$, which is affine in the
observations.

\end{proof}
}

\subsection{Proof of Proposition~\ref{prop:nasheqparity}}
\begin{proof}[Proof of Proposition~\ref{prop:nasheqparity}]
First, suppose that the following statement holds true:
\begin{equation}\label{eq:conditionprop}
\textrm{``$K_{\epsilon,P^\star_\infty}$ is invertible for the limit $P^\star_\infty$ of some weakly convergent subsequence of $P^\star_N$''}.
\end{equation}
Then by Theorem~\ref{thm:existnash}, the sequence $P_N^\star$ weakly converges to $P^\star_\infty$ as a whole. \ssw{Complete the sentence.} Since the sequence $P^\star_N$ is Gaussian, we have $\textrm{det}(K^{(N)}_\epsilon)\to\textrm{det}(K_{\epsilon,P^\star_\infty})\neq0$ as $N\to\infty$.

Next, suppose that condition~\eqref{eq:conditionprop} does not hold,
% Next, suppose $K_\epsilon$ is not invertible (equivalently, $\textrm{det}(K_\epsilon)=0$) for the limit $P^\star_\infty$ of every weakly convergent subsequence of $P^\star_N$
\ssw{This is a really bad way to denote $P_\infty^\star$ -- haven't we already established that this is unique and that the whole sequence converges? Then we can just call it the limit!} and suppose that $\textrm{det}(K_\epsilon^{(N)})$ does not converge to $0$ as $N\to\infty$. Then there exists a subsequence of $P_N^\star$, denoted by $P^\star_{N_l},l\geq1$, such that the sequence of determinants $\textrm{det}(K^{(N_l)}_\epsilon)$ is uniformly bounded away from $0$. Furthermore, $P^\star_{N_l}$ has a weak limit (upon passing into a further subsequence, e.g., by Proposition~\ref{prop:weakcompact}). It is clear that condition~\eqref{eq:conditionprop} applies to this subsequence $P^\star_{N_l}$, which is a contradiction!
\end{proof}

\subsection{Proof of Lemma~\ref{lem:smalldelta}}
\begin{proof}[Proof of Lemma~\ref{lem:smalldelta}]
Let $P$ be such that $\W(P,P_0)\leq\delta$. To show that the matrix in question is strictly positive definite, it suffices to show that 
\[
\inf_{\xi\in\mathbb{R}^d:\|\xi\|_2=1} \xi^\top \big(\mathbb E_P[Y_iY_j]\big)_{ij} \xi = \inf_{\xi\in\mathbb{R}^d:\|\xi\|_2=1}\mathbb{E}_P\left[\left(\sum_{i=1}^m \xi_i(u(x_i)+\epsilon_i)\right)^2\right]> 0.
\]
To this end, let $\pi$ be a coupling such that 
\[
\mathbb{E}_{\pi}\left[c((b,\epsilon),(b^0,\epsilon^0))\right] = \mathbb{E}_{\pi}\left[\|\epsilon-\epsilon^0\|^2_2 + \|b- b^0\|^2_{\mathcal{H}_w}\right]\leq\delta^2,
\]
where the marginal distribution of $(b,\epsilon)$ and $( b^0,\epsilon^0)$ are $P$ and $P_0$, respectively. We denote $u=T(b)$ and $u^0=T(b^0)$. Let $C_m$ be the constant such that (due to the RKHS property of $\mathcal{H}_{\tilde w}$)\ssw{Why is this constant $C_m$ and not just $C$? Is the $m$-dependence important?{\color{red}:xz: $C_m$ is not bounded as $m\to\infty$, otherwise we conclude $\|\cdot\|_\infty\leq\|\cdot\|_{\tilde w}$}}
\[
|f(x_i)|\leq C_m \|f\|_{\mathcal{H}_{\tilde w}}\qquad\forall f\in\mathcal{H}_{\tilde w},~i=1,\ldots,m.
\]
By Fatou's lemma\ssw{Say that $C_w$ comes from assumption bla on $T$. In line 5, could we not strengthen this to small o?{\color{red}:xz:here the index $k$ runs from $1$ to $N$, if the index $k$ runs from $N$ to $\infty$ then we could strengthen to small o}} \ssw{If it's big-O, I wouldn't use the $O(1)$ notation here. It's just a universal constant like the other ones, isn't it??}
\begin{align}
\mathbb{E}_\pi\left[(u(x_i)-u^0(x_i))^2\right]&\leq \liminf_{N\to\infty} \mathbb{E}_\pi\left[\left(\sum_{k=1}^N\langle b-b^0,e_k\rangle T(e_k)(x_i)\right)^2\right]\notag\\
& = \liminf_{N\to\infty} \mathbb{E}_\pi\left[\left(T\left(\sum_{k=1}^N\langle b-b^0,e_k\rangle e_k\right)(x_i)\right)^2\right]\notag\\
&\leq C_m^2\liminf_{N\to\infty} \mathbb{E}_\pi\left[\left\|T(\sum_{k=1}^N\langle b-b^0,e_k\rangle e_k)\right\|^2_{\mathcal{H}_{\tilde w}}\right]\notag\\
&\leq C_m^2C_{\tilde w}^2\liminf_{N\to\infty} \mathbb{E}_\pi\left[\left\|\sum_{k=1}^N\langle b-b^0,e_k\rangle e_k\right\|^2_{\mathcal{H}_{\tilde w}}\right]\notag\\
&\leq C_m^2C_{\tilde w}^2 \tilde C\liminf_{N\to\infty}\mathbb{E}_\pi\left[\left\|\sum_{k=1}^N\langle b-b^0,e_k\rangle e_k\right\|^2_{\mathcal{H}_w}\right]\notag\\
&\leq  C_m^2C_{\tilde w}^2 \tilde C \mathbb{E}_\pi\left[\left\|b-b^0\right\|^2_{\mathcal{H}_w}\right]\label{eq:sufficondiineq},
\end{align}
where the constant $C_{\tilde w}$ is from Assumption~\ref{assmp:operator}(ii), and the constant $\tilde C$ is due to Assumption~\ref{assmp:operator2}.

For any $\xi\in\mathbb{R}^d$ with $\|\xi\|_2=1$, 
\begin{align*}
&\mathbb{E}_\pi\left[\left(\sum_{i=1}^m \xi_i( u^0(x_i)+\epsilon^0_i)\right)^2\right]\\
& = \mathbb{E}_{\pi}\left[\left(\sum_{i=1}^m \xi_i(u(x_i)+\epsilon_i) + \sum_{i=1}^m \xi_i( u^0(x_i)-u(x_i) +\epsilon^0_i-\epsilon_i)\right)^2\right]\\
& \leq 2\mathbb{E}_\pi\left[\left(\sum_{i=1}^m \xi_i(u(x_i)+\epsilon_i)\right)^2\right] + 4 \mathbb{E}_{\pi}\left[\left(\sum_{i=1}^m \xi_i(u^0 (x_i)-u(x_i))\right)^2\right] \\
&\qquad +4 \mathbb{E}_{\pi}\left[\left(\sum_{i=1}^m \xi_i(\epsilon^0_i-\epsilon_i)\right)^2\right]\\
& \leq 2\mathbb{E}_P\left[\left(\sum_{i=1}^m \xi_i(u(x_i)+\epsilon_i)\right)^2\right] +4 \mathbb{E}_{\pi}\left[\sum_{i=1}^m ( u^0(x_i)-u(x_i))^2\right] +4\mathbb{E}_{\pi} [\|\epsilon^0-\epsilon\|_2^2]\\
& \leq 2\mathbb{E}_P\left[\left(\sum_{i=1}^m \xi_i(u(x_i)+\epsilon_i)\right)^2\right]  +O(1)\delta^2, 
\end{align*}
where we used the previous inequality~\eqref{eq:sufficondiineq} to estimate $\mathbb{E}_{\pi}\left[\sum_{i=1}^m ( u^0(x_i)-u(x_i))^2\right]$. Note that $
 \left(\mathbb{E}_{P_0}[Y_iY_j]\right)_{ij}\succeq\sigma^2 I_{m\times m}$, since $\epsilon_i^0$ are independent $\mathcal{N}(0,\sigma^2)$ noise under $P_0$. Therefore,
\[
\inf_{\xi\in\mathbb{R}^d,\|\xi\|_2=1}\mathbb{E}_P\left[\left(\sum_{i=1}^m \xi_i(u(x_i)+\epsilon_i)\right)^2\right]\geq \frac{1}{2}\sigma^2 -O(1)\delta^2>0,
\]
for all $\delta<\delta_0$, where $\delta_0$ is a positive constant that depends on $T,m,(x_i)_i,\mathcal{H}_w,\mathcal{H}_{\tilde w}$ and $\sigma^2$. This completes the proof.
\end{proof}

\subsection{Proof of Theorem~\ref{thm:swapinverse}}\label{sec:regressionproof_revision}
The proof for the inverse problem works verbatim as the proof to Theorem~\ref{thm:swap} with minor modifications. \review{Instead of considering $\coef(\phi)\in\Span\{e_n\}_{n=1}^N$, we consider $\coef(\phi)\in\Span\{T(e_n)\}_{n=1}^N$. Additionally, in the proof to Theorem~\ref{thm:swap} Claim 2, we modify $\{\tilde e_k\}_{k=1}^{\infty}$ to be the Gram-Schmidt orthonormalization of $\{T(e_k)\}_{k=1}^\infty$\ssw{Is this the finite-dimensional space or the infinite dimensional one??}. In the process of orthonormalizing $\{T(e_k)\}_{k=1}^\infty$, for any $n$ such that the function $T(e_n)$ is already included in the span of $\{T(e_k)\}_{k=1}^{n-1}$, we define $\tilde e_n$ as the zero function.}
%Note that the functions $\{\tilde e_k\}_{k=1}^{\infty}$ do not necessarily constitute a basis of $L_2(\mathcal{D})$.}

\subsection{Details of Frank-Wolfe Algorithm}\label{sec:algodetail_revision}\review{

To numerically compute the Nash equilibrium of the min-max problem in the finite-dimensional approximation, we borrow the computational technique developed in \cite{ref:shafieezadeh2018wasserstein}. In particular, we note that the finite-dimensional problem~\eqref{eq:finitedimdualitymain} has the same structure as~\cite[Equation (2)]{ref:shafieezadeh2018wasserstein}. To see this, denoting
\[
Y = (Y_1,\ldots,Y_m),
\]
we can rewrite problem~\eqref{eq:finitedimdualitymain} as 
\begin{subequations}
\begin{align}
&\min_{\phi\in\mathcal{M}} \max_{Q^{(N)}:\W_N(Q^{(N)},Q_0^{(N)})\le \delta}~\Obj(\phi,Q^{(N)})\label{eq:computeleft}\\
& \qquad = \min_{\psi\in\mathcal{M}_N}\max_{Q^{(N)}:\W_N(Q^{(N)},Q_0^{(N)})\le \delta} \mathbb{E}_{Q^{(N)}}[\sum_{n=1}^N(\langle b, e_n \rangle  - \psi_n(Y))^2]\label{eq:computeright},
\end{align}
\end{subequations}
where $\mathcal{M}_N$ denotes the family of all measurable functions from $\mathbb{R}^m$ to $\mathbb{R}^N$\label{nomencl:MNfi}\nomenclature{$\mathcal{M}_N$}{measurable functions from $\mathbb{R}^m$ to $\mathbb{R}^N$}{nomencl:MNfi}, and $\psi = (\psi_1,\ldots,\psi_N)$. We note that the vector $(\{\langle b,e_n\rangle\}_{n=1}^N,Y)$ can be obtained by the matrix transform
\[
(\{\langle b,e_n\rangle\}_{n=1}^N,Y) = A \cdot (\{\langle b,e_n\rangle\}_{n=1}^N,\epsilon),
\]
where $\cdot$ denotes matrix multiplication, and
\[
A = \begin{pmatrix} I_{N} & 0 \\ A_{21} & I_m \end{pmatrix},\quad (A_{21})_{jn} = T(e_n)(x_j)~\forall 1\leq n \leq N, 1\leq j \leq m.
\]
% \[
% A = \begin{pmatrix} A_{11} & 0 \\ A_{21} & I\end{pmatrix}, A_{11} = diag(\lambda_1,\ldots,\lambda_N), (A_{21})_{ij} = \lambda_je_j(x_i).
% \]
Hence the covariance of $(\{\langle b,e_n\rangle\}_{n=1}^N,Y)$ under the nominal measure is given by
\[
A \Sigma_0 A^\top,
\]
where
\begin{align*}
\Sigma_0 &= \mathbb{C}\text{ov}_{Q_0^{(N)}}\left((\{\langle b,e_n\rangle\}_{n=1}^N,\epsilon)\right).
\end{align*}
By~\cite[Theorem 2.5]{ref:shafieezadeh2018wasserstein} and its proof,  the finite-dimensional problem~\eqref{eq:computeright} is thus equivalent to the finite convex program
\begin{align}
\sup &\quad \textrm{tr}\left[S_{11} - S_{12} S_{22}^{-1}S_{21}\right]\notag\\
\textrm{s.t. } &\quad S = \begin{pmatrix} S_{11} & S_{12} \\ S_{21} & S_{22}\end{pmatrix} \in \PD^{N+m},\quad S_{11}\in\PD^{N}, \quad S_{22} \in \PD^m,\quad S_{12} = S_{21}^\top \in\mathbb{R}^{N\times m},\notag\\
& \quad \textrm{tr}\left[BSB^\top+\tilde\Sigma - 2\left(\tilde\Sigma^{1/2} BSB^\top\tilde\Sigma^{1/2}\right)^{1/2}\right]\leq\delta^2, \quad S  \succeq \overline\sigma I_{N+m}.\label{eq:finiteconvex_revision}
\end{align}
where $B=\sqrt{W_N}A^{-1},\tilde \Sigma = \sqrt{W_N}\Sigma_0\sqrt{W_N}$, and $\overline\sigma$ is the minimum eigenvalue of $A \Sigma_0 A^\top$.
% \begin{align*}
% \Sigma &= Cov_{Q_0^{(N)}}\left((\langle b, e_1\rangle, \langle b, e_2\rangle,...,\langle b, e_N\rangle, \epsilon_1,\epsilon_2,\ldots, \epsilon_m)\right).
% \end{align*}
% and
% \[
% A = \begin{pmatrix} A_{11} & 0 \\ A_{21} & I\end{pmatrix}, A_{11} = diag(\lambda_1,\ldots,\lambda_N), (A_{21})_{ij} = \lambda_je_j(x_i).
% \]
Assuming the finite convex program has optimal solutions $S^\star,S^\star_{11},S^\star_{22}$ and $S^\star_{12}$. Then the optimal solution $\psi^\star$ to equation~\eqref{eq:computeright} is $\psi^\star(Y) = S^\star_{12}(S^\star_{22})^{-1}Y$, and hence the optimal solution $\phi^\star$ to equation~\eqref{eq:computeleft} is $\phi^\star = \sum_{n=1}^N \psi_n^*(Y) e_n$. The worst-case prior on $(\{\langle b, e_n \rangle\}_{n=1}^N)$ is given by a finite-dimensional Gaussian with zero-mean and covariance matrix $S^*_{11}$ and the worst-case posterior on $(\{\langle b, e_n \rangle\}_{n=1}^N)$ (conditional on $Y)$ is given by a finite-dimensional Gaussian with zero-mean and covariance matrix $S^\star_{11} - S_{12}^\star (S_{22}^\star)^{-1}S_{21}^\star$.

The finite convex problem~\eqref{eq:finiteconvex_revision} is a nonlinear semi-definite program (SDP) and to solve it tractably~\cite{ref:shafieezadeh2018wasserstein} proposed a tailored first-order method (Frank-Wolfe algorithm). The crucial observation in ~\cite[Definition 3.1]{ref:shafieezadeh2018wasserstein} is that the function $f(S)= \textrm{tr}\left[S_{11} - S_{12} S_{22}^{-1}S_{21}\right]$ satisfies the relation
\[
f(S) = \langle S, \nabla f(S)\rangle,
\]
because
\[
\langle S, \nabla f(S)\rangle = \left\langle \begin{pmatrix} S_{11} & S_{12} \\ S_{21} & S_{22}\end{pmatrix}, \begin{pmatrix} I_N & -S_{12}S_{22}^{-1} \\ -S_{22}^{-1}S_{21} & S_{22}^{-1} S_{21}S_{12} S_{22}^{-1}\end{pmatrix}\right\rangle = f(S). 
\]
Hence, the objective in problem~\eqref{eq:finiteconvex_revision} can be replaced with a linear approximation. In particular, following the iterative algorithm in~\cite[Equations 7(a) and 7(b)]{ref:shafieezadeh2018wasserstein} we consider
\begin{equation}\label{eq:fwalgo_revision}
S^{(k+1)} = a_k F(S^k) + (1-a_k) S^k, \forall k \geq 0,
\end{equation}
with $S^{(0)} = B^{-1}\tilde\Sigma B^{-\top} = A \Sigma_0 A^\top$ is the starting solution, $a_k$ represents a judiciously chosen step-size, and
\begin{equation}\label{eq:bisection_revision}
F(S) = \left\{
\begin{array}{ccl}
\arg&\max & \langle L, \nabla f(S)\rangle, \\
&\textrm{s.t.} & L  \succeq \bar\sigma I_{N+m}, ~\textrm{tr}\left[BLB^\top+\tilde\Sigma - 2\left(\tilde\Sigma^{1/2} BLB^\top\tilde\Sigma^{1/2}\right)^{1/2}\right]\leq\delta^2.
 \end{array}
\right.
\end{equation}
Problem~\eqref{eq:bisection_revision} can be further rewritten by considering the transformation $\tilde L \leftarrow BLB^\top$:
\begin{equation}\label{eq:bisection_revision_reduction}
F(S) = \left\{ 
\begin{array}{ccl}
\arg&\max &\langle \tilde L, B^{-\top}\nabla f(S)B^{-1}\rangle, \\
&\textrm{s.t.}& \tilde L  \succeq \tilde{\bar\sigma} I_{N+m},~\textrm{tr}\left[\tilde L+\tilde\Sigma - 2\left(\tilde\Sigma^{1/2} \tilde L\tilde\Sigma^{1/2}\right)^{1/2}\right]\leq\delta^2,
\end{array}
\right.
\end{equation}
which is exactly~\cite[Equations 7(b)]{ref:shafieezadeh2018wasserstein}. 

In conclusion, the finite-dimensional min-max problem~\eqref{eq:computeleft} can be solved using the iterations~\eqref{eq:fwalgo_revision} with the oracle~\eqref{eq:bisection_revision_reduction}, which is exactly~\cite[Equations 7(a) and 7(b)]{ref:shafieezadeh2018wasserstein}. Therefore, the Frank-Wolfe algorithm in~\cite{ref:shafieezadeh2018wasserstein} can be applied to numerically compute the distributionally robust predictor and the worst-case prior. We refer interested readers to~\cite[Algorithms 1 and 2]{ref:shafieezadeh2018wasserstein} for the concrete algorithm and the analysis thereof, which is beyond the scope of this work.

% \[
% F(S) = \begin{cases}
% \arg\max_{L  \succeq \bar\sigma I_{N+m}} \langle L, B^{-T}\nabla f(S)B^{-1}\rangle, \\
% \textrm{s.t. 
%  Tr}\left[BSB^\top+\tilde\Sigma - 2\left(\tilde\Sigma^{1/2} BSB^\top\tilde\Sigma^{1/2}\right)^{1/2}\right]\leq\rho^2.
% \end{cases}
% \]

% Thus the Nash equilibrium satisfies $G^\star = S^\star_{xy}(S^\star_{yy})^{-1}$ and $S^\star$ satisfies $S^\star = B^{-1}\tilde S^\star(B^\top)^{-1}$ and $\tilde S^\star$ solves
% \[
% \begin{array}{cl}
% \sup\limits_{\tilde S\succeq 0} & \left\langle (B^{-1})^\top \begin{pmatrix} I & -G^\star \\ -(G^\star)^\top & (G^\star)^\top G^\star\end{pmatrix} B^{-1}, \tilde S\right\rangle \\
% \st & \textrm{Tr}\left[\tilde S+\tilde\Sigma - 2\left(\tilde\Sigma^{1/2} \tilde S\tilde\Sigma^{1/2}\right)^{1/2}\right]\leq\rho^2.
% \end{array}
% \]
% By Viet's paper, Lemma A.2, we have such an $\tilde S^\star$ satisfies the equation
% \[
% \tilde S^\star = ( I - \gamma^{-1} D)^{-1} \tilde\Sigma (I - \gamma^{-1} D)^{-1}
% \]
% where $\gamma$ is some positive number and 
% \[
% D =(B^{-1})^\top \begin{pmatrix} I & -G^\star \\ -(G^\star)^\top & (G^\star)^\top G^\star\end{pmatrix}B^{-1}
% \]

}

\review{

%% This is the right place, to print the list of used acronyms.
\renewcommand{\nomname}{B List of Notations}

\printnomenclature[2cm]

}

\end{appendix}
%%%%%%%%%%%%%%%%%%%%%%%%%%%%%%%%%%%%%%%%%%%%%%
%% Multiple Appendixes:                     %%
%%%%%%%%%%%%%%%%%%%%%%%%%%%%%%%%%%%%%%%%%%%%%%
%\begin{appendix}
%\section{???}
%
%\section{???}
%
%\end{appendix}

%%%%%%%%%%%%%%%%%%%%%%%%%%%%%%%%%%%%%%%%%%%%%%
%% Support information, if any,             %%
%% should be provided in the                %%
%% Acknowledgements section.                %%
%%%%%%%%%%%%%%%%%%%%%%%%%%%%%%%%%%%%%%%%%%%%%%
\textbf{Acknowledgements.} 
Material in this paper is based upon work supported by the Air Force Office of Scientific Research under award number FA9550-20-1-0397. J.~Blanchet acknowledges additional support from NSF grants 1915967, 2118199, 2229012, 2312204, and ONR 13983111.
\newpage
\bibliographystyle{siam}
\bibliography{bibliography}

\end{document}